\documentclass[10pt]{amsart}
\usepackage{amsmath,bbm}
\usepackage{amsfonts}
\usepackage{amsthm}
\usepackage{amssymb}
\usepackage{url}
\usepackage{verbatim}
\usepackage{subcaption}
\usepackage[usenames,dvips]{color}
\usepackage{graphicx}
\usepackage[dvipsnames]{xcolor}
\usepackage{mathrsfs}
\usepackage{lmodern}
\usepackage{bm}
\usepackage{math}
\usepackage{mathtools}
\usepackage{mathrsfs}
\usepackage{geometry}
\numberwithin{equation}{section}
\numberwithin{subsection}{section}
\newtheorem{theorem}{Theorem}
\newtheorem{lemma}{Lemma}[section]

\newtheorem{definition}[lemma]{Definition}
\newtheorem{remark}[lemma]{Remark}

\newtheorem{proposition}[lemma]{Proposition}

\newcommand{\ubmmu}{\underline{\bmmu}}

\newcommand{\bmkappa}{\bm{\kappa}}

\newcommand{\bmmu}{\bm{\mu}}

\newcommand{\bmxi}{\bm{\xi}}

\newcommand{\bmc}{\bm{c}}

\newcommand{\bmD}{\bm{D}}

\newcommand{\bmf}{\bm{f}}

\newcommand{\bmj}{\bm{j}}
\newcommand{\bmk}{\bm{k}}

\newcommand{\bmu}{\bm{u}}

\newcommand{\bmW}{\bm{W}}
\newcommand{\bmx}{\bm{x}}

\newcommand{\bmz}{\bm{z}}

\newcommand{\bmzero}{\bm{0}}

\newcommand{\sL}{{\mathscr{L}}}
\newcommand{\fI}{\mathfrak{I}}

\renewcommand{\fL}{\mathfrak{L}}

\DeclareMathOperator{\gradp}{\nabla^{\perp}}

\newcommand{\ubmk}{{\underline{\bmk}}}

\newcommand{\uw}{{\underline{w}}}

\newcommand{\Res}[2]{\mathtt{R}\begingroup 
\setlength\arraycolsep{0pt}{\scriptsize \begin{matrix} #2 \\[0.7mm] #1 \end{matrix} }\endgroup}

\newcommand{\white}[1]{{\textcolor{white}{#1}}}
\newcommand{\ent}[6]{\begingroup 
\setlength\arraycolsep{-2pt}\begin{matrix}{\tB_{#1}^{[#2]}} &\white{{|}^{|}}_{#3,#4}^{#5,#6}\end{matrix}\endgroup}

\newcommand{\entL}[6]{\begingroup 
\setlength\arraycolsep{-2pt}\begin{matrix}{\tL_{#1}^{[#2]}} &\white{{|}^{|}}_{#3,#4}^{#5,#6}\end{matrix}\endgroup}

\usepackage[colorlinks=true]{hyperref}

\title[]{ Instabilities of internal gravity waves \\in the two-dimensional  Boussinesq system }

\author{R. Bianchini$^{(\dagger)}$}
\address[$\dagger$]{IAC, Consiglio Nazionale delle Ricerche (CNR), via dei Taurini 19, 00185, Rome, Italy}
\email[$\dagger$]{roberta.bianchini@cnr.it}

\author{A. Maspero$^{(\star)}$}
\address[$\star$]{International School for Advanced Studied (SISSA), via Bonomea 265, 34136, Trieste, Italy}
\email[$\star$]{amaspero@sissa.it}

\author{S. Pasquali$^{(\ast)}$}
\address[$\ast$]{International School for Advanced Studied (SISSA), via Bonomea 265, 34136, Trieste, Italy}
\email[$\ast$]{stefano.pasquali@sissa.it}

\begin{document}

\begin{abstract}
We consider a two-dimensional, incompressible, inviscid fluid with variable density, subject to the action of gravity. Assuming a stable equilibrium density profile, we adopt the so-called Boussinesq approximation, which neglects density variations in all terms except those involving gravity. This model is widely used in the physical literature to describe internal gravity waves.

In this work, we prove a modulational instability result for this system. Specifically, we show that the Floquet operator associated with the linearization around a small-amplitude traveling wave admits at least one eigenvalue with positive real part, bifurcating from a double eigenvalue of the unperturbed linear operator. This can be regarded as the first rigorous justification of the Parametric Subharmonic Instability (PSI) of inviscid internal waves, wherein energy is transferred from an initially excited primary wave to two secondary waves with different frequencies.
Our approach uses Floquet–Bloch decomposition and Kato's similarity transformations to compute rigorously the perturbed eigenvalues without requiring boundedness of the perturbed operator - differing fundamentally from prior analyses involving viscosity. 

Notably, the inviscid setting is especially relevant in oceanographic applications, where viscous effects are often negligible.
\end{abstract}

\maketitle

\tableofcontents

\section{Introduction}
We consider two-dimensional, incompressible fluids \emph{stratified} under the Boussinesq approximation, where density variation $\rho=\rho(t, x, y)$ is a small fluctuation of its (constant) average $\rho_0$, and does not affect inertial terms; see, for example, \cite[Section 2]{dauxois2018instabilities} for an introduction and \cite{long1965boussinesq, benjamin1966internal} for a discussion on the Boussinesq approximation. 
In the simplest setting, density at the hydrostatic equilibrium is given by an affine, stable profile $\bar{\rho}(y)$ with $\bar{\rho}'(y)=\beta<0$, and density variation takes the form $\rho= \bar \rho (y) + \frac{\rho_0}{\mathfrak{g}} b$, where $\mathfrak{g}$ is gravity acceleration constant. 
The two-dimensional \emph{inviscid} Boussinesq equations for (stably) stratified fluids in $\R^2$ are given below (see (1.5) in \cite{bianchini2024reflection}):
\begin{align}\label{eq:bouss}
\partial_t b - N^2 v + (\mathbf u \cdot \nabla) b &=0, \notag \\
\partial_t u + \partial_x P +  (\mathbf u \cdot \nabla) u &=0, \notag \\
\partial_t v + b + \partial_y P +  (\mathbf u \cdot \nabla) v &=0, \notag \\
\partial_x u + \partial_y v&=0,
\end{align}
where $\mathbf u=(u,v)=(u(t, \bmx), v(t,\bmx)): \R^+ \times \R^2 \to \R^2$, $\bmx=(x, y)$, is the (divergence-free) velocity field,  $P=P(t, \bmx): \R^+ \times \R^2 \to \R$ is the  pressure, \begin{align} \label{eq:relNrho}
N^2 &:= - \frac{\mathfrak{g}\bar{\rho}'(y)}{\rho_0}
\end{align}
is the so-called Brunt--V\"ais\"al\"a frequency, which represents the maximum frequency of oscillations of internal gravity waves supported by the above system. Introducing the vorticity
\begin{align} \label{eq:vort}
\omega &:= \partial_x v - \partial_y u \, = \, - \Delta \psi,
\end{align}
system \eqref{eq:bouss} can be written in  vorticity-stream formulation (see Sec. 2 of \cite{abarbanel1986nonlinear}; see also \cite{benjamin1986boussinesq}):
\begin{equation} \label{eq:BIISysNew}
\begin{cases}
\partial_t b =  - \, N^2 \, \partial_x \psi -  \{b,\psi\}, \\
\partial_t \omega = - \partial_x b -  \{\omega,\psi\},
\end{cases} \quad -\Delta\psi = \omega ,
\end{equation} 
where the velocity can be recovered from the Biot-Savart law
\begin{equation}
    \label{eq:BS}
    \mathbf{u}=\nabla^\perp(-\Delta)^{-1}\omega,
\end{equation}

and for any couple of functions $f, g$ we denote the Poisson brackets $\{f, g\}=\partial_x f \partial_y g - \partial_y f \partial_x g$.

\subsection*{State of the art on the Boussinesq equations}
From a mathematical perspective, the analysis of the Boussinesq system \eqref{eq:bouss} has recently attracted considerable attention. Among the outstanding open problems is the question of whether solutions originating from smooth, finite-energy initial data can develop finite-time singularities in the case \(N=0\); see the recent blow-up result for \(C^{1,\alpha}\) velocity fields in the unstable setting \cite{pasqualotto2023}. Even when the linearized system at the origin is spectrally stable — as is the case for \(N^2 \neq 0\) \cite{desjardins2021normal} — the issue of global regularity versus finite-time blow-up for smooth solutions remains unresolved; see \cite{wid2024} for recent progress in this direction.

Beyond the question of global existence versus blow-up, another important aspect of the Boussinesq system concerns the possible growth of Sobolev norms over time, which is closely connected to the formation of small scales and turbulent behaviors; see \cite{kiselev2025}. In the inviscid setting considered here, algebraic-in-time growth of the density gradient \(\|\nabla b\|_{L^2}\) or the vorticity \(\|\omega\|_{L^2}\), with \(\omega = \nabla \cdot \mathbf{u}^{\perp}\), has been demonstrated under suitable symmetry assumptions on the initial data \cite{kiselev2025}, and also in the presence of a background shear flow \cite{BBCZD2021}.

\subsection*{Internal waves and the scope of this work}
The Boussinesq equations for stratified fluids describe the propagation of internal gravity waves, which play a crucial role in stratified fluids, contributing significantly to ocean mixing \cite{dauxois2018instabilities}. The propagation and reflection of internal waves are by now fairly well understood, with several mathematical results available; see, for instance, \cite{de2020attractors, bianchini2021near, desjardins2021normal}.
Nonlinear phenomena, such as  energy transfers, dissipation mechanisms, and related turbulent behaviors remain active areas of investigation; see the recent review \cite{dauxois2023}.

This paper provides the first rigorous mathematical study of instabilities in internal gravity waves in the \emph{inviscid} setting, focusing on growth mechanisms generated by resonant interactions under small-amplitude perturbations. A preliminary mathematical study in the \emph{viscous} case was provided in \cite{bianchini2023triadic}. 

Here, we investigate the \emph{Parametric Subharmonic Instability} (PSI) of inviscid internal waves - the inviscid counterpart of the \emph{Triadic Resonant Instability} (TRI); see~\cite[Section~3.5]{dauxois2018instabilities}. Remarkably, the \emph{inviscid} setting is particularly relevant for oceanographic applications, where viscosity is often negligible \cite{dauxois2018instabilities}.
In the 1970s, internal gravity waves were discovered to be unstable to infinitesimal perturbations, which can grow to form temporal and spatial resonant triads \cite{davis1967stability}. 
The convection term in the Boussinesq system plays a crucial role by enabling energy transfer among waves of different frequencies. Specifically, under the effect of nonlinearities, a primary wave perturbed by small disturbances can transfer part of its energy to two secondary waves of lower frequency (subharmonics). Thus, TRI and PSI represent direct energy transfer mechanisms, without the need for turbulent cascade.
TRI or PSI arises when three waves of wavevectors  $\mathbf{k}_j$ and time frequency $\Omega (\mathbf{k}_j)$, $j \in \{1, 2, 3\}$ interact resonantly, satisfying the relations \( \sum_{j=1}^3 \Omega (\mathbf{k}_j) = 0 \) for the frequencies and \( \sum_{j=1}^3 \mathbf{k}_j = \bmzero \) for the wavevectors  \cite{dauxois2018instabilities}. This instability has been observed in laboratory experiments \cite{benielli1998excitation, joubaud2012experimental} and confirmed by oceanic field measurements \cite{hibiya2002nonlinear, mackinnon2013parametric}. Interestingly, recent physics studies \cite{akylas2023stability} have used Floquet methods to study PSI and to identify a novel ``broadband'' instability at large Floquet parameters.

Although TRI has previously been studied in the presence of viscosity, our focus here is its inviscid counterpart, PSI~\cite[Section~3.1]{dauxois2018instabilities}. The earlier work \cite{bianchini2023triadic} used classical perturbation theory to construct approximate eigenvalues (quasi-modes) in the \emph{viscous} case. However, that approach requires the perturbation to be relatively compact or bounded - a condition violated in the inviscid case, where the unperturbed operator is of zeroth order and the perturbation (convection term) is of first order. To overcome this, we develop a new framework for the inviscid setting, based on Floquet-Bloch decomposition and Kato’s similarity transformations, which allows to rigorously compute perturbed eigenvalues without assuming boundedness of the perturbation. This approach fundamentally differs from prior work involving viscosity~\cite{bianchini2023triadic}, and provides the first rigorous justification of PSI for inviscid internal waves.

 We denote by $\dot{H}^1(\mathbb{R}^2)
:= \left\{ f \in \mathcal{S}'(\mathbb{R}^2) : \nabla f \in L^2(\mathbb{R}^2) \right\}$ and $\dot H^{-1}(\mathbb R^2):=(\dot{H}^1(\mathbb{R}^2))'$ the corresponding dual space. 

Our main result reads as follows.
\begin{theorem}
Let $\ubmk = (\tm,\tn)^T \in \mathbb{Z}^2$ with $\tm,\tn > 0$, and let $\epsilon>0$. Consider the plane wave solution to \eqref{eq:BIISysNew} given by
\begin{align*}
\begin{pmatrix}
b_\e(t,\bmx) \\ \omega_\e(t,\bmx)
\end{pmatrix}
= \e \begin{pmatrix}
   1\\
   |\ubmk|/N
\end{pmatrix} \sin (\ubmk \cdot \bmx - \Omega(\ubmk)t), \quad 
\quad \Omega(\ubmk) = N \frac{\tm}{|\ubmk|}.
\end{align*}
Let $\mathcal{L}_\epsilon$ denote the operator obtained by linearizing \eqref{eq:BIISysNew} around this solution. Then there exists $\epsilon_0 > 0$ such that, for all $0 < \epsilon \le \epsilon_0$, the operator $\mathcal{L}_\epsilon$, viewed as an (unbounded) operator on $X:={H}^1(\mathbb{R}^2) \times (L^2\cap \dot H^{-1})(\mathbb{R}^2)$, 
has unstable spectrum, i.e.
\[
\sigma_{X}(\mathcal{L}_\epsilon)
\cap \{ \lambda \in \mathbb{C} : \Re \lambda > 0 \}
\neq \emptyset.
\]
\end{theorem}

\noindent The choice of the space $X$ is motivated by the fact that the  linear evolution 
$\partial_t 
\big(
\begin{smallmatrix}
b \\ \omega
\end{smallmatrix}\big)
= \mathcal{L}_\epsilon 
\big(
\begin{smallmatrix}
b \\ \omega
\end{smallmatrix}\big)
$ 
is well-posed on it.

Historically, wave-wave interactions have been studied via formal weakly nonlinear expansions, which allow for the formal computation of the real part of unstable eigenvalues~\cite{dauxois2018instabilities}.

Here, we perform a rigorous analysis of the spectrum of the inviscid Boussinesq equations linearized around a background (primary) wave, and we rigorously prove the existence of eigenvalues with strictly positive real part. This confirms the expressions previously obtained in the literature, as discussed in Section~\ref{sec:comparison}. Our approach is based on modulational instability, a ubiquitous phenomenon whereby traveling wave solutions of nonlinear dispersive equations become unstable under long-wave perturbations.

In fluid dynamics, modulational instability traces back to the pioneering works of Benjamin and Feir, Zakharov, Lighthill, and McLean~\cite{benjamin1967instability, benjamin1967disintegration, zakharov1967instability, lighthill1965contributions, mclean1982instabilities} on the instability of Stokes waves in irrotational water waves.

This problem has recently attracted renewed interest due to the development of new analytical tools that allow for the rigorous computation of portions of the unstable spectrum of the linearized water wave operator. These advancements have been made both near the origin~\cite{bridges1995proof, nguyen2023proof, berti2022full, berti2023benjamin, berti2024stokes, hur2023unstable}, away from the origin~\cite{Creedon_Deconinck_Trichtchenko_2022, berti2025first, berti2024infinitely, hur2023unstable}, and even in the context of 3D water waves~\cite{creedon2025proof, creedon2026transverse, jiao2025small}. We further emphasize that modulational instability can also be interpreted in terms of Krein signatures of eigenvalues; see \cite{kapitula2004counting,kapitula2013spectral,noble2023spectral}.

In this paper, we adapt the ideas in~\cite{berti2022full, berti2024infinitely}, with some important differences. First, the Boussinesq equations in vorticity form admit a non-canonical Hamiltonian structure~\cite{abarbanel1986nonlinear}; as a result, the linearized Boussinesq operator around an internal wave (see~\eqref{eq:Leps}) does not appear to be linearly Hamiltonian. Therefore, we use only spectral projectors to compute the matrix representing the action of the linearized operator on a two-dimensional invariant subspace, rather than employing Kato's transformation operator (see Lemma~\ref{lem:MatLmueps}). However, the linearized operator turns out to be reversible, and we exploit this structure.

Another important difference with respect to~\cite{berti2022full, berti2024infinitely} concerns the spectral structure of the unperturbed linearized operator $\mathcal{L}_{\bmmu,0}$ defined in~\eqref{eq:Lmu0}. Specifically, we take the Floquet parameter $\bmmu \in \mathbb{R}^2 \setminus \{ \bmzero \}$ to enforce the resonance condition
\begin{align}\label{res.intro}
    \Omega(\ubmk) - \Omega(\ubmk + \bmmu) = \Omega(\bmmu),
\end{align}
which is inspired by the TRI and PSI phenomena, where energy is transferred from a primary wave to two secondary waves. In Lemma~\ref{lem:Rk}, we analytically characterize the set $\mathcal{R}_{\ubmk}$ of values of $\bmmu$ that satisfy~\eqref{res.intro}.

The key observation is that whenever $\bm{\mu} \in \mathcal{R}_{\ubmk} \cap ( \mathbb{R}^2 \setminus \mathbb{Z}^2 )$
the operator $\mathcal{L}_{\bm{\mu}, 0}$ possesses, as an unbounded operator on
$H^1(\mathbb{T}^2) \times L^2(\mathbb{T}^2)$, a purely imaginary eigenvalue $\lambda^+$
with multiplicity \emph{at least} $2$. However, due to the higher dimensionality of $\mathbb{T}^2$,
this multiplicity may be greater than two, complicating the bifurcation analysis.
To address this issue, we exploit the fact that the Boussinesq traveling wave solution only involves
Fourier coefficients at $\pm \ubmk$. Therefore, the operator $\mathcal{L}_{\bm{\mu},0}$ leaves
invariant the subspace $X^1_{\ubmk} \subseteq H^2(\mathbb{T}^2) \times H^1(\mathbb{T}^2)$,
consisting of functions supported only on wavevectors proportional to $\ubmk$ (see~\eqref{Hlk}).
We then prove that, when restricted to $X^1_{\ubmk}$, and for any
$\bm{\mu} \in \mathcal{R}_{\ubmk} \cap ( \mathbb{R}^2 \setminus \mathbb{Z}^2 )$
(modulo a discrete set), the eigenvalue $\lambda^+$ is isolated and has algebraic multiplicity
exactly $2$ (see~\eqref{eq:DecompSpecLmu0}). This reduction enables us to apply the perturbative theory of isolated eigenvalues (see Section~\ref{sec:Kato}), ultimately reducing the problem to the computation of the eigenvalues of a $2 \times 2$ matrix $\mathtt{L}(\bmmu,\epsilon)$ (see~\eqref{eq:MatRepr}). In Theorem~\ref{thm:expansion}, we show that these eigenvalues are given by
\begin{equation*}
    \lambda^\pm(\bmmu, \epsilon) = \lambda^+ + \im \, \mathcal{O}(\epsilon^2) \pm \epsilon \sqrt{\mathtt{e}(\bmmu) + \mathcal{O}(\epsilon)},
\end{equation*}
and we provide an explicit expression for the real-valued function $\mathtt{e}(\bmmu)$, together with its asymptotics for $|\bmmu| \ll 1$ and for $|\bmmu| \gg 1$, showing that it is strictly positive in both regimes.

In Section~\ref{sec:Taylor}, we compute the small-$\epsilon$ expansions of the entries of the matrix $\mathtt{L}(\bmmu,\epsilon)$ using so-called entanglement coefficients, based on the jets of the operator  $\sL_{\bmmu,\epsilon} := \mathcal{L}_{\bmmu,\epsilon}|_{X^1_{\ubmk}}$. By generalizing the approach of~\cite{berti2024infinitely}, we define, for each $\bmkappa \in \ubmk \mathbb{Z}$ and for any operator $A$ acting on  $X^1_{\ubmk}$ , the projected operators $A^{[\bmkappa]}$ whose matrix coefficients are supported on the bands where $(j_2 - j_1) \ubmk = \bmkappa$, for $j_1, j_2 \in \mathbb{Z}$ (see~\eqref{band.def}). This framework allows us to exploit structural properties of the jets of $\sL_{\bmmu,\epsilon}$ and leads to a more efficient computation of the entanglement coefficients.

\subsection*{Notation and conventions}
We list the notation and conventions used throughout the paper:
for a given vector $\bmf=(f_1,f_2)^T$, we use the notation $\bmf^{\perp} \coloneqq (f_2,-f_1)^T$; similarly, $\gradp \coloneqq (\partial_y,-\partial_x)^T$; we denote by \( \mathcal{O} (\e^{n}) \) an analytic function that,  for some \( C > 0 \) and for any small \( \e \), satisfy
$
|\mathcal{O} (\e^n) |  \leq C |\e|^n$; 
we denote by $r(\e^n)$ a function $\mathcal{O}(\e^n)$ which is real-valued;
for any $ k\in\mathbb{N} $,  $\cO_k$ is an operator from  ${H}^2(\mathbb{T}^2) \times {H}^1(\mathbb{T}^2)$  to  ${H}^1(\mathbb{T}^2) \times L^2(\mathbb{T}^2)$, of size $\e^k $; for any couple of functions $f, g$, the Poisson brackets $\{f, g\}=\partial_x f \partial_y g - \partial_y f \partial_x g$. \\

\subsection*{Structure of the paper} 
In Sec. \ref{sec:LinSys} we prove the existence of internal gravity waves, we write the linearized system around an internal gravity wave, and we study the spectrum of the operator $\cL_{\e}$ appearing in the linearized system. Next, in Sec. \ref{sec:Kato} we construct rigorously the perturbed spectral subspaces $\cV_{\bmmu,\e}$, we compute the matrix representing the action of the operator $\cL_{\bmmu,\e}|_{X^1_{\ubmk}}$ on a suitable basis $\cF$ of $\cV_{\bmmu,\e}$; moreover, we provide an expansion in $\e$ of the entries of such matrix. Then, in Sec. \ref{sec:Taylor} we introduce the so-called entanglement coefficients, which are used to prove the expansion of Proposition \ref{prop.exp}. In Sec. \ref{sec:proof.prop} we prove Proposition \ref{prop.exp}. In Sec. \ref{sec:comparison} we compare our main result with the instability results found in physical literature. Finally, in Appendix \ref{app:Rk} we prove Lemma \ref{lem:Rk}, and in Appendix \ref{app:de} we prove Proposition \ref{prop:wBound}.

\section{Internal plane waves and the linearized system} \label{sec:LinSys}

System \eqref{eq:BIISysNew} supports the propagation of traveling wave solutions - called \emph{internal gravity waves} - of the form 
\begin{align} \label{eq:travel}
b(t,\bmx) = \check{b}( \bmx - \bmc t) , &\;\; \omega(t,\bmx) = \check{\omega}( \bmx- \bmc t) , \quad \bmc = (c_1, c_2)^\top \in \R^2 \ 
\end{align}
where the profiles $\check{b}(\bmx)$ and $\check{\omega}(\bmx)$ are
plane waves. The
following result can be easily verified. 
\begin{proposition}\label{prop:trav}
For any  $\epsilon>0$ and $\underline{\bmk} = (\tm, \tn)^T \in \Z^2 \setminus \{ \bmzero \}$, 
\begin{align}\label{eq:primary_wave}
b_{\epsilon}(\bmx) &:= \epsilon  \, \sin(\underline{\bmk} \cdot \bmx) , \qquad 
\omega_{\epsilon}(\bmx) := \epsilon \, \frac{|\ubmk|}{N} \, \sin(\underline{\bmk} \cdot \bmx) ,  \; \; \\
\underline{\bmc} & := \frac{N \,\tm }{|\underline{\bmk}|^3} \, \underline{\bmk} , 
\qquad 
\psi_{\epsilon}(\bmx) := \epsilon \, \frac{1}{N \, |\ubmk|} \, \sin(\underline{\bmk} \cdot \bmx) 
\;\;  \label{eq:trav}
\end{align}
are solutions to the following system:
\begin{equation} \label{eq:SteadyTranslBIISysNew}
\begin{cases}
(\underline{\bmc} \cdot \nabla) b  - N^2 \, \partial_x \psi -  \{b,\psi\}= 0, \\
(\underline{\bmc} \cdot \nabla) \omega - \partial_x b -  \{\omega,\psi\}= 0.
\end{cases}
\end{equation} 
\end{proposition}

\begin{remark}
Notice that \eqref{eq:travel}-\eqref{eq:trav} are internal (plane) wave solutions to \eqref{eq:BIISysNew} of the form 
\begin{align} \label{eq:travel2}
\begin{pmatrix}
b \\ \bmu
\end{pmatrix}
&= \frac{1}{2} \left[ \bmW \, e^{\mathrm{i} \left( \ubmk \cdot \bmx - \frac{N \tm}{|\ubmk|}  t \right)} + \overline{\bmW} \, e^{\mathrm{i} \left( \frac{N \tm}{|\ubmk|} t - \ubmk \cdot \bmx \right)} \right] , \;\; \bmW \coloneqq 
\begin{pmatrix}
\mathrm{i} \frac{|\ubmk| N}{\tn} u_0 \\ u_0 \\ - \frac{\tn}{\tm} u_0 
\end{pmatrix}
, \;\; u_0 = \epsilon \, \frac{\tn}{N \, |\ubmk|} ,
\end{align}
see also \cite[(2.4)]{bianchini2023triadic}.
\end{remark}

\subsection{The linearized operator }
We now linearize \eqref{eq:BIISysNew}, in a moving frame, near the traveling wave $\left( b_{\epsilon},\omega_{\epsilon}  \right)$ of wavevector $\ubmk$ in \eqref{eq:trav}. 
The linearized system is 
\begin{align}\label{eq:linearized}
\partial_t \,
\begin{pmatrix}
\tilde{b} \\ \tilde{\omega}
\end{pmatrix}
&=  \mathcal{L}_{\epsilon}
\begin{pmatrix}
\tilde{b} \\ \tilde{\omega}
\end{pmatrix}.
\end{align}
where $\mathcal{L}_{\epsilon}$ is the real operator
\begin{equation}\label{eq:Leps}
\mathcal{L}_{\epsilon}:=   
(\underline{\bmc} \cdot \nabla) \mathbb{I}_2 +
\begin{pmatrix}
- \gradp \psi_{\epsilon} \cdot \nabla & \left( -N^2 \partial_x - \nabla b_{\epsilon} \cdot \gradp \right) (-\Delta)^{-1} \\
-\partial_x  & - \gradp \psi_{\epsilon} \cdot \nabla  - \nabla \omega_{\epsilon} \cdot \gradp  (-\Delta)^{-1} 
\end{pmatrix}.
\end{equation}

For concreteness, we consider $\cL_\e$ as an operator on 
\begin{equation}\label{X.def}
X:=  H^1(\R^2) \times (L^2\cap \dot H^{-1})(\R^2)
\quad \mbox{ with dense domain } D(\cL_\e):= \{ U \in X \colon \cL_\e U \in X \} \ .
\end{equation}

To justify the choice of the space $X$, we prove the following energy estimate for \eqref{eq:linearized}:
\begin{equation}\label{energy_est}
\|b(t)\|_{H^1} + \|\omega(t)\|_{L^2 \cap \dot H^{-1}}
\lesssim 
\|b(0)\|_{H^1} + \|\omega(0)\|_{L^2 \cap \dot H^{-1}}
\exp\Big(Ct (\| b_\e\|_{W^{2,\infty}}+\| \mathbf u_\e\|_{W^{2, \infty}})\Big),
\end{equation}
where $(b_\e(\bx), \omega_\e(\bx))$ is the primary wave \eqref{eq:primary_wave}, 
$\mathbf{u}_\e :=\nabla^\perp \psi_\e$ in \eqref{eq:primary_wave}
and $C>0$ is a universal constant.
To prove such energy estimate is more convenient to write the linearized system \eqref{eq:linearized} in velocity form, namely, 
recalling from the Biot-Savart law \eqref{eq:BS} that $\mathbf{u}=\nabla^\perp(-\Delta)^{-1}\omega$, we consider the linear system 
\begin{align}\label{eq:bouss-lin-vel}
\partial_t b - N^2 v + (\mathbf u_\e \cdot \nabla) b + (\mathbf u \cdot \nabla) b_\e &=0, \notag \\
\partial_t u + \partial_x P +  (\mathbf u_\e \cdot \nabla)u + (\mathbf u \cdot \nabla) u_\e &=0, \notag \\
\partial_t v + b + \partial_y P +  (\mathbf u_\e \cdot \nabla) v + (\mathbf u \cdot \nabla) v_\e &=0, \notag \\
\partial_x u + \partial_y v&=0.
\end{align}
By incompressibility we have the zeroth-order estimate:
\begin{align}
\frac{d}{dt}\mathcal{E}_0(t)
:= \frac{d}{dt}(\|b\|_{L^2}^2+N^2\|\mathbf{u}\|_{L^2}^2)
&\lesssim  \|\mathbf{u}\|_{L^2}(\|\nabla b_\e\|_{L^\infty} \|b\|_{L^2}+\|\nabla \mathbf{u}_\e\|_{L^\infty}\|\mathbf{u}\|_{L^2})\notag\\
&\lesssim (\|b_\e\|_{W^{1, \infty}}+\|\mathbf{u}_\e\|_{W^{1, \infty}}) \mathcal{E}_0(t) \ . 
\end{align}
Recalling  $\mathbf{u}=\nabla^\perp(-\Delta)^{-1}\omega$, one has control over 
$\|\mathbf{u}\|_{L^2} \simeq \|\omega\|_{\dot H^{-1}}$.
Then, denoting $\textbf{D}:=\frac{1}{\im} \nabla$, we compute
\begin{align}
\frac{d}{dt}\mathcal{E}_1(t)
:= \frac{d}{dt}(\||\mathbf{D}|b\|_{L^2}^2+N^2\|\omega\|_{L^2}^2)
&\lesssim \|\mathbf u_\e\|_{W^{1, \infty}}\|b\|_{\dot H^1}^2 \notag\\
&+ (\|b_\e\|_{W^{2, \infty}}+\|\mathbf u_\e\|_{W^{2, \infty}})
(\|b\|_{\dot H^1}+\|\omega\|_{L^2\cap {\dot H^{-1}}})\|\omega\|_{L^2\cap {\dot H^{-1}}} \notag\\ 
& \lesssim (\|b_\e\|_{W^{2, \infty}}+\|\mathbf u_\e\|_{W^{2, \infty}})(\mathcal{E}_0(t)+\mathcal{E}_1(t)).
\end{align}
In the above estimate, we used the Kato-Ponce commutator bound
\[
\|[|\mathbf{D}|, f]g\|_{L^2} \lesssim \|\nabla f\|_{L^\infty}\|g\|_{L^2}.
\]
In view of \eqref{energy_est}, the operator $\cL_\e$ generates a strongly continuous semigroup of bounded operators on $X$ and is therefore closed in its domain.

Since the operator $\cL_\e$ depends on functions bi-periodic in the $\bmx$-variable, it is natural to study its spectrum by means of the 
 Bloch-Floquet transform, that we now introduce: for any $v \in \cS(\R^2)$, we define 
\[
Bv(\bmmu,\bmx)
:= \sum_{\bmk\in\mathbb{Z}^2}
   e^{\,\mathrm{i}\,\bmk\cdot\bmx}\,
   \widehat v(\bmk+\bmmu),
\qquad \bmmu\in\R^2,
\]
where $\mathfrak{F}(v)(\bmxi) \coloneqq \hat v(\bmxi) \coloneqq \frac{1}{(2\pi)^2} \int_{\R^2} e^{- \im \bmxi \cdot \bmx} v(\bmx) \, \di \bmx$ is the Fourier transform of $v$. For a vector function the Bloch-Floquet transform is defined component-wise.

Writing every $\bmmu\in\R^2$  uniquely as 
$\bmmu=\tilde\bmmu+\bm n$ with 
$\tilde\bmmu\in[-\frac 12, \frac 12)^2$ and $\bm n\in\Z^2$, the Bloch-Floquet  transform satisfies the Floquet periodicity relation
\begin{equation}\label{eq:periodicity}
(Bv)(\tilde\bmmu+\bm n,\bmx)
= e^{-\,\mathrm{i}\,\bm n\cdot\bmx}\,(Bv)(\tilde\bmmu,\bmx),
\qquad \bm n\in\Z^2 .
\end{equation}

On the fundamental Brillouin zone $\tilde\bmmu\in[-\frac 12, \frac 12)^2$, 
the inversion formula allows to recover $v(\bmx)$ from its Bloch-transform $(Bv)(\tilde\bmmu,\bmx)$ by 
\begin{equation}\label{eq:inverse}
v(\bmx)
=
\int_{ [-\frac 12, \frac 12)^2} \quad
e^{\,\mathrm{i}\,\tilde\bmmu\cdot\bmx}\,
(Bv)(\tilde\bmmu,\bmx)\,
\mathrm{d}\tilde\bmmu ,
\end{equation}
and a direct computation shows the isometry property
\begin{equation}\label{isom}
\|v\|_{L^2(\R^2)}
=
\|B v \|_{L^2( [ -\tfrac12,\tfrac12 )^2;\,L^2(\T^2))},
\end{equation}
which  allows to  extend 
$B\colon 
L^2(\R^2)
\to
L^2([-\frac 12, \frac 12)^2;\,L^2(\T^2))$ as a continuous bijection.
In view of the identity  $B[\grad g](\bmmu, \bmx) = (\grad + \im \bmmu ) (Bg)(\bmmu, \bmx)$ and the isometry 
\eqref{isom}, one has  
\begin{equation}\label{iso.dotH}
\norm{g}_{ \dot H^1(\R^2)} = 
\norm{({\bf D} + \bmmu )(Bg)}_{L^2([-\frac 12, \frac 12)^2;\,  L^2(\T^2))}^2 
\end{equation}
  and, by duality,  
\begin{equation}\label{iso.dotH-1}
\norm{g}_{ \dot H^{-1}(\R^2)} = 
\norm{|{\bf D} + \bmmu|^{-1}(Bg)}_{L^2([-\frac 12, \frac 12)^2;\,  L^2(\T^2))}^2. 
\end{equation}

If $U \in X$ (cf. \eqref{X.def}) is such that its Bloch-Fourier transform    $B U$  is supported  in the  $\bmmu$-variable  away from any vector of integers $\mathbb Z^2$, then \eqref{iso.dotH} and \eqref{isom} yield
\begin{equation}
    \label{B.X}
    \norm{U}_{X} \simeq \norm{B U}_{L^2([-\frac 12, \frac 12)^2;\,  H^1(\T^2)\times L^2(\T^2))}
\end{equation}

 The advantage of the transform $B$ is to conjugate
the operator $\cL_\e$ to the fiber operator $\cL_{\bmmu, \e}$, namely
\begin{equation} \label{eq:BlochFloquet}
B[\cL_\e U ](\bmmu, \bmx) = \cL_{\bmmu, \e}[B U(\bmmu, \cdot)](\bmx) \ , \quad \forall U \in X\ , 
\end{equation}

where $\cL_{\bmmu, \e}$ is defined by
\begin{equation}
\begin{aligned}
\label{eq:Lmueps}
\mathcal{L}_{\bmmu, \epsilon}  &\coloneqq
\underline{\bmc} \cdot (\nabla + \im \bmmu) \mathbb{I}_2 +
\begin{pmatrix}
- \gradp \psi_{\epsilon} \cdot (\nabla + \im \bmmu)  & \left( -N^2(\partial_x + \im \mu_1)- \nabla b_{\epsilon} \cdot (\gradp  + \im \bmmu^\perp)\right) \, \vert {\bf D} + \bmmu \vert^{-2} \\
-\partial_x - \im \mu_1 & 
 - \gradp \psi_{\epsilon} \cdot (\nabla + \im \bmmu)
 - \nabla \omega_{\epsilon} \cdot (\gradp  + \im \bmmu^\perp)  \vert {\bf D} + \bmmu \vert^{-2}
\end{pmatrix} \ .
\end{aligned}
\end{equation}
For any $\bmmu=\tilde\bmmu+\bm n\in\R^2$, relation \eqref{eq:periodicity}
yields the identity
\[
e^{\,\mathrm{i}\,\bmmu\cdot\bmx}\,(Bv)(\bmmu,\bmx)
=
e^{\,\mathrm{i}\,\tilde\bmmu\cdot\bmx}\,(Bv)(\tilde\bmmu,\bmx),
\]
so the inversion formula \eqref{eq:inverse} applies verbatim to all
$\bmmu\in\R^2$. We shall restrict only to values of  $\bmmu \in \mathbb R^2\setminus \mathbb Z^2$, so we  regard $\mathcal{L}_{\bmmu, \epsilon}$ as an  unbounded operator on 

\begin{align}\label{Xper.def}
    X_{\rm per} &\coloneqq H^1(\T^2) \times L^2(\T^2) 
    \quad \mbox{ with domain  } \quad  D(\mathcal{L}_{\bmmu,\e}) = \{  U \in X_{\rm per} \colon \mathcal{L}_{\bmmu,\e}U \in X_{\rm per} \} \ .
\end{align}

We collect some properties of the operator $\cL_{\bmmu, \e}$ in the following. 

\begin{lemma} \label{lem:BlochFloquet}
For any $\e \geq 0$ and for any $\bmmu \in \R^2 \setminus \Z^2$, the Bloch-Floquet operator $\mathcal{L}_{\bmmu,\e}$: $D(\mathcal{L}_{\bmmu,\e}) \to 
X_{\rm per}$ in \eqref{eq:Lmueps} is closed and densely defined.  
In addition, 
  \begin{itemize}
    \item[(i)] 
    $\sigma_{X_{\rm per}}  (\mathcal{L}_{-\bmmu, \epsilon}) = \overline{ \sigma_{X_{\rm per}} (\mathcal{L}_{\bmmu, \epsilon})}.$
        \item[(ii)]  $\mathcal{L}_{\bmmu, \epsilon}$ is complex reversible, namely 
\begin{equation}\label{Lmue.rev}
    \mathcal{L}_{\bmmu, \epsilon} \circ \bar \varrho = - \bar \varrho \circ \mathcal{L}_{\bmmu, \epsilon}
\end{equation}
where $\bar \varrho$ is the complex involution 
\begin{equation}\label{cinvolution}
   \bar  \varrho \begin{bmatrix} b(\bmx) \\ \omega(\bmx) \end{bmatrix} := 
    \begin{bmatrix}   \bar b(-\bmx) \\ \bar \omega(-\bmx) \end{bmatrix}.
\end{equation}
\item[(iii)] The following inclusion holds:
\begin{equation}\label{spectrum}
\sigma_X(\cL_\e) \supseteq \bigcup_{\bmmu \in \R^2 \setminus \Z^2 } \sigma^p_{X_{\rm per}}(\cL_{\bmmu, \e}), \ 
\end{equation}
where $\sigma^p_{X_{\rm per}}(\cL_{\bmmu, \e})$ denotes the point spectrum of $\cL_{\bmmu, \e}$. 
\end{itemize}
\end{lemma}

\begin{proof}

We write $U= (b,\omega)^T \in X_{\rm per}$ and take $\bmmu \in \R^2 \setminus \Z^2$. Adapting the energy estimate \eqref{energy_est} to $\cL_{\bmmu, \e}$, one shows that 
$$
\norm{U(t)}_{X_{\rm per}} \lesssim \norm{U(0)}_{X_{\rm per}} \, \exp\left(C(\bmmu) t (\| b_\e\|_{W^{2, ^\infty}}+\| \mathbf u_\e\|_{W^{2, \infty}}) \right),
$$
where $C(\bmmu)>0$ is a positive constant,  
so that the  linear system
  $  \partial_t U = \cL_{\bmmu,\e} U$
defines a strongly continuous one-parameter semigroup  on $X_{\rm per}$ and then 
$\cL_{\bmmu,\e}: D(\cL_{\bmmu,\e}) \subset X_{\rm per} \to X_{\rm per}$ is densely defined and closed (see e.g. \cite[Theorem 4, Chapter 34.1]{lax2014functional}).

Item (i) follows from the general property that, if $\mathcal{A}$ is a real operator, then $\overline{\mathcal{A}_{\bmmu}} = \mathcal{A}_{- \bmmu}$, whereas 
item (ii) is a direct computation. \\
To prove item (iii) we argue as in Appendix A.2 of \cite{melinand2024phase} (see also Sec.2 of \cite{jiao2025small}). 
Recall that we consider $\mathcal{L}_{\e}$ as an operator on $X = \dot H^1(\R^2) \times L^2(\R^2)$ with domain $D(\mathcal{L}_{\e})$.
Assume that $\lambda_0 \in \sigma^p_{X_{\rm per}}({\cL_{\bmmu_0,\e}})$ for some $\bmmu_0 \in \R^2 \setminus \Z^2$, and denote by $U_0$ a corresponding eigenvector in 
$D(\cL_{\bmmu_0, \e})$.
Take   $\delta>0$  sufficiently small so that the ball $B(\bmmu_0,\delta)$ of center $\bmmu_0$ and radius $\delta$ does not intersect $\Z^2$. 
Note that,  by \eqref{eq:Lmueps},  one has 
\begin{equation}\label{diff.Lmu}
\sup_{\bmmu \in B(\bmmu_0,\delta)} \norm{\cL_{\bmmu,\e}- \cL_{\bmmu_0,\e}}_{X_{\rm per} \to X_{\rm per} }  \lesssim \delta \ .   
\end{equation}
This implies that
$U_0 \in D(\cL_{\bmmu, \e})$ for any $\bmmu \in B(\bmmu_0,\delta)$, and 
 we define 
\begin{equation*}
    U_{0,\delta}(\bmx) \coloneqq \int_{B(\bmmu_0,\delta)} e^{\mathrm{i} \bmmu\cdot\bmx} U_0(\bmx) \mathrm{d}\bmmu \ . 
\end{equation*}
By the above definition, the vector 
$U_{0,\delta}$ has Bloch-Floquet transform 
$ \mathbbm{1}_{B(\bmmu_0,\delta)} U_0$ supported, in the $\bmmu$-variable, away from any vector of integers as $\bmmu_0 \in \mathbb R^2\setminus \mathbb Z^2$, so that by \eqref{B.X}
$$\norm{U_{0,\delta}}_X \simeq \norm{\mathbbm{1}_{B(\bmmu_0,\delta)} U_0}_{L^2([-\frac12, \frac12)^2; X_{\rm per})} \lesssim \norm{U_0}_{X_{\rm per}} < \infty.$$
Analogously, by \eqref{eq:BlochFloquet},  the vector 
$\mathcal{L}_{\e} U_{0,\delta}$ has Bloch-Fourier transform supported, in the $\bmmu$-variable, away from zero, so that by applying first \eqref{B.X}, using that $U_0$ is an eigenvector of $\cL_{\bmmu_0, \e}$ and then  \eqref{diff.Lmu}, we have
\begin{align*}
\frac{\| ( \mathcal{L}_{\e} - \lambda_0 ) U_{0,\delta} \|_{ X } }{ \| U_{0,\delta} \|_{ X } } 
\approx 
\frac{
\norm{\mathbbm{1}_{B(\bmmu_0,\delta)} \left( \cL_{\bmmu,\e}- \cL_{\bmmu_0,\e} \right) U_0 }_{L^2([-\frac12, \frac12)^2; X_{\rm per})}
}
{\norm{\mathbbm{1}_{B(\bmmu_0,\delta)}U_0}_{L^2([-\frac12, \frac12)^2; X_{\rm per})}}
\stackrel{\delta \to 0}{\to} 0 ,
\end{align*}
hence $U_{0,\delta} \in D(\cL_\e)$ and $\lambda_0 \in \sigma_{ X }(\mathcal{L}_{\e})$.

\end{proof}


\begin{remark} \label{rem:FloquetFourier}

As an alternative to \eqref{eq:BlochFloquet}, since the operator $\cL_\e$ depends only on the variable $\ubmk \cdot \bmx$,  one may apply a  Fourier transform in the transverse direction $\ubmk^{\perp}$ and a Bloch--Floquet decomposition in the direction $\ubmk$, 
following Sec.~2.3 of \cite{jiao2025small}. 
This leads to an operator $\mathcal{L}(\alpha,\mu,\epsilon)$, for $(\alpha,\mu) \in \mathbb{R}^2 \setminus (\{0\} \times \mathbb{Z})$ and $\epsilon \ge 0$, whose spectrum in the unperturbed case $\epsilon=0$ is given by a formula analogous to \eqref{eq:lambdasigmak-restr}. Further details are provided in Appendix~\ref{app:FBF}.
We prefer to keep the current formulation since it allows us to perform the computations in Sec. 4 by using the more familiar form \eqref{eq:BIISysNew} of the Boussinesq system in Cartesian coordinates.

\end{remark}


This reformulation reduces the problem to analyze the point spectrum of $\mathcal{L}_{\bmmu,\epsilon}$ on $X_{\rm per}$ for different values of $\bmmu \in \mathbb R^2\setminus \mathbb Z^2$. We start with the simpler case $\e =0$.

\medskip
{\bf The spectrum of the linear operator $\mathcal{L}_{\bmmu, 0}$.} When $\e=0$ and $\bmmu \in \R^2   \setminus  \Z^2$, the operator $\cL_{\bmmu,0}$ in \eqref{eq:Lmueps} reduces to the matrix of Fourier multipliers
\begin{align} \label{eq:Lmu0}
\cL_{\bmmu,0} &\coloneqq
\begin{pmatrix}
\underline{\bmc} \cdot  (\nabla+ \im \bmmu) & -N^2 \,  (\partial_x + \im \mu_{1}) \, \left| \nabla + \im \bmmu \right|^{-2} \\
- (\partial_x+ \im \mu_1) & \underline{\bmc} \cdot (\nabla + \im \bmmu)
\end{pmatrix}
.
\end{align}
Such operator is complex Hamiltonian, meaning
\begin{equation}
  \cL_{\bmmu,0} = \cJ({\bmmu}) \,  \cB({\bmmu})
\end{equation}
where $\cJ({\bmmu})$ is the skew-adjoint operator
\begin{equation}\label{Jmu}
    \cJ({\bmmu}):=
    \begin{pmatrix}
        0 & \pa_x + \im \mu_1 \\
        \pa_x + \im \mu_1 & 0 
    \end{pmatrix}  \ , \quad \cJ({\bmmu})^* = -\cJ({\bmmu}) \ , 
\end{equation}
and $ \cB({\bmmu})$ is the selfadjoint operator
\begin{equation}
     \cB({\bmmu}):=
    \begin{pmatrix}
        -1 & (\pa_x + \im \mu_1)^{-1} ( \underline{\bmc} \cdot ( \grad + \im \bmmu ) ) \\
       (\pa_x + \im \mu_1)^{-1} ( \underline{\bmc} \cdot ( \grad + \im \bmmu ) )   &  -N^2 |\grad + \im \bmmu|^{-2}
    \end{pmatrix}  \ , \quad \cB({\bmmu})^* = \cB({\bmmu}) \ .
\end{equation}
Note that the operator $\cJ({\bmmu})$ is reversible, namely
\begin{equation}\label{Jmu.rev}
    \cJ({\bmmu}) \circ \bar \varrho = - \bar \varrho \circ  \cJ({\bmmu}) , 
\end{equation}
where $\bar \varrho$ is the complex involution defined in \eqref{cinvolution}.

The spectrum $\sigma_{X_{\rm per}}(\mathcal{L}_{\bmmu,0})$ is purely imaginary and composed  by the eigenvalues 
\begin{align}
\lambda^{\sigma}_{\bmk}(\bmmu) &\coloneqq \mathrm{i} \, w^{\sigma}_{\bmk}(\bmmu) , \quad w^{\sigma}_{\bmk}(\bmmu) \coloneqq \underline{\bmc} \cdot (\bmk + \bmmu) + \sigma \, \Omega(\bmk+\bmmu) , 
\quad \bmk \in \Z^2 , \quad \sigma = \pm ,
\label{eq:lambdasigmak}
\end{align}
where the dispersion relation 
\begin{align} \label{eq:DispRel}
 \Omega(\bmk) \, := \,  \, N \frac{ k_1}{|\bmk|} , \quad \bmk = (k_1, k_2)^T \in \Z^2 \setminus \{ \bmzero \} , \qquad 
 \Omega(\bmzero) := 0 \ .
 \end{align}
The corresponding eigenvectors 
$\{  f_{\bmk}^\sigma(\bmmu) \}_{\bmk \in \Z^2, \sigma = \pm}$, forming a basis of $X_{\rm per}$, are given by 
\begin{equation}\label{eq:v20}
    f_{\bmk}^\sigma(\bmmu) :=
    \frac{1}{\sqrt{2}}
    \begin{pmatrix}
    N \\ -\sigma  | \bmk + \bmmu |
    \end{pmatrix}
    e^{\mathrm{i}{\bmk} \cdot \bmx} , 
\qquad \cL_{\bmmu, 0}  f_{\bmk}^\sigma(\bmmu)  = \im w_{\bmk}^\sigma(\bmmu) \,   f_{\bmk}^\sigma(\bmmu) \ , \quad \bmk \in \Z^2, \sigma = \pm \ .
\end{equation}
Note that the eigenvalues \eqref{eq:lambdasigmak} satisfy the following relations
\begin{align} 
\lambda^{\sigma}_{\bmk}(\bmzero) &= - \lambda^{\sigma}_{-\bmk}(\bmzero) , \;\; \forall \sigma \in \{\pm \}, \;\; \forall \bmk \in \mathbb{Z}^2 \setminus \{ \bmzero \} , \label{eq:lambdaSym} \\
\lambda^{-}_{\ubmk}(\bmzero) &= \lambda^{-}_{-\ubmk}(\bmzero) = 0 . \label{eq:TrivColl} 
\end{align}

Moreover, denoting the complex scalar product 
\begin{align*}
    \langle f,g \rangle &\coloneqq \frac{1}{(2\pi)^2} \int_{\mathbb{T}^2} (f_1 \overline{g_1} + f_2 \overline{g_2} ) \mathrm{d}\bmx , \quad \forall f = (f_1,f_2)^T, g = (g_1,g_2)^T \in {L^2(\T^2)\times L^2(\T^2)},
\end{align*}
and  $\| f \| \coloneqq \langle f,f \rangle^{1/2}$, 
the eigenvectors  $\{  f_{\bmk}^\sigma(\bmmu) \}_{\bmk \in \Z^2, \sigma = \pm}$ satisfy the orthogonality conditions
\begin{align*}
    \langle f_{\bmk}^\sigma(\bmmu) ,   f_{{\bmk}'}^{\sigma'}(\bmmu) \rangle &= 0 , \quad \mbox{ if } \bmk \neq \bmk' , \quad \forall \sigma, \sigma' \in \{ \pm \}, \\
    \langle f_{\bmk}^\sigma(\bmmu) ,   f_{{\bmk}}^{\sigma'}(\bmmu) \rangle &= \frac{1}{2} ( N^2 + \sigma \sigma' |\ubmk + \bmmu|^2 ) , \quad \forall \bmk \in \Z^2 , \quad \forall \sigma, \sigma' \in \{ \pm \}     \ , 
\end{align*}

the symplectic conditions
\begin{equation}\label{symplfjs}
    \langle   \cJ(\bmmu) f_{\bmk}^\sigma(\bmmu) ,   f_{{\bmk}'}^{\sigma'}(\bmmu) \rangle 
    = \begin{cases}
        -  \im |\bmk + \bmmu|^2 \Omega(\bmk + \bmmu) & \mbox{ if } \bmk = {\bmk}' \mbox{ and } \sigma = \sigma'  = +\ , \\
       \ \  \im |\bmk + \bmmu|^2 \Omega(\bmk + \bmmu)  & \mbox{ if } \bmk = {\bmk}' \mbox{ and } \sigma = \sigma'  = -\ , \\
        0 & \mbox{ otherwise }
    \end{cases}
\end{equation}
and  the reversibility property
\begin{equation}\label{base.rev}
    \bar \varrho f_{\bmk}^\sigma(\bmmu) = f_{\bmk}^\sigma(\bmmu) , \qquad \bar \varrho \mbox{ in } \eqref{cinvolution} \ . 
\end{equation}


Since $\cL_{\bmmu,0}$ is complex Hamiltonian, it is natural to  look for unstable 
eigenvalues of $\mathcal{L}_{\bmmu,\epsilon}$ bifurcating from multiple eigenvalues of
$\cL_{\bmmu,0}$.
In the setting of  Hamiltonian system this is a classical approach,  connected with the instability associated to the Krein signature (see e.g. \cite{kapitula2004counting,kapitula2013spectral} and reference therein). 
Although in our setting we could not find a nice Hamiltonian structure for the operator $\mathcal{L}_{\bmmu,\e}$ when $\e \neq 0$, we take in any case this approach and look for  multiple eigenvalues  \eqref{eq:lambdasigmak} of the form 
\begin{align} \label{eq:Rk01}
\lambda^{-}_{\underline{\bmk}}(\bmmu) = \lambda^{+}_{\bmzero}(\bmmu), \quad \text{or equivalently,} \quad \Omega(\underline{\bmk}) - \Omega(\underline{\bmk}+\bmmu) = \Omega(\bmmu).
\end{align}  
Such condition is  directly  inspired by the so-called Triadic Resonant Instability (TRI)  and   Parametric Subharmonic Instability (PSI), where the energy is transferred from a primary wave to two secondary waves, which in the physics literature is  identified as the key resonance mechanism leading to instability, see e.g.  \cite[Section 3]{dauxois2018instabilities} and \cite[Section 5]{bianchini2023triadic} for a discussion of such a phenomenon for a Boussinesq system with viscosity, and references therein.

Imposing the resonance \eqref{eq:Rk01} is tantamount to take $\bmmu$ in the  resonant set  
\begin{align}\label{resonant.set}
\mathcal{R}_{\ubmk} &\coloneqq \{ \bmzero \} \cup \left\{ \bmmu \in  \mathbb{R}^2 \setminus \{ \bmzero \}  : \Omega(\ubmk) - \Omega(\ubmk+\bmmu) = \Omega(\bmmu) \right\} \ , 
\end{align}  

which enforces  $\cL_{\bmmu,0}$ to have  an eigenvalue with algebraic multiplicity  {\em at least} 2. 
However, such multiplicity might be higher and moreover such eigenvalues might not be isolated, which makes the analysis more complicated. To avoid this, in the next paragraph we shall restrict $\cL_{\bmmu,0}$ to an invariant subspace where $\lambda^{-}_{\underline{\bmk}}(\bmmu) = \lambda^{+}_{\bmzero}(\bmmu)$ is an isolated eigenvalue with multiplicity exactly 2. 
Before doing that, it is useful to study analytically the resonant set $\mathcal{R}_{\ubmk}$: in the next lemma we show that, under a technical assumption on $\ubmk$, $\cR_{\ubmk}$ consists of two analytic branches through $\bmmu=\bmzero$, with a possible non-smooth point on the negative branch if $2\frac{\tm}{|\ubmk|} < 1$.

\begin{figure}[h!]
    \centering
     \begin{subfigure}[b]{0.45\textwidth}
    \includegraphics[width=\linewidth]{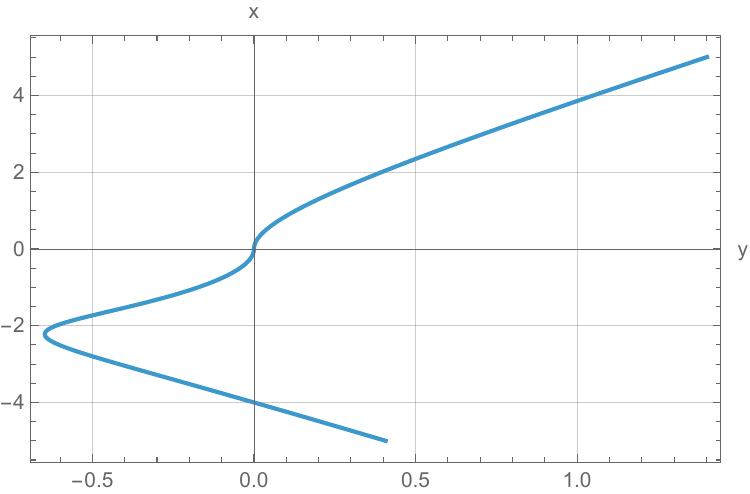}
    \caption{The set $\cR_{\ubmk}$ in case $\ubmk = (2,2)$}
    \label{fig:subfig1}
  \end{subfigure}
  \hfill
  \begin{subfigure}[b]{0.45\textwidth}
    \includegraphics[width=\linewidth]{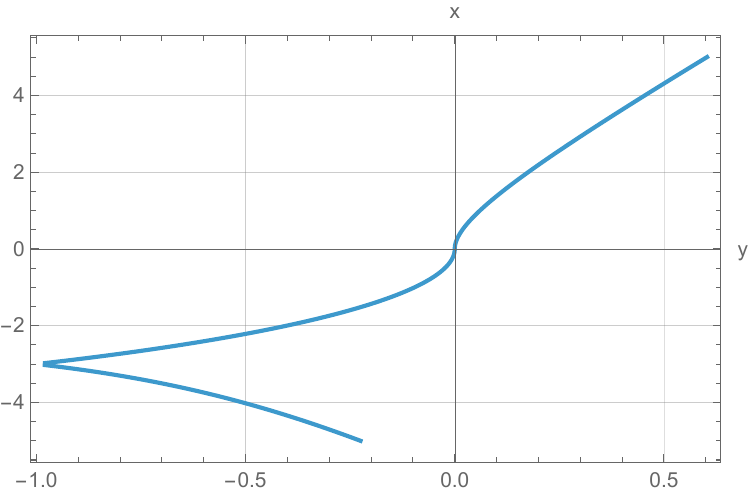}
    \caption{The set $\cR_{\ubmk}$ in case $\ubmk = (1,3)$}
    \label{fig:subfig2}
  \end{subfigure}
\label{fig}
\end{figure}

\begin{lemma}\label{lem:Rk}
Let $\ubmk = (\tm, \tn)^T \in \Z^2$ with $\tm > 0$, $\tn > 0$. The set $\cR_{\ubmk}$ in \eqref{resonant.set} decomposes as
\[
\cR_{\ubmk} = \cR_{\ubmk}^+ \cup \{\mathbf{0}\} \cup \cR_{\ubmk}^-,
\]
with 
\[
\cR_{\ubmk}^\pm = \{ (\varphi_\pm(y), y) \colon y \gtrless 0 \},
\]
i.e., both $\cR_{\ubmk}^\pm$ are the graphs of functions $\varphi_\pm(y)$. Moreover:
\begin{itemize}
    \item $\varphi_+$ is real analytic in $(0, +\infty)$ and strictly increasing;
    
    \item $\varphi_-$ satisfies $\varphi_-(-2\tn) = 0$, and:
    \begin{itemize}
        \item If $2\tfrac{\tm}{|\ubmk|} > 1$, then $\varphi_-$ is real analytic in $(-\infty, 0)$;
        \item If $2\tfrac{\tm}{|\ubmk|} < 1$, then $\varphi_-$ is analytic on $(-\infty, 0)\setminus\{-\tn\}$, Lipschitz at $y = -\tn$, and satisfies the asymptotic behavior
        \begin{equation}\label{varpi.exp2}
            \varphi_-(y) \sim -\tm + a_* |\tn + y| \quad \text{as } y \to -\tn, \quad
            a_* := \sqrt{\frac{4\tm^2}{\tn^2 - 3\tm^2}}.
        \end{equation}
    \end{itemize}
\end{itemize}

Moreover, the functions $\varphi_\pm$ satisfy the following asymptotics as $y \to 0^\pm$ and $y \to \pm\infty$:
\begin{equation}\label{varpi.exp}
    \varphi_\pm(y) \sim \pm \frac{\tm \tn}{|\ubmk|^3} y^2 \quad \text{as } y \to 0^\pm, \qquad 
    \varphi_\pm(y) \sim \pm \frac{\tm}{\sqrt{3\tm^2 + 4\tn^2}} \left(\frac{\tn}{2} + y \right) - \frac{\tm}{2} \quad \text{as } y \to \pm\infty.
\end{equation}

See Figure \ref{fig:subfig1} and \ref{fig:subfig2}.

Finally, one has:
\begin{equation}\label{linR}
    \ell \ubmk \in \cR_{\ubmk},\ \ell \in \Z \quad \Rightarrow \quad \ell \in \{-1, 0\}.
\end{equation}
\end{lemma}
We postpone the proof to Appendix~\ref{app:Rk}.\\

\begin{remark}
    Notice that in the regime of high Floquet parameter $y \rightarrow \pm \infty $ considered in Lemma~\ref{lem:Rk}, we have
\begin{align*}
    \Omega (\ubmk + \bmmu) &= \frac{\Omega (\ubmk)}{2} + \mathcal{O} \left( \frac{1}{y} \right), \quad
    \Omega (\bmmu) = \frac{\Omega (\ubmk)}{2} + \mathcal{O} \left( \frac{1}{y} \right),
\end{align*}
so that, to leading order, the oscillation frequency of both secondary waves - corresponding to wavevectors  $\ubmk + \bmmu$ and $- \bmmu$ - is approximately half that of the primary wave. This corresponds precisely to the regime of Parametric Subharmonic Instability (PSI), which is particularly relevant for inviscid fluids (see~\cite[Section~3.1]{dauxois2018instabilities}).
\end{remark}

\medskip
{\bf The invariant subspace $X_{\ubmk}^0$.}
As already commented, for $\bmmu \in \cR_{\ubmk}$ we cannot guarantee  the eigenvalue in \eqref{eq:Rk01} to be isolated with multiplicity  2. 
Therefore, inspired by Remark \ref{rem:FloquetFourier}, we exploit that 
the operator $\cL_{\bmmu, \e}$ depends only on the variable $\ubmk \cdot \bmx$ to perform a further reduction and introduce an invariant lower-dimensional subspace where we prove 
that \eqref{eq:Rk01} is an isolated eigenvalue with multiplicity  2. 
For $\ell \in \N\cup\{0\}$ we define the spaces

 \begin{equation}\label{Hlk}
    X^{\ell}_{\ubmk}:= \left\lbrace \vect{f_1}{f_2} \in {H}^{\ell+1}(\T^2) \times {H}^{\ell}(\T^2) \colon f_j(\bmx) = \sum_{n \in \Z  } f_{j,n} e^{\im n \ubmk \cdot \bmx} \right \rbrace \ ,
\end{equation}

 consisting of functions whose Fourier coefficients are supported only on wavevectors  proportional to the primary traveling wave wavevector \( \ubmk \), i.e., the wavevectors  \( \ubmk\Z :=\{ n \ubmk \colon n \in \Z\} \).

We also define the spaces $X^{\ell,\perp}_{\ubmk}$ via the orthogonal decompositions

\begin{equation*}
    {H}^{\ell+1}(\T^2) \times {H}^{\ell}(\T^2) =  X^{\ell}_{\ubmk} \oplus X^{\ell,\perp}_{\ubmk} \ , \quad X^{\ell,\perp}_{\ubmk}:= \left( X^{\ell}_{\ubmk} \right)^\perp.
\end{equation*}

Recall that, given  an operator $A$ on its domain
$D(A) \subset H$, $H$ a Hilbert space, 
a closed subspace $V\subseteq  H$ is said to be {\em invariant} if $A$ maps $D(A) \cap V$ into $V$; in this case we denote by $A\vert_V$ the restriction  $A\vert_V\colon D(A)\cap V \to V$. We also denote by $P_V\colon H \to V$ the projection of $H$ on the closed subspace $V$.

\begin{lemma}\label{Lmue.prop}
    For any $\e \geq 0$ and  $\bmmu \in \R^2 \setminus \Z^2$, the closed subspaces
    $X^0_{\ubmk}, \ X^{0, \perp}_{\ubmk} \subset X_{\rm per}$ are invariant, and 
\begin{equation}\label{spec.dec}
\sigma_{X_{\rm per}}(\cL_{\bmmu, \e}) = \sigma_{X^0_{\ubmk} }( \sL_{\bmmu,\e}) \cup \sigma_{X^{0,\perp}_{\ubmk} }(\sL_{\bmmu,\e}^\perp)  \ 
\end{equation} 
where 
  \begin{equation}\label{L0}
    \sL_{\bmmu,\e} \coloneqq \cL_{\bmmu, \e}\vert_{X^0_{\ubmk}}\colon D(\cL_{\bmmu, \e}) \cap X^0_{\ubmk} \to X^0_{\ubmk} \quad 
   \mbox{ and }
   \quad
   \sL_{\bmmu,\e}^\perp \coloneqq \cL_{\bmmu, \e}\vert_{X^{0,\perp}_{\ubmk}}\colon D(\cL_{\bmmu, \e}) \cap X^{0,\perp}_{\ubmk} \to X^{0,\perp}_{\ubmk} \ . 
  \end{equation}
\end{lemma}
\begin{proof}
Since  the functions \( \psi_\e, \omega_\e, b_\e \) in Proposition~\ref{prop:trav} depend only on $\ubmk \cdot \bmx$,
   for any vector $U \in  D(\cL_{\bmmu, \e})\cap X_{\ubmk}^0$, one has
$   \cL_{\bmmu, \e}U \in X^0_\ubmk $, hence  $X^0_\ubmk $ is invariant.  Similarly, $X^{0, \perp}_\ubmk $ is invariant. 
We verify now that the  projection $P_0:= P_{X^0_\ubmk }\colon X_{\rm per}\to X^0_\ubmk $ maps the domain $D(\cL_{\bmmu, \e})$ in \eqref{Xper.def} to itself. 
Then let $U \in  D(\cL_{\bmmu, \e})$. Again exploiting that 
 \( \psi_\e, \omega_\e, b_\e \) depend only on $\ubmk \cdot \bmx$,
 one has $\cL_{\bmmu, \e}P_0 U = P_0 \cL_{\bmmu, \e} U \in X^0_\ubmk \subseteq X_{\rm per} $, showing that $P_0 U \in D(\cL_{\bmmu, \e})$. 
 It follows by \cite[Proposition 1.15]{S} that $X^0_\ubmk$ is a reducing subspace for $\cL_{\bmmu, \e}$, namely  
   $\cL_{\bmmu, \e}$ decomposes as $\cL_{\bmmu, \e}=  \sL_{\bmmu,\e} \oplus  \sL_{\bmmu,\e}^\perp$ and therefore 
  \eqref{spec.dec} follows. 
\end{proof}

Our next steps are to show that $(i)$  the eigenvalues in \eqref{eq:Rk01} are isolated eigenvalues of multiplicity exactly two of  $\sL_{\bmmu, 0}$ 
whenever $\bmmu \in \cR_{\ubmk}$ (recall \eqref{resonant.set}), and $(ii)$ that for $\e\neq 0$ sufficiently small, such double eigenvalue bifurcates in two distinct eigenvalues $\lambda^\pm(\bmmu, \e)$, with $\lambda^+(\bmmu, \e)$ with strictly positive real part. 
In addition, we shall prove that the corresponding eigenfunction belongs to $D(\cL_{\bmmu, \e})$, 
therefore
$\lambda^\pm(\bmmu, \e) \in \sigma^p_{X_{\rm per}}(\cL_{\bmmu, \e})$. 
Hence, our instability result follows by the spectral inclusion in \eqref{spectrum}.

\medskip

{\bf Double isolated eigenvalue of $\sL_{\bmmu, 0}$.} 
Remark that  $\sigma_{X^0_{\ubmk} }(\sL_{\bmmu, 0})$ consists in those eigenvalues in \eqref{eq:lambdasigmak}
having eigenvectors $f_{\bmk}^\sigma(\bmmu) \in D(\cL_{\bmmu, \e})\cap X_\ubmk^{(0)}$, namely 
\begin{align}
\lambda^{\sigma}_{n \ubmk}(\bmmu) &\coloneqq \mathrm{i} \, w^{\sigma}_{n \ubmk}(\bmmu) , \quad w^{\sigma}_{n \ubmk}(\bmmu) \coloneqq \underline{\bmc} \cdot (n \ubmk + \bmmu) + \sigma \, \Omega(n \ubmk+\bmmu) , 
\quad n \in \Z , \sigma = \pm.
\label{eq:lambdasigmak-restr}
\end{align}
We shall now prove that whenever $\bmmu \in \cR_{\ubmk} \cap ( \R^2 \setminus \Z^2 ) $ (except for a discrete set)
the double eigenvalue $\lambda^{-}_{\underline{\bmk}}(\bmmu) = \lambda^{+}_{\bmzero}(\bmmu)$ is isolated and with  multiplicity 
\emph{exactly} 2.

\begin{proposition} \label{prop:wBound}
For any wavevector $\ubmk=(\tm, \tn)^T \in \Z^2$ with $\tm, \tn >0  $, there exists a discrete set $\cD_\ubmk \subset \cR_\ubmk $, 
with $\bmzero \in \cD_{\ubmk}$,   such that
for any $\bmmu \in \cR_\ubmk\setminus \cD_\ubmk $, 
there exists $\alpha=\alpha (\bmmu)>0$ such that:
\begin{align} \label{eq:wEst}
\inf_{\ell \in \mathbb{Z}, \ell \neq 0} | w^{+}_{\ell \ubmk}( \bmmu ) - w^{+}_{\bmzero}( \bmmu )| &\geq \alpha N ,  
\qquad 
\inf_{ \ell \in \mathbb{Z}, \ell \neq 1} | w^{-}_{\ell \ubmk}( \bmmu ) - w^{-}_{\underline{\bmk}}( \bmmu ) | \geq \alpha N , 
\end{align}
where $N$ is the Brunt--V\"ais\"al\"a frequency defined in \eqref{eq:relNrho}.

In particular, whenever
$\bmmu \in ( \cR_\ubmk\setminus \cD_\ubmk ) \cap ( \R^2 \setminus \Z^2 ) $, the operator $\sL_{\bmmu, 0}$ in \eqref{L0}  has the isolated double eigenvalue 
$\lambda^{-}_{ \ubmk}(\bmmu) = \lambda^+_{\bmzero}(\bmmu)$.
\end{proposition}
\begin{proof}
See Appendix \ref{app:de}.
\end{proof}

By Proposition \ref{prop:wBound}, for any $\bmmu \in ( \cR_\ubmk\setminus \cD_\ubmk ) \cap ( \R^2 \setminus \Z^2 )$, the spectrum of $\sL_{\bmmu, 0}$ decomposes in two separated parts
\begin{align} \label{eq:DecompSpecLmu0}
\sigma_{X^0_{\ubmk} }( \sL_{\bmmu, 0} )  &= \sigma^{\prime}( \sL_{\bmmu, 0} ) \cup \sigma^{\prime \prime}( \sL_{\bmmu, 0} ) , \quad \sigma^{\prime}( \sL_{\bmmu, 0} ) \coloneqq \{ \lambda^{+} \} , \quad \lambda^+:=\lambda^{-}_{\underline{\bmk}}(\bmmu) \equiv \lambda^{+}_{\bmzero}(\bmmu). 
\end{align}
With the notation of~\eqref{eq:v20}, we denote by $\mathcal{V}_{\ubmmu, 0}  \subset D(\cL_{\bmmu, 0})  \cap X_\ubmk^{0} $ the spectral subspace generated by the eigenvectors $\{f_{\ubmk}^-, f_{\mathbf{0}}^+\}$ associated with the eigenvalue $\lambda^+$ in~\eqref{eq:lambdasigmak-restr}.

\medskip

 By Kato's Similarity Transformation (Section~\ref{sec:Kato}, Lemma~\ref{lem:Kato}), we will show that whenever $\e>0$ is small enough, the spectrum of  $\sL_{\bmmu, \e}$ decomposes as well in two 
 parts
 \begin{align} \label{eq:DecompSpecLmueps}
 \sigma_{X^0_{\ubmk} }( \sL_{\bmmu, \e} )  = \sigma^{\prime}\big( \sL_{\bmmu, \e} \big) \cup \sigma^{\prime\prime}\big( \sL_{\bmmu, \e}\big),
\end{align}
where $\sigma^{\prime}\big( \sL_{\bmmu, \e} \big)$ consists of two eigenvalues close to $\lambda^+$.
 Denoting by $\mathcal{V}_{\bmmu, \epsilon}$ the spectral subspace associated with $\sigma^{\prime}\big( \sL_{\bmmu, \e}\big)$, we provide the following matrix representation of the operator
$\sL_{\bmmu, \e} \colon \mathcal{V}_{\bmmu, \epsilon} \to \mathcal{V}_{\bmmu, \epsilon}$.

\begin{theorem}\label{thm:expansion}
Let $\bmmu \in ( \cR_\ubmk\setminus \cD_\ubmk ) \cap ( \R^2 \setminus \Z^2 )$.
There exists  $\e_0 \equiv \e_0(\bmmu) >0$ such that for any $\e \in (0, \e_0)$, the operator
 $ \sL_{\bmmu,\e} \colon \cV_{\bmmu, \e} \to \cV_{\bmmu, \e}$ is represented by the 
 $2\times 2$ matrix:
 \begin{align}\label{}
   \begin{pmatrix}
       \lambda^+ + \im \, r_1(\e^2) & \im\,  \tb_0(\bmmu) \e + \im \, r_2(\e^2) \\
       \im \, \tb_1(\bmmu) \e +\im  r_3(\e^2) &  \lambda^+ + \im \, r_4(\e^2)
   \end{pmatrix}
 \end{align}
 where  $\lambda^+$ is in \eqref{eq:DecompSpecLmu0}  and 
 \begin{equation}
     \begin{aligned}
        & \tb_1(\bmmu):=  (\ubmk^{\perp} \cdot {\bmmu})  \frac{ \big( |\ubmk+{\bmmu}| -|\ubmk|\big) \, (|\ubmk + {\bmmu}| + |\ubmk| - |{\bmmu}| \big)}{4 N\,|\ubmk|\, |{\bmmu}|\,  |\ubmk+{\bmmu}|}, \\
& \tb_0(\bmmu):= - (\ubmk^{\perp} \cdot {\bmmu}) \frac{\big( |\ubmk| +|{\bmmu}|\big) \, \big( |\ubmk + {\bmmu}| + |\ubmk| - |{\bmmu}|\big)   }{4 N \, |\ubmk| \, |{\bmmu}| \, |\ubmk + {\bmmu}|}. 
     \end{aligned}
 \end{equation}
 Its eigenvalues are given by
 \begin{equation}\label{lambdapm}
     \lambda^\pm(\bmmu, \e)=  \lambda^+ + \im r(\e^2)  \pm   \e  \sqrt{\mathtt{e}(\bmmu) + r(\e)},
 \end{equation}
and
 \begin{align}\label{temu}
    \te(\bmmu)&= - \tb_1(\bmmu)\tb_0(\bmmu)=
    (\ubmk^\perp \cdot \bmmu)^2 \,  \frac{ (|\ubmk+\bmmu| - |\ubmk|)\, (|\ubmk| + |\bmmu|) \, (|\ubmk + \bmmu| + |\ubmk| - |\bmmu|)^2 }{16 N^2\,  |\bmmu|^2\,  |\ubmk + \bmmu|^2\, |\ubmk|^2}. 
\end{align}
In particular, 
 \begin{itemize}
     \item[(1)] in the regime of small Floquet parameter $|\bmmu| \ll 1$, writing 
     \begin{equation}\label{paranear0}
     \bmmu (\tau) = \begin{cases}
         (\varphi_+(\tn\tau), \tn \tau)^T & \mbox{ if } \tau >0 \\
         0 & \mbox{ if } \tau  = 0 \\
          (\varphi_-(\tn \tau), \tn \tau)^T & \mbox{ if } \tau <0
     \end{cases}
     \end{equation}
  we have
     \begin{equation}\label{eto0}
    \mathtt{e}(\bmmu(\tau)) = \frac{\tm^2 \tn^2}{4 N^2 |\ubmk|^2} \tau + \cO(\tau^2) \ , \quad \mbox{ as } \tau \to 0;
    \end{equation}
    \item[(2)] in the regime of high Floquet parameter $|\bmmu| \gg 1$, 
    writing 
     $$
     \bmmu (\tau) = 
         (\varphi_+(\tau), \tau)^T 
     $$
     we have 
     \begin{equation}\label{etoinf}
    \mathtt{e}(\bmmu(\tau)) = \frac{\tm^2}{256 N^2 |\ubmk|^6} (\tn - \sqrt{3 \tm^2+4\tn^2})^2 (3 \tm^2 + 2\tn^2 + \tn \sqrt{3 \tm^2 + 4\tn^2})^2 + \mathcal{O}(\tau^{-1}), \quad \mbox{ as } \tau \to \infty.
\end{equation}
\end{itemize}
\end{theorem}
The proof is in Section \ref{proof}.

\begin{remark} \label{rem:epsmu}

    The threshold $\e_0(\bmmu)$ in the statement of Theorem \ref{thm:expansion} comes from the statement of Lemma \ref{lem:Kato}, and it depends on the constant $\alpha=\alpha(\bmmu)$ in Proposition \ref{prop:wBound} describing the separation between the eigenvalue $\lambda^{+}$ and the rest of the spectrum. 
    From the proof of Proposition \ref{prop:wBound} (see Appendix \ref{app:de}) we have that the constant $\alpha=\alpha(\bmmu)$ in Proposition \ref{prop:wBound} is given by the infimum of a finite family of functions which are Lipschitz (actually, a.e. analytic) on $\cR_{\ubmk} \setminus \mathcal{D}_{\ubmk}$. This implies that the map
    \begin{align*}
        \cR_{\ubmk} \setminus \mathcal{D}_{\ubmk} &\to \mathbb{R}_{+} \, \quad 
        \bmmu \mapsto \alpha(\bmmu)
    \end{align*}
    is continuous. As a consequence, the thresholds $\e_0(\bmmu)$ in the statement of Theorem \ref{thm:expansion} can be chosen uniformly on compact subsets of $( \cR_\ubmk\setminus \cD_\ubmk ) \cap ( \R^2 \setminus \Z^2 )$.

\end{remark}

\section{Perturbative approach to separated eigenvalues} \label{sec:Kato}
The property \eqref{eq:DecompSpecLmueps} states that the spectrum of $ \sL_{\bmmu,\e}$ splits into two disjoint components. 
Following the approach in Kato's Similarity Transformation (see \cite{kato1966perturbation} and Sec. 3 of \cite{berti2022full}), we construct transformation operators which map isomorphically the unperturbed spectral subspace $\cV_{\bmmu,0}$ - generated by the two eigenvectors associated with the double eigenvalue $\lambda^+$ - into the perturbed one $\cV_{\bmmu,\e}$.
\begin{lemma} \label{lem:Kato}

Fix $\bmmu \in ( \cR_\ubmk\setminus \cD_\ubmk ) \cap ( \R^2 \setminus \Z^2 ) $. 
Let $\Gamma$ be a simple, closed, counterclockwise-oriented circle around $\lambda^{+}$ (see \eqref{eq:DecompSpecLmu0}) in the complex plane, separating 
$\sigma'({\sL}_{\bmmu,0}) = \{\lambda^{+}\}$ from the rest of the spectrum 
$\sigma''({\sL}_{\bmmu,0})$ in \eqref{eq:DecompSpecLmu0}. 
We further assume that the region enclosed by $\Gamma$ has diameter strictly less than $\alpha N$, 
where $\alpha$ is the positive constant in \eqref{eq:wEst} and $N$ is the Brunt--V\"ais\"al\"a frequency defined in \eqref{eq:relNrho}.

Then there exists $\e_0 \equiv \e_0(\bmmu)>0$ such that, for any $ \e \in (0, \e_0)$, the following statements hold true:
\begin{enumerate}
\item[1.] Recalling \eqref{L0}, the domain  of the operator $\sL_{\bmmu, \e}$ is given by
\begin{equation}\label{dom.sL}
    D(\sL_{\bmmu, \e})  =  X^1_\ubmk  \ , 
\end{equation}
and the curve $\Gamma$ belongs to the resolvent set of the operator $\sL_{\bmmu , \epsilon}:  X^1_\ubmk  \to X^0_\ubmk$ in \eqref{L0} ;

\item[2.] the operators
\begin{align} \label{eq:SpecProj}
P_{\bmmu , \epsilon} &\coloneqq -\frac{1}{2\pi\mathrm{i}} \oint_{\Gamma} ( {\sL}_{\bmmu , \epsilon} - \lambda )^{-1} \mathrm{d}\lambda : X^0_\ubmk \to X^1_\ubmk
\end{align}
are well-defined projectors commuting with ${\sL}_{\bmmu , \epsilon}$, namely
\begin{equation} \label{eq:propSpecProj}
P_{\bmmu , \epsilon}^2 = P_{\bmmu , \epsilon} , \;\; P_{\bmmu, \epsilon} {\sL}_{\bmmu , \epsilon} = \mathcal{\sL}_{\bmmu , \epsilon} P_{\bmmu, \epsilon} ,
\end{equation}
the map $ \epsilon \mapsto P_{\bmmu , \epsilon}$ is analytic from $ B(\epsilon_0)$ to $L(X^0_\ubmk,X^1_\ubmk)$ and
the projectors are reversibility preserving, i.e. $\bar \rho P_{\bmmu, \e} = P_{\bmmu, \e} \bar \rho$;
\item[3.] the domain $X^1_\ubmk$ of the operator ${\sL}_{\bmmu , \epsilon}$ decomposes as the direct sum
\begin{align} \label{eq:SpaceDecomp}
X^1_\ubmk &= \mathrm{Ker}( P_{\bmmu , \epsilon} ) \oplus \mathcal{V}_{\bmmu , \epsilon} , \;\; \mathcal{V}_{\bmmu , \epsilon} \coloneqq \mathrm{Rg}( P_{\bmmu , \epsilon} ) = \mathrm{Ker}( \mathrm{id} - P_{\bmmu , \epsilon} ) ,
\end{align}
where the closed subspaces $\mathrm{Ker}( P_{\bmmu , \epsilon} )$ and $\mathcal{V}_{\bmmu , \epsilon}$ are invariant under ${\sL}_{\bmmu , \epsilon}$,
\begin{equation*}
{\sL}_{\bmmu , \epsilon} : \mathcal{V}_{\bmmu , \epsilon} \to \mathcal{V}_{\bmmu , \epsilon} , \;\; {\sL}_{\bmmu , \epsilon} : \mathrm{Ker}( P_{\bmmu , \epsilon} ) \to \mathrm{Ker}( P_{\bmmu , \epsilon} ) .
\end{equation*}
Moreover,
\begin{align}
\sigma_{ X^0_\ubmk } ( {\sL}_{\bmmu , \epsilon} )  \cap \{ z \in \mathbb{C} \;\; \text{inside} \;\; \Gamma \} &= \sigma_{ X^0_\ubmk } (  {\sL}_{\bmmu , \epsilon} |_{\mathcal{V}_{\bmmu , \epsilon}} ) = \sigma^{\prime}( {\sL}_{\bmmu , \epsilon} ) , \nonumber \\
\sigma_{ X^0_\ubmk } ( {\sL}_{\bmmu , \epsilon} ) \cap \{ z \in \mathbb{C} \;\; \text{outside} \;\; \Gamma \} &= \sigma_{ X^0_\ubmk } (  {\sL}_{\bmmu , \epsilon} |_{ \mathrm{Ker}({P}_{\bmmu , \epsilon}) } ) = \sigma^{\prime \prime}( {\sL}_{\bmmu , \epsilon} ) , \label{eq:semicont}
\end{align}
proving the ``semicontinuity property" \eqref{eq:DecompSpecLmueps} of separated parts of the spectrum;

\item[4.] the subspaces $\mathcal{V}_{\bmmu , \epsilon}$ are isomorphic one to another.
In particular, $\mathrm{dim}( \mathcal{V}_{\bmmu , \epsilon} ) = \mathrm{dim}( \mathcal{V}_{{\bmmu} , 0} ) = 2$ for all $ \epsilon \in B(\epsilon_0)$.
\end{enumerate}
\end{lemma}
\begin{proof}

\begin{enumerate}
\item[1.] 
For any $\lambda \in \mathbb{C}$ we write $\sL_{\bmmu , \epsilon} - \lambda = \sL_{\bmmu , 0 } - \lambda + \mathcal{R}_{ \bmmu , \epsilon}$, where the operator 
$\mathcal{R}_{ \bmmu , \epsilon} \coloneqq \sL_{\bmmu , \epsilon} - \sL_{\bmmu , 0 }$ is explicitly given by 
\begin{align}
\mathcal{R}_{ \bmmu , \epsilon}= 
\begin{pmatrix}
- \gradp \psi_{\epsilon} \cdot (\nabla + \im \bmmu)  & - \nabla b_{\epsilon} \cdot (\gradp  + \im \bmmu^\perp) \, \vert {\bf D} + \bmmu \vert^{-2} \\
0 & 
 - \gradp \psi_{\epsilon} \cdot (\nabla + \im \bmmu)
 - \nabla \omega_{\epsilon} \cdot (\gradp  + \im \bmmu^\perp)  \vert {\bf D} + \bmmu \vert^{-2}
\end{pmatrix}
\label{eq:Rformula}
\end{align}
and  bounded  $ X^1_\ubmk \to X^0_\ubmk$ quantitatively as 
\begin{equation}\label{R.est}
    \norm{\mathcal{R}_{ \bmmu , \epsilon}}_{L(X^1_\ubmk,  X^0_\ubmk)} \lesssim \epsilon \ . 
\end{equation}

We now show that, for any $\lambda \in \Gamma$, the operator 
\[
\sL_{\bmmu , 0 } - \lambda \colon X^1_\ubmk \to X^0_{\ubmk}
\]
is invertible. We use that $(\sL_{\bmmu , 0 } - \lambda )^{-1}$ is a matrix of Fourier multipliers with associated symbol
\begin{equation*}
\left[(\mathrm{i} \underline\bmc\cdot (j \ubmk + \bmmu) -\lambda)^2 + \frac{N^2(j \tm+\mu_1)^2}{|j \ubmk + \bmmu|^2}\right]^{-1}
\begin{pmatrix}
\im \underline{\bmc} \cdot (j \ubmk + \bmmu)-\lambda  & \frac{\im N^2(j \tm +  \mu_1)}{|j\ubmk + \bmmu|^{2}} \\
\im (j \tm +  \mu_1) & 
\im  \underline{\bmc} \cdot (j \ubmk + \bmmu)-\lambda
\end{pmatrix},
\quad j \in \Z.
\end{equation*}
Using the notation
\[
c_j := \mathrm{i} \underline\bmc\cdot (j \ubmk + \bmmu) -\lambda, 
\qquad 
d_j := c_j^2 + \frac{N^2(j \tm+\mu_1)^2}{|j \ubmk + \bmmu|^2},
\]
and writing $U= (b,\omega)^T$, we prove that
\begin{align}
    \label{resolvent1}
    \|( \sL_{\bmmu , 0 } - \lambda )^{-1}U\|_{X^1_{\ubmk}} \le C \|U\|_{X^0_{\ubmk}}.
\end{align}
We treat the first component:
\begin{align*}
\|( (\sL_{\bmmu , 0 } - \lambda )^{-1}U)_1\|_{H^2(\mathbb T^2)}^2
&=
\sum_j \left|d_j^{-1} \left(c_j \widehat b_j + \frac{N^2(j\tm + \mu_1)}{|j\ubmk+\bmmu|^2}\widehat\omega_j\right)\right|^2 |j \ubmk+\bmmu|^4 \\
& \lesssim \sum_j \left(\left|c_j\right|^2 |\widehat b_j|^2 +|j\ubmk+\bmmu|^{-2}|\widehat\omega_j|^2\right) |d_j|^{-2}|j \ubmk+\bmmu|^4.
\end{align*}
Estimating the coefficients of $|\widehat b_j|^2$ and $|\widehat \omega_j|^2$, since $|d_j| \ge |c_j|^2$, we have
\[
|d_j|^{-2}|j\ubmk+\bmmu|^4|c_j|^2 \le \frac{|j\ubmk+\bmmu|^4}{|c_j|^2}; \quad \frac{|j\ubmk+\bmmu|^2}{|d_j|^2} \le \frac{|j\ubmk+\bmmu|^2}{|c_j|^4}.
\]
Thus it suffices to show that, uniformly in $j$,
\[
|c_j|^{-1} |j\ubmk+\bmmu|^2 \le C|j|; \quad |c_j|^{-2} |j\ubmk+\bmmu| \le C.
\]
Since $\lambda \in \Gamma$ and $\Gamma$ is contained in a ball of radius $\delta < \alpha N$ around $\lambda^+=\im \underline{\bmc} \cdot \bmmu + \im \Omega(\bmmu)$, we have, for a suitable $\delta=\delta(\bmmu)>0$, 
\begin{align*}|c_j|&=|\im \underline{\bmc} \cdot(j\ubmk+\bmmu)-\lambda| \ge ||j\underline{\bmc} \cdot \ubmk-\Omega(\bmmu)|-\delta| =| |j\Omega (\ubmk)-\Omega(\bmmu)|-\delta|\\
&=||(j-1)\Omega (\ubmk)+\Omega(\ubmk+\bmmu)|-\delta|, 
\end{align*}
where in the last equality we used \eqref{resonant.set}. 
%
Consequently, for all $j$,
the desired inequalities hold for some constant $C=C(\bmmu)>0$. This proves the estimate of the first component, while the estimate of the second component is similar and therefore is omitted.

Now, by \eqref{resolvent1} and \eqref{R.est}, we obtain that there exists a sufficiently small $\epsilon_0 >0$ such that for any $\epsilon \in  B(\epsilon_0)$, uniformly on the compact set $\Gamma$, the operator 
 \begin{equation*}
 ( \sL_{\bmmu , 0 } - \lambda )^{-1} \mathcal{R}_{\bmmu , \epsilon} : X^1_\ubmk \to X^1_\ubmk
 \end{equation*}
 is bounded, and has small operatorial norm.
Hence $\sL_{\bmmu , \epsilon} - \lambda$ is invertible by Neumann series, and
\begin{equation} \label{eq:InverseOp}
( \sL_{\bmmu , \epsilon} - \lambda )^{-1} = \left[ \mathrm{id} + ( \sL_{\bmmu , 0 } - \lambda )^{-1} \mathcal{R}_{\bmmu , \epsilon} \right]^{-1} ( \sL_{\bmmu , 0 } - \lambda )^{-1} : X^0_\ubmk \to X^1_\ubmk ,
\end{equation}
which proves that $\Gamma$ belongs to the resolvent set of $\sL_{\bmmu , \epsilon} \colon X^1_\ubmk \to X^0_\ubmk$.\\
We prove now \eqref{dom.sL}. Recall that $D(\sL_{\bmmu, \e}) := D(\cL_{\bmmu, \e}) \cap X^0_\ubmk$ (cf. \eqref{L0}), and that 
$D(\cL_{\bmmu, \e})$ is defined in \eqref{Xper.def}. 
Clearly if $U \in X^1_\ubmk$ then $\cL_{\bmmu, \e} U \in X^0_\ubmk\subset X_{\rm per}$, hence $X^1_\ubmk \subseteq D(\sL_{\bmmu, \e}) $. 
Take now $U \in D(\sL_{\bmmu, \e})  $, $\lambda \in \Gamma$ and write
$U = (\sL_{\bmmu , \epsilon} - \lambda)^{-1} V$, with $V := (\sL_{\bmmu , \epsilon} - \lambda) U$. Because  $U \in D(\sL_{\bmmu, \e})  $, $V \in X_{\rm per} \cap X^0_\ubmk = X^0_\ubmk$. Then \eqref{eq:InverseOp} implies that $U \in X^1_\ubmk$, concluding the proof of \eqref{dom.sL}. 

\item[2.] 

By part 1. of the statement the operator $P_{\bmmu , \epsilon}: X^0_\ubmk \to X^1_\ubmk$ is well-defined and bounded, and it clearly commutes with $\sL_{\bmmu , \epsilon}$; the projection property $P_{\bmmu , \epsilon}^2 = P_{\bmmu , \epsilon}$ is a classical result (see \cite{kato1966perturbation}). The map 
\begin{equation*}
\epsilon \mapsto ( \sL_{\bmmu , 0 } - \lambda )^{-1} \mathcal{R}_{ \bmmu , \epsilon} \in L(X^1_\ubmk)
\end{equation*}
is analytic; since for $\| T \|_{L(X^1_\ubmk)} < 1$ the map $T \mapsto (\mathrm{id}+T)^{-1}$ is analytic in $L(X^1_\ubmk)$, then also the maps
\begin{equation*}
\epsilon \mapsto ( \sL_{\bmmu , \epsilon } - \lambda )^{-1} \in L(X^0_\ubmk,X^1_\ubmk) , \;\;  \epsilon \mapsto  P_{\bmmu , \epsilon }  \in L(X^0_\ubmk,X^1_\ubmk)
\end{equation*}
are analytic.

\item[3.] 

Since the operator $P_{\bmmu , \epsilon}$ is continuous in $L(X^1_\ubmk)$, we can deduce the decomposition \eqref{eq:SpaceDecomp}. The invariance of the subspaces $\mathrm{Ker}( P_{\bmmu , \epsilon} )$ and $\mathcal{V}_{\bmmu , \epsilon}$ follows because $P_{\bmmu , \epsilon}$ and $\sL_{\bmmu , \epsilon}$ commute. To prove property \eqref{eq:semicont} we define for any $\lambda_0 \notin \Gamma$ the operator
\begin{align*}
R_{\bmmu , \epsilon}(\lambda_0) &\coloneqq -\frac{1}{2\pi\mathrm{i}} \oint_{\Gamma} \frac{1}{\lambda-\lambda_0} ( \sL_{\bmmu , \epsilon} - \lambda )^{-1} \mathrm{d}\lambda : X^0_\ubmk \to X^1_\ubmk .
\end{align*}
If $\lambda_0$ is outside $\Gamma$, we have 
\begin{align*}
R_{\bmmu , \epsilon}(\lambda_0)  ( \sL_{\bmmu ,\epsilon} - \lambda_0 ) &= ( \sL_{\bmmu ,\epsilon} - \lambda_0 ) R_{\bmmu , \epsilon}(\lambda_0) \,= \, P_{\bmmu , \epsilon} ,
\end{align*}
hence $\lambda_0 \notin \sigma(  \sL_{\bmmu , \epsilon} |_{\mathcal{V}_{\bmmu , \epsilon}} )$; on the other hand, if $\lambda_0$ is inside $\Gamma$, we have 
\begin{align*}
R_{\bmmu , \epsilon}(\lambda_0)  ( \sL_{\bmmu ,\epsilon} - \lambda_0 ) &= ( \sL_{\bmmu ,\epsilon} - \lambda_0 ) R_{\bmmu , \epsilon}(\lambda_0) \,= \, P_{\bmmu , \epsilon} - \mathrm{id} ,
\end{align*}
hence $\lambda_0 \notin \sigma(  \sL_{\bmmu , \epsilon} |_{ \mathrm{Ker}(\mathcal{P}_{\bmmu , \epsilon}) } )$.

\item[4.] 

By the definition \eqref{eq:SpecProj} of spectral projectors, by the formula \eqref{eq:Rformula} and by the identity $A^{-1} - B^{-1} = A^{-1} (B-A) B^{-1}$, we obtain
\begin{align*}
P_{\bmmu , \epsilon} - P_{\bmmu , 0} &= \frac{1}{2\pi\mathrm{i}} \oint_{\Gamma} ( \sL_{\bmmu , \epsilon } - \lambda )^{-1} \mathcal{R}_{ \bmmu , \epsilon} ( \sL_{\bmmu , 0 } - \lambda )^{-1} \mathrm{d}\lambda .
\end{align*}
Then for any $\epsilon \in B(\epsilon_0)$ we have $\| P_{\bmmu , \epsilon} - P_{\bmmu ,0} \|_{L(X^1_\ubmk)} < 1$, so that the transformation operators 
\begin{align} \label{eq:TransOp}
U_{\bmmu , \epsilon} &\coloneqq \left[ \mathrm{id} - (P_{\bmmu , \epsilon} - P_{\bmmu , 0})^2 \right]^{-1/2} \, \left[ P_{\bmmu , \epsilon}  P_{\bmmu , 0} + (\mathrm{id} - P_{\bmmu , \epsilon})(\mathrm{id} - P_{\bmmu , 0}) \right]
\end{align}
are well defined in $L(X^1_\ubmk)$; moreover, they are bounded and invertible in $X^1_\ubmk$ and in $X^0_\ubmk$, with inverse given by
\begin{align} \label{eq:InvTransOp}
U_{\bmmu , \epsilon}^{-1} &\coloneqq \left[ P_{\bmmu , 0}  P_{\bmmu , \epsilon} + (\mathrm{id} - P_{\bmmu , 0})(\mathrm{id} - P_{\bmmu , \epsilon}) \right] \, \left[ \mathrm{id} - (P_{\bmmu , \epsilon} - P_{\bmmu , 0})^2 \right]^{-1/2} , 
\end{align}
and they satisfy the property
\begin{equation} \label{eq:TransSpecProj}
U_{\bmmu , \epsilon} P_{\bmmu , 0} U_{\bmmu , \epsilon}^{-1} = P_{\bmmu , \epsilon} , \;\; U_{\bmmu , \epsilon}^{-1} P_{\bmmu , \epsilon} U_{\bmmu , \epsilon} = P_{\bmmu , 0} 
\end{equation}
(see Chapter I, Sec. 4.6 in \cite{kato1966perturbation} for a proof of \eqref{eq:InvTransOp} and \eqref{eq:TransSpecProj}). We also notice that since the map $\epsilon \mapsto P_{\bmmu , \epsilon} \in L(X^1_\ubmk)$ is analytic, then the map $\epsilon \mapsto U_{\bmmu , \epsilon} \in L(X^1_\ubmk)$ is also analytic. The thesis now follows from \eqref{eq:TransSpecProj}.

\end{enumerate}

\end{proof}

By Lemma~\ref{lem:Kato} we ensure that the corresponding perturbed eigenvalues remain isolated from the rest of the spectrum. We now consider the following basis for the perturbed spectral subspace
$\mathcal{V}_{{\bmmu},\epsilon}$:
\begin{align} \label{eq:basisF}
\mathcal{F} &\coloneqq \{ f_{\ubmk}^-(\bmmu , \epsilon) , f_{\bmzero}^+(\bmmu,\epsilon) \}, \quad f_{\bmj}^\sigma(\bmmu , \epsilon) \coloneqq P_{\bmmu , \epsilon} f_{\bmj}^\sigma(\bmmu) , \; \bmj=\ubmk,\bmzero,
\end{align}
given by the projection  $P_{\bmmu,\epsilon}$ in \eqref{eq:SpecProj} of $f_{\ubmk}^-(\bmmu), f_{\bmzero}^+(\bmmu)$  in \eqref{eq:v20} (forming a basis of $\mathcal{V}_{{\bmmu},0}$).

Below, we compute the $2\times 2$ matrix representing the action  of the operator ${\sL}_{\bmmu,\epsilon}: \mathcal{V}_{\bmmu,\epsilon} \to \mathcal{V}_{\bmmu,\epsilon}$ on the basis $\cF$.

\begin{lemma} \label{lem:MatLmueps}
Fix $\bmmu \in ( \cR_\ubmk\setminus \cD_\ubmk ) \cap ( \R^2 \setminus \Z^2 ) $.  
The operator ${\sL}_{\bmmu,\epsilon}: \mathcal{V}_{\bmmu,\epsilon} \to \mathcal{V}_{\bmmu,\epsilon}$ on the basis $\mathcal{F}$ in \eqref{eq:basisF} admits the following matrix representation:
\begin{align} \label{eq:MatRepr}
\tL({\bmmu,\e}) &= \tI({\bmmu, \e})^{-1} \tA({\bmmu, \e}),
\end{align}
where
\begin{align}
&  \tI({\bmmu, \e}):= 
\begin{pmatrix}
    \langle \fI(\bmmu,\epsilon) f_{\ubmk}^- , f_{\ubmk}^-\rangle & 
    \langle \fI(\bmmu,\epsilon) f_{\bmzero}^+, f_{\ubmk}^- \rangle \\
    \langle \fI(\bmmu,\epsilon) f_{\ubmk}^- , f_{\bmzero}^+ \rangle & 
    \langle \fI(\bmmu,\epsilon) f_{\bmzero}^+ , f_{\bmzero}^+ \rangle
\end{pmatrix}  = 
\begin{pmatrix}
    \iota_{11}(\bmmu,\epsilon)  & 
   \iota_{01}(\bmmu,\epsilon)\\
    \iota_{01}(\bmmu, \epsilon) & 
    \iota_{00}(\bmmu, \epsilon)
\end{pmatrix}  
\ , \label{eq:matrI}\\
&  \tA({\bmmu, \e}):=
\begin{pmatrix}
    \langle \fL({\bmmu,\epsilon})f_{\ubmk}^- ,  f_{\ubmk}^- \rangle  & 
    \langle \fL({\bmmu,\epsilon}) f_{\bmzero}^+ ,  f_{\ubmk}^- \rangle \\
    \langle \fL({\bmmu,\epsilon}) f_{\ubmk}^- ,  f_{\bmzero}^+ \rangle &
    \langle \fL({\bmmu,\epsilon}) f_{\bmzero}^+ ,  f_{\bmzero}^+ \rangle 
\end{pmatrix}
= 
\im 
\begin{pmatrix}
     a_{11}(\bmmu,\epsilon)  & 
   a_{01}(\bmmu,\epsilon)\\
    a_{10}(\bmmu, \epsilon) & 
    a_{00}(\bmmu, \epsilon)
\end{pmatrix}
\label{eq:matrA}
\end{align}
and $f_{\ubmk}^-\equiv f_{\ubmk}^-(\bmmu), f_{\bmzero}^+ \equiv f_{\bmzero}^+(\bmmu)$ in \eqref{eq:v20}, and 
\begin{equation}\label{fIfL}
    \fI(\bmmu, \e) := P_{\bmmu, 0}^* \, P_{\bmmu, \e}^* \, P_{\bmmu, \e} \, P_{\bmmu, 0}\ , \qquad 
    \fL(\bmmu, \e) := P_{\bmmu, 0}^* \, P_{\bmmu, \e}^* \,  \sL_{\bmmu, \e} \, P_{\bmmu, \e} P_{\bmmu, 0}\ .
\end{equation}
The functions $\e\to \iota_{ij}(\bmmu, \e)$ and $\e \to a_{ij}(\bmmu, \e)$ are real valued and  analytic. Furthermore
\begin{equation}
    \tI(\bmmu, 0) = \begin{pmatrix}
        \iota_{11}(\bmmu,0) & 0 \\
        0 &  \iota_{00}(\bmmu,0)
    \end{pmatrix} \ , \qquad 
    \tA(\bmmu, 0) = 
    \begin{pmatrix}
   \lambda_{\underline{\bmk}}^-(\bmmu) \, \iota_{11}(\bmmu,0)   & 0\\
    0 &   \lambda_{\bmzero}^+(\bmmu) \, \iota_{00}(\bmmu,0) 
\end{pmatrix}
,
\end{equation}
and 
\begin{equation}\label{n11n00}
\begin{aligned}
\iota_{11}(\bmmu,0)  &= \langle f_{\ubmk}^{-} , f_{\ubmk}^{-} \rangle \, = \, \frac{1}{2}  \left( N^2 + |\ubmk+\ubmmu|^2 \right) , \quad 
\iota_{00}(\bmmu,0)  &= \langle f_{\bmzero}^{+} , f_{\bmzero}^{+} \rangle \, = \, \frac{1}{2} \left( N^2 + |\ubmmu|^2 \right) .
\end{aligned}
\end{equation}

\end{lemma}

\begin{proof}
We write 
\begin{align*}
\mathcal{L}_{\bmmu,\epsilon} f_{\ubmk}^{-}(\bmmu,\e) = \alpha f_{\ubmk}^{-}(\bmmu,\e) + \beta f_{\bmzero}^{+}(\bmmu,\e), &\;\; \; \;  \mathcal{L}_{\bmmu,\epsilon} f_{\bmzero}^{+}(\bmmu,\e) = \gamma f_{\ubmk}^{-}(\bmmu,\e) + \delta f_{\bmzero}^{+}(\bmmu,\e) ,
\end{align*}
for some $\alpha, \beta, \gamma, \delta \in \mathbb{C}$. Taking the scalar product with $f_{\ubmk}^{-}(\bmmu,\e)$, $f_{0}^{+}(\bmmu,\e)$ and using \eqref{eq:basisF} give \eqref{eq:MatRepr} \eqref{eq:matrI} and \eqref{eq:matrA}.

The reality of the entries of \( \tI(\bmmu, \e) \) and the purely imaginary nature of \( \tA(\bmmu, \e) \) follow from the reversibility property \eqref{base.rev}, the reversibility of the projector \( P_{\bmmu, \e} \) (Lemma \ref{lem:Kato}, point 2), the complex reversibility of \( \sL_{\bmmu, \e} \) (see \eqref{Lmue.rev}), and the identity  
\begin{equation}
    ( f, g ) = \overline{( \bar \varrho f, \bar \varrho g)} \ , \quad
    \forall f, g \in L^2(\R^2, \C) \ .
\end{equation}   
\end{proof}

In the following, we provide an asymptotic expansion of the entries of $\tL(\bmmu, \e)$.
\begin{proposition}\label{prop.exp}
Fix $\bmmu \in ( \cR_\ubmk\setminus \cD_\ubmk ) \cap ( \R^2 \setminus \Z^2 ) $.   The real analytic entries of the matrix $\tL(\bmmu, \e)$ in \eqref{eq:MatRepr} admit the following expansions as $\e \to 0$:
\begin{align}
\label{iota11}
    \iota_{11}(\bmmu, \e) & = \iota_{11}(\bmmu, 0) + \mathcal{O}(\e^2) , \quad 
    \iota_{10}(\bmmu, \e)   = \gamma_1(\bmmu) \e + \mathcal{O}(\e^2)  , \quad 
    \iota_{00}(\bmmu, \e) = \iota_{00}(\bmmu,0) + \mathcal{O}(\e^2) , \\
       \label{a00}
    a_{11}(\bmmu, \e) & =  \underline{w} \, \iota_{11}(\bmmu, 0) +  \mathcal{O}(\e^2)  \ , \quad 
    a_{00}(\bmmu, \e)  = \underline{w} \, \iota_{00}(\bmmu, 0) +  \mathcal{O}(\e^2) , \   \\
    \label{a10}
    a_{10}(\bmmu, \e) & = \beta_1(\bmmu) \e + \mathcal{O}(\e^2)   \ , \quad 
    a_{01}(\bmmu, \e)  = \beta_0(\bmmu) \e + \mathcal{O}(\e^2) \ .
\end{align}
where $\iota_{11}(\bmmu,0)$ and $\iota_{00}(\bmmu,0)$ are given in \eqref{n11n00}, 
\begin{equation}\label{def:undw}
 \uw:= \uw(\bmmu):= -\im \lambda^+ =   \underline{\bmc}\cdot \bmmu  + \Omega(\bmmu)   
  = \underline{\bmc}\cdot (\ubmk + \bmmu) - \Omega(\ubmk + \bmmu)
\end{equation}
with $\underline{\bmc}$ in \eqref{eq:trav}, 
and finally the functions $\beta_1(\bmmu), \beta_0(\bmmu) $ and $\gamma_1(\bmmu)$ fulfill
\begin{align}\label{b1-g1}
\beta_1(\bmmu) - \gamma_1(\bmmu) \uw(\bmmu) &  = 
\frac{1}{8} (\ubmk^{\perp} \cdot {\bmmu})  \frac{ \big( |\ubmk+{\bmmu}| -|\ubmk|\big) \, (|\ubmk + {\bmmu}| + |\ubmk| - |{\bmmu}| \big) \, \big(N^2 + |{\bmmu}|^2\big) }{N\,|\ubmk|\, |{\bmmu}|\,  |\ubmk+{\bmmu}|} \\
\label{b0-g0}
\beta_0(\bmmu) - \gamma_1(\bmmu) \uw(\bmmu) & = 
- \frac{1}{8} (\ubmk^{\perp} \cdot {\bmmu}) \frac{\big( |\ubmk| +|{\bmmu}|\big) \, \big( |\ubmk + {\bmmu}| + |\ubmk| - |{\bmmu}|\big) \, \big(|\ubmk + {\bmmu}|^2 + N^2 \big)  }{N \, |\ubmk| \, |{\bmmu}| \, |\ubmk + {\bmmu}|}
\end{align}
\normalsize

\end{proposition}
 We stress  that the quantities $\beta_0(\bmmu)$, $\beta_1(\bmmu)$ and $\gamma_1(\bmmu)$ do not depend on $\e$. The proof of Proposition \ref{prop.exp} relies on machinery introduced in the next section, with the proof itself deferred to Section \ref{sec:proof.prop}. Meanwhile, assuming Proposition \ref{prop.exp}, we can prove Theorem \ref{thm:expansion}.

\subsection{Proof of Theorem \ref{thm:expansion}}\label{proof}
By \eqref{eq:MatRepr} and the expansions of Proposition \ref{prop.exp}, 
\begin{equation}
    \tL(\bmmu, \e) =  
       \begin{pmatrix}
       \lambda^+ + \im r_1(\e^2) &  \im \dfrac{\beta_0(\bmmu) - \gamma_1(\bmmu)\und{w}(\bmmu)}{\iota_{11}(\bmmu,0)} \e + \im r_2(\e^2) \\
      \im \dfrac{\beta_1(\bmmu) - \gamma_1(\bmmu)\und{w}(\bmmu)}{\iota_{00}(\bmmu,0)}\e + \im r_3(\e^2) &  \lambda^+ + \im r_4(\e^2)
   \end{pmatrix}
\end{equation}
Next,  
$\tb_1(\bmmu) := \frac{\beta_1(\bmmu) - \gamma_1(\bmmu)\und{w}(\bmmu)}{\iota_{00}(\bmmu,0)}$ is obtained by \eqref{b1-g1}  and \eqref{n11n00}. Similarly $\tb_0(\bmmu) := \frac{\beta_0(\bmmu) - \gamma_1(\bmmu)\und{w}(\bmmu)}{\iota_{11}(\bmmu,0)}$ is computed via  \eqref{b0-g0}  and \eqref{n11n00}.
The eigenvalues \eqref{lambdapm}
are readily computed. 

The asymptotics in \eqref{eto0}-\eqref{etoinf} are obtained by computing  \( \te(\bmmu) \) along the corresponding parametrizations, using the asymptotics of the functions $\varphi_\pm$ in Lemma \ref{lem:Rk}.

\section{Taylor expansion of the operators and entanglement coefficients } \label{sec:Taylor}

This section introduces the technical machinery to prove Proposition \ref{prop.exp}, providing the entanglement coefficients \eqref{entcoeff}.
To begin with, we introduce the {\em jets} of  the operator 
$\sL_{\bmmu, \e}$ in \eqref{L0}, obtained by Taylor-expanding the operator $\sL_{\bmmu, \e}$ at the point $(\bmmu, 0)$:
\begin{equation}
\sL_{j} :=  \frac{1}{ j!} \big(\pa^j_\e \sL_{\bmmu,\e} \big)(\bmmu,0)\, 
\e^j \,  \ , 
\qquad j \in \N_0 \ , 
\end{equation}
so that 
\begin{equation} 
\sL_{ \bmmu,\e } = \sL_0 + \sL_1 +  \cO_2 \, .
\end{equation}
Explicitly, by \eqref{L0},  \eqref{eq:Lmueps} and  \eqref{eq:trav} we have that $\sL_0$ is the matrix of Fourier multipliers
\begin{align}\label{cL0}
    & \sL_0 =     
    \begin{pmatrix}
\underline{\bmc} \cdot  (\nabla+ \im \bmmu) & -N^2 \,  (\partial_x + \im \mu_{1}) \, \left| \nabla + \im \bmmu \right|^{-2} \\
- (\partial_x+ \im \mu_1) & \underline{\bmc} \cdot (\nabla + \im \bmmu)
\end{pmatrix}
\end{align}
whereas $\sL_1$ is the operator
\small
\begin{align}\label{cL1}
 & \sL_1   =
\e \cos(\ubmk \cdot \bmx) 
\begin{pmatrix}
- (N \, |\ubmk|)^{-1}\,   \ubmk^\perp  \cdot (\nabla + \im \bmmu)   &   -   \ubmk \cdot (\gradp + \im \bmmu^\perp) \, \vert {\bf D + \bmmu}\vert^{-2} \\
0  & 
- (N \, |\ubmk|)^{-1}\, \ubmk^\perp \cdot ( \nabla +\im \bmmu)
 - N^{-1} \, |\ubmk|\,  \ubmk \cdot (\gradp+ \im \bmmu^\perp)  \vert {\bf D} + \bmmu  \vert^{-2}
\end{pmatrix}  \ . 
\end{align}
\normalsize
Analogously we Taylor-expand the projectors $P_{\bmmu,\e}$ in  \eqref{eq:SpecProj} at the point $(\bmmu, 0)$ as
\begin{equation*}
P _{ \bmmu,\e } =P_0 + P_1 +  \cO_2 \, 
\end{equation*}
where  
\begin{equation}\label{Psani}
\begin{aligned}
 &P_0 := P_{ \bmmu,0 }= -\frac{1 }{2\pi \im } \oint_\Gamma (\sL_{ \bmmu,0}-\lambda)^{-1} {\rm d}\lambda  \ , 
 \end{aligned} 
 \end{equation}
and  
\begin{equation}\label{hP}
 \begin{aligned}
  P_1 := \mathcal{P} \big[\sL_1\big] \ , \quad 
 \mathcal{P} & \big[A\big] :=
  \frac{1 }{2\pi \im } \oint_\Gamma (\sL_{ \bmmu,0}-\lambda)^{-1}  \, A \, (\sL_{ \bmmu,0}-\lambda)^{-1} {\rm d}\lambda \, ,
 \end{aligned}
\end{equation}
where $\Gamma$ is the curve of Lemma \ref{lem:Kato}.

We therefore obtain the following lemma about the jets of the  operators $\mathfrak{L} ({\bmmu,\e}) $ and $\mathfrak{I}(\bmmu,\e)$:
\begin{lemma} \label{expansionthm}
 The operators $\mathfrak{L} (\bmmu,\e) $  and $\fI(\bmmu, \e)$ in \eqref{fIfL} have the Taylor expansions 
 \begin{equation}
 \begin{aligned}
& \mathfrak{L} ({\bmmu,\e}) =
\mathfrak{L}_0 +
\mathfrak{L}_1 +\cO_2 \ , \qquad
 \mathfrak{I} ({\bmmu,\e}) =
\mathfrak{I}_0 +
\mathfrak{I}_1 +\cO_2
\end{aligned}
\end{equation}
where
 \begin{subequations}\label{upto4exp} 
\begin{align}\label{ordini012}
\mathfrak{L}_0 & := P_0^*{\sL}_0 P_0 \, , 
 \quad 
 \mathfrak{L}_1 := P_0^*\left( P_1^* \sL_0 + \sL_0 P_1  + {\sL}_1 \right) P_0 \, , \\
\mathfrak{I}_0 & := P_0^*  P_0 \, , 
 \quad 
 \mathfrak{I}_1 := P_0^*(P_1 + P_1^*) P_0 \, ,  
\end{align}
\end{subequations}
with ${\sL}_j$, $j=0,1$, in \eqref{cL0}, \eqref{cL1} and $P_j$, $j=0,1$, in \eqref{Psani}.
\end{lemma}
\begin{proof}
This readily follows Taylor-expanding the explicit formula \eqref{fIfL}, using the fact that \( P_j^2 = P_j \).
\end{proof}

\subsection{Entanglement coefficients.} Next we introduce the  {\em entanglement} coefficients, that describe the action of the jets $\sL_\ell$ and $\cJ(\bmmu) \sL_{\ell}$ over the vectors
$\{f_{\bmk}^\sigma\}_{\bmk \in \Z^2, \sigma = \pm}$, $f_{\bmk}^\sigma:= f_{\bmk}^\sigma(\bmmu)$ in \eqref{eq:v20}, forming a basis of $L^2(\T^2, \C^2)$ and fulfilling the symplectic relation  \eqref{symplfjs}.

First, given  a linear 
operator $A  = \begin{pmatrix}
    A_1 & A_2 \\
    A_3 & A_4
\end{pmatrix}  \colon   X^1_{\ubmk} \to X^0_{\ubmk} $ (recall definition \eqref{Hlk})   
we define for any $j_1, j_2 \in \Z$ its {\em matrix elements} 
\begin{equation}\label{def.matrix.elements}
A^{j_1}_{j_2} := 
\begin{pmatrix}
        ( A_1 [e^{\im j_1 \ubmk\cdot \bmx}] , \ e^{\im j_2 \ubmk\cdot \bmx  } )_{L^2(\T^2)} & ( A_2[e^{\im j_1 \ubmk\cdot \bmx}] , \ e^{\im j_2 \ubmk\cdot \bmx  } )_{L^2(\T^2)}\\
        ( A_3[e^{\im j_1 \ubmk\cdot \bmx}] , \ e^{\im j_2 \ubmk\cdot \bmx  } )_{L^2(\T^2)} & ( A_4[e^{\im j_1 \ubmk\cdot \bmx}] , \ e^{\im j_2 \ubmk\cdot \bmx  } )_{L^2(\T^2)}
\end{pmatrix} \ , 
\end{equation}
where for $f, g \in L^2(\T^2, \C)$ we denoted 
$$ (f,g) = \frac{1}{(2\pi)^2} \int_{\T^2} f(\bmx) \, \overline{g(\bmx)} \, d \bmx \  . $$
In this way  the action of $A\colon  X^1_{\ubmk} \to X^0_{\ubmk} $  is given by  
$$
 X^1_{\ubmk}  \ni \vec{h}(\bmx) = \sum_{j_1 \in \Z} \vec{h}_{j_1} e^{\im j_1 \ubmk\cdot \bmx} \mapsto  
(A \vec{h})(\bmx) = \sum_{j_2\in \Z} 
\Big( \sum_{j_1 \in \Z} A^{j_1}_{j_2} \, \vec{h}_{j_1} \Big) \, e^{\im j_2 \ubmk\cdot \bmx}  \in  X^0_{\ubmk} \, .
$$
Next, given $\bmkappa \in \ubmk\Z$, we define the projected operator $A^{[\bmkappa]}$ 
whose matrix coefficients are those of $A$ supported on the 
{\em ``band"} $ (j_2 -j_1) \ubmk = \bmkappa $, i.e.
\begin{equation}\label{band.def}
A^{[\bmkappa]} \vec{h} = \sum_{j_1, j_2 \in \Z \atop (j_2-j_1)\ubmk = \bmkappa } A^{j_1}_{ j_2} \, \vec{h}_{j_1} \, e^{\im j_2 \ubmk \cdot \bmx} \, . 
\end{equation}
The following identities, for operators $X^\infty_\ubmk \to X^\infty_\ubmk$, where $X^\infty_\ubmk:= \cap_{\ell} X^\ell_\ubmk$, are easy to check: 
\begin{equation}\label{A.exp.band}
    A = \sum_{\bmkappa \in \ubmk\Z} A^{[\bmkappa ]} \ , 
    \qquad
    \big[A^{[\bmkappa]}\big]^* = (A^*)^{[-\bmkappa]} \ ,  
    \qquad
      (A\circ B)^{[\bmkappa]} = \sum_{
      \bmkappa_1,  \bmkappa_2 \in \ubmk\Z \atop 
      \bmkappa_1 + \bmkappa_2 = \bmkappa} A^{[\bmkappa_1]} \circ B^{[\bmkappa_2]} \ . 
\end{equation}

Following \cite{berti2024infinitely}, we introduce the spaces $\mathfrak{F}_\ell$, which highlight important structural properties of the jets of the operator $\mathcal{L}_{\bmmu, \e}$. Roughly speaking, these are spaces of operators whose jets $A_j$ have ``bands'' of order at most $j$ and share the same parity as $j$.

\begin{definition}[{\bf Spaces $\mathfrak{F}_\ell$}]
\label{defFell}
For any $\ell \in \mathbb{N}_0$, we define $\mathfrak{F}_\ell$ as the space of operators $X^\infty_{\ubmk} \to X^\infty_{\ubmk}$ satisfying
\begin{equation} \label{TFj}
A^{[n \ubmk]} = 0 \quad \text{if } |n| > \ell \quad \text{or} \quad n \not\equiv \ell \pmod{2} \, .
\end{equation}
\end{definition}

The following lemma, proved as in \cite[Lemma 5.5]{berti2024infinitely}, states some  properties of the spaces 
$ \mathfrak{F}_\ell$ which will be used repeatedly. 
\begin{lemma}\label{Fj}
Let $ A \in \mathfrak{F}_{\ell}$, $ B \in \mathfrak{F}_{\ell'}$.  Then 
\begin{itemize}
    \item[(i)] {\bf Composition:} $A \circ B \in \mathfrak{F}_{\ell+\ell'}$ with 
    $    (A\circ B)^{[\pm(\ell+\ell')\ubmk]}= A^{[\pm \ell\ubmk]} \circ B^{[\pm \ell'\ubmk]}$;
    \item[(ii)] {\bf Adjoint:} $A^* \in \mathfrak{F}_{\ell}$.
    \item[(iii)]  {\bf Finite range interactions:} For any $ \vec{v}_1, \vec{v}_2 \in \C^2$,
    \begin{equation}\label{AFlele}
( A \,  \vec{v}_1 e^{\im j_1 \ubmk \cdot \bmx} , 
\vec{v}_2 e^{\im j_2 \ubmk \cdot \bmx } ) = 0 \ \ \mbox{ if } |j_1 - j_2| > \ell  \mbox{ or } j_1 - j_2 \not\equiv \ell \mbox{ mod } 2 \ . 
\end{equation}
\end{itemize}
\end{lemma}
\begin{proof}
(i) It follows from the third of  \eqref{A.exp.band}, using that if $\bmkappa = n \ubmk$, $\bmkappa_1 = j_1 \ubmk$, $\bmkappa_2 = j_2 \ubmk$, then necessarily $j_1 \equiv \ell \mod 2$, $j_2 \equiv \ell' \mod 2$, hence $j_1 + j_2 \equiv \ell + \ell' \mod 2$, and if  $\ell + \ell' < |n| \leq |j_1| + |j_2|$, then either $|j_1| > \ell$
hence $A^{[\bmkappa_1]} = 0$, 
or $|j_2| > \ell'$ hence  $ B^{[\bmkappa_2]} = 0$. 
Moreover, if $ |n| = \ell + \ell'$, being $j_1 + j_2 = n$ and $|j_1| \leq \ell$, $|j_2| \leq \ell'$, 
the only possibility is that $j_1 = \ell$, $j_2 = \ell'$.\\
(ii) It follows from the second of \eqref{A.exp.band}. \\
(iii) By the first of \eqref{A.exp.band} and the definition \eqref{band.def}, we have
$$
( A \,  \vec{v}_1 e^{\im j_1 \ubmk \cdot \bmx} , 
\vec{v}_2 e^{\im j_2 \ubmk \cdot \bmx } ) = \sum_{\bmkappa \in \ubmk\Z} ( A^{[\bmkappa]} \,  \vec{v}_1 e^{\im j_1 \ubmk \cdot \bmx} , 
\vec{v}_2 e^{\im j_2 \ubmk \cdot \bmx } ) = ( A^{[(j_2 - j_1)\ubmk]} \,  \vec{v}_1 e^{\im j_1 \ubmk \cdot \bmx} , 
\vec{v}_2 e^{\im j_2 \ubmk \cdot \bmx } )
$$
and the claim follows from  property \eqref{TFj}.
\end{proof}

In the next lemma, we state useful properties of the jets of the operators $\sL_{\bmmu,\e}$, $P_{\bmmu,\e}$, $\mathfrak{L}(\bmmu, \e)$, $\mathfrak{J}(\bmmu, \e)$.
\begin{lemma}\label{jetsinFl}
The jets of the 
operators $\sL_{\bmmu, \e}$ in \eqref{L0},   $P_{\bmmu, \e}$ in \eqref{eq:SpecProj}, $\mathfrak{L}(\bmmu, \e)$  and 
 $\mathfrak{J}(\bmmu, \e)$ in \eqref{fIfL}
fulfill: 
\begin{align}
\sL_\ell \ , \ \ P_\ell \ , \ \ \mathfrak{L}_\ell \ , \ \ \mathfrak{I}_\ell \in \mathfrak{F}_\ell \ , \qquad \ell  = 0, 1 \ . 
\end{align}
In addition $\sL_0^{[\bmzero]} = \sL_0$, the $\pm\ubmk$-th band of the operator $\sL_1$ in \eqref{cL1} is 
\begin{align}\label{Lk1}
    \sL_{1}^{[\pm\ubmk]} &= \frac{\e}{2} e^{\pm\im \ubmk \cdot \bmx} \,
    \begin{pmatrix}
        - (N \, |\ubmk|)^{-1}\,   \im \ubmk^\perp  \cdot \bmmu   &   \im   \ubmk^{\perp} \cdot \bmmu \, \vert \bmD + \bmmu \vert^{-2} \\
0  & 
- (N \, |\ubmk|)^{-1}\, \im \ubmk^\perp \cdot \bmmu
 + N^{-1} \, |\ubmk|\, \im \ubmk^{\perp} \cdot \bmmu  \vert \bmD + \bmmu  \vert^{-2}
    \end{pmatrix}  \ , 
\end{align}
and, recalling \eqref{hP}, 
\begin{align}\label{proj.Fl}
    P_1^{[\pm \ubmk]} = \cP[\sL_1^{[\pm\ubmk]}] \ . 
\end{align}
\end{lemma}
\begin{proof}
The operator $\sL_0$ in \eqref{cL0} is a Fourier multiplier, so that its matrix elements $(\sL_0)^{q_1}_{q_2} = 0$ for $q_1 \neq q_2$, hence the only non-trivial band is the zeroth.

We now compute the matrix elements $(\sL_{1})_{q_2}^{q_1}$
   with $q_1 , q_2 \in \Z$. 
Using the explicit expression of $\sL_1$ in \eqref{cL1}, we get
\begin{equation}\label{Bl.el}
    (\sL_{1})_{q_2}^{q_1}
    = 
    \frac{\e}{2} \ \delta(q_2 - q_1, \pm 1)\ 
  \begin{pmatrix}
        - (N \, |\ubmk|)^{-1}\,   \im \ubmk^\perp  \cdot \bmmu   &   \im   \ubmk^{\perp} \cdot \bmmu \, \vert  q_1\ubmk + \bmmu \vert^{-2} \\
0  & 
- (N \, |\ubmk|)^{-1}\, \im \ubmk^\perp \cdot \bmmu
 + N^{-1} \, |\ubmk|\, \im \ubmk^{\perp} \cdot \bmmu  \, \vert q_1\ubmk + \bmmu  \vert^{-2}
    \end{pmatrix}  \ . 
\end{equation}
Then, for $q_2 - q_1  \not \in \{\pm 1\} $, then
 $(\sL_{1})_{q_2}^{q_1} = 0$, proving that $\sL_1 \in \mathfrak{F}_1$. 
 In addition, we compute the bands $\pm 1$: 
writing  $\vec{h} = \sum_q \vec{h}_q e^{\im q \ubmk \cdot \bmx} \in X^1_{\ubmk}$, we have
\begin{align*}
\sL_{1}^{[\ubmk]}\vec{h} & = \sum_{q_2 - q_1 = 1 }    (\sL_{1})_{q_2}^{q_1} \vec{h}_{q_1} \, e^{\im q_2 \ubmk \cdot \bmx} 
 =  \sum_{q_1 }    (\sL_{1})_{q_1+1}^{q_1} \, \vec{h}_{q_1} \, e^{\im (q_1+1)\ubmk\cdot \bmx} \\
&  \stackrel{\eqref{Bl.el}}{ = }
 \frac{\e}{2}\, e^{ \im \ubmk \cdot \bm x} 
 \begin{pmatrix}
        - (N \, |\ubmk|)^{-1}\,   \im \ubmk^\perp  \cdot \bmmu   &   \im   \ubmk^{\perp} \cdot \bmmu \, \vert \bmD + \bmmu \vert^{-2} \\
0  & 
- (N \, |\ubmk|)^{-1}\, \im \ubmk^\perp \cdot \bmmu
 + N^{-1} \, |\ubmk|\, \im \ubmk^{\perp} \cdot \bmmu  \vert \bmD + \bmmu  \vert^{-2}
    \end{pmatrix} \vec{h},
\end{align*}
proving \eqref{Lk1}. For $\ell = -1$ the proof is analogous.

The operator $(\sL_{\bmmu,0}- \lambda)^{-1} \in \mathfrak{F}_0$, being a Fourier multiplier,  as well as $P_0$ in \eqref{Psani}. 
The operator $P_1$ in \eqref{hP} belongs to $\mathfrak{F}_1$ by the composition property in Lemma \ref{Fj}, which also gives \eqref{proj.Fl}. 
Finally, $\mathfrak{L}_\ell$ and $\mathfrak{I}_\ell$ belongs to $\mathfrak{F}_\ell$, $\ell = 0,1$, in view of their expressions in \eqref{upto4exp} and Lemma \ref{Fj}.
\end{proof}
We are now ready to introduce the entanglement coefficients.
For a fixed $\bmmu \in ( \mathcal{R}_{\ubmk} \setminus \cD_{\ubmk} ) \cap ( \R^2 \setminus \Z^2 )$, consider $\{f_{\bmk}^\sigma\}_{\bmk \in \Z^2, \sigma = \pm}$, $f_{\bmk}^\sigma:= f_{\bmk}^\sigma(\bmmu)$ in \eqref{eq:v20}, forming a basis of $H^1(\T^2) \times L^2(\T^2)$ such that \eqref{symplfjs} holds.

\begin{definition}
Let $\bmmu \in ( \mathcal{R}_{\ubmk} \setminus \cD_{\ubmk} ) \cap ( \R^2 \setminus \Z^2 )$,
and consider  $f_{\bmk}^\sigma:= f_{\bmk}^\sigma(\bmmu)$ in \eqref{eq:v20} with  $\bmk \in \Z^2$.
The {\em entanglement coefficients} are:
\begin{equation}\label{entcoeff}
\begin{aligned}
    & \ent{\ell}{\bmkappa }{\bmj'}{\bmj}{\sigma'}{\sigma}
:= \big (\cJ(\bmmu) {\sL}_{\ell}^{[\bmkappa ]} f_{\bmj}^\sigma, f_{\bmj'}^{\sigma'}  \big)   \ , \qquad \bmkappa \in \{ n \ubmk, \ n\in  \Z\} , \\
& \entL{\ell}{\bmkappa }{\bmj'}{\bmj}{\sigma'}{\sigma}
:= \big ( {\sL}_{\ell}^{[\bmkappa]} f_{\bmj}^\sigma, f_{\bmj'}^{\sigma'}  \big) , \qquad \qquad \;\;\, \bmkappa \in \{ n \ubmk, \ n\in  \Z\} ,
\end{aligned}
\end{equation}
where $\cJ(\bmmu)$ is as in \eqref{Jmu}, and ${\sL}_{\ell}^{[\bmkappa]}$ is the $\kappa$-band of the jet $\sL_\ell$.
\end{definition}

Since $\cJ(\bmmu)$,  $\sL_{\bmmu, \e}$ and the basis are reversible, see \eqref{Lmue.rev}, \eqref{base.rev}, 
\begin{equation}
     \overline{\ent{\ell}{\bmkappa }{\bmj'}{\bmj}{\sigma'}{\sigma}} =  \ent{\ell}{\bmkappa }{\bmj'}{\bmj}{\sigma'}{\sigma} \ , \quad 
     \overline{\entL{\ell}{\bmkappa }{\bmj'}{\bmj}{\sigma'}{\sigma}} = - \entL{\ell}{\bmkappa }{\bmj'}{\bmj}{\sigma'}{\sigma}.
\end{equation}

\noindent{\bf Computation of the entanglement coefficients.}
We now compute the  entanglement coefficients that will appear in Section \ref{sec:proof.prop}.
\begin{lemma}\label{lem:explicit-ent-coeff}
    For $\sigma = \pm$, $w^\sigma_{n \ubmk}, n=0,1$ as in \eqref{eq:lambdasigmak-restr}, $f^\sigma_{n\ubmk}, n=0,1$ in \eqref{eq:v20}, we have
    \begin{align}
        & \ent{1}{-\ubmk}{\bmzero}{\ubmk}{-}{-} = \frac{\e}{4} \, (\ubmk^{\perp} \cdot \bmmu) \, \mu_1 \, \frac{ ( |\ubmk+\bmmu| - |\bmk| ) ( |\ubmk+\bmmu| + |\bmk| + |\bmmu| ) }{ |\ubmk| \, |\ubmk+\bmmu| } , \label{eq:B1mk} \\
        & \ent{1}{\ubmk}{\ubmk}{\bmzero}{+}{+} = \frac{\e}{4} \, (\ubmk^{\perp} \cdot \bmmu) (\tm+\mu_1) \, \frac{ ( |\bmmu|+ |\ubmk| )( |\ubmk| - |\bmmu| - |\ubmk + \bmmu| )  }{ |\ubmk| |\bmmu| } , \label{eq:B1pk} \\
        & \entL{0}{\bmzero}{\bmzero}{\bmzero}{+}{-} = \im 
        w_{\bmzero}^-(\bmmu) \, (f_{\bmzero}^-, f_{\bmzero}^+)
        , \label{eq:L000p0m} \\
        &  \entL{0}{\bmzero}{\bmzero}{\bmzero}{-}{+} = \im 
        w_{\bmzero}^+(\bmmu) \, (f_{\bmzero}^-, f_{\bmzero}^+), \label{eq:L000m0p} \\
        &  \entL{0}{\bmzero}{\ubmk}{\ubmk}{+}{-}  = 
        \im w_{\ubmk}^-(\bmmu) \, (f_{\ubmk}^-, f_{\ubmk}^+)
        , \label{eq:L00kpkm} \\
        &  \entL{0}{\bmzero}{\ubmk}{\ubmk}{-}{+}  =   \im w_{\ubmk}^+(\bmmu) \, (f_{\ubmk}^-, f_{\ubmk}^+)
        , \label{eq:L00kmkp} \\
        & \entL{1}{\ubmk}{\ubmk}{\bmzero}{-}{+} = - \frac{\im}{4} \e \, (\ubmk^{\perp} \cdot \bmmu) \frac{ (|\ubmk| + |\bmmu|) \, \left[ (|\ubmk|-|\bmmu|) |\ubmk+\bmmu| + N^2 \right] }{ |\ubmk| |\bmmu| N  }  , \label{eq:L1kk-0+} \\
        & \entL{1}{-\ubmk}{\bmzero}{\ubmk}{+}{-} = \frac{\im}{4} \e \, (\ubmk^{\perp} \cdot \bmmu) \frac{ (|\ubmk+\bmmu| - |\ubmk|) \, \left[ (|\ubmk+\bmmu|+|\ubmk|) |\bmmu| - N^2 \right] }{ |\ubmk| |\ubmk + \bmmu| N  }  . \label{eq:L1mk0+km}
    \end{align}
\end{lemma}
\begin{proof}

We report the explicit computations for $\ent{1}{-\ubmk}{\bmzero}{\ubmk}{-}{-}$ and for $\entL{0}{\bmzero}{\bmzero}{\bmzero}{+}{-}$; 
the others one are obtained in a similar way, and we omit the computation. A code for computing all entanglement coefficients using Mathematica can be found at the link in footnote \footnote{https://git-scm.sissa.it/spasqual/2D-boussinesq-system/}.

By \eqref{entcoeff} we have
\begin{align*}
    \ent{1}{-\ubmk}{\bmzero}{\ubmk}{-}{-} &= ( \cJ(\bmmu) \sL_{1}^{[-\ubmk]} f_{\ubmk}^{-} , f_{\bmzero}^{-} ) ,
\end{align*}
so by \eqref{Jmu},    \eqref{Lk1}  and \eqref{eq:v20} we get
\begin{align*}
    \mathcal{J}(\bmmu) \sL_{1}^{[-\ubmk]} f_{\ubmk}^{-} &= - \e \frac{\mu_1}{2\sqrt{2}} (\ubmk^{\perp} \cdot \bmmu) |\ubmk|^{-1} |\ubmk+\bmmu|^{-1}
    \begin{pmatrix}
        N^{-1} ( |\ubmk|^2 - |\ubmk+\bmmu|^2 ) \\ |\ubmk| - |\ubmk+\bmmu|
    \end{pmatrix}
    ,
\end{align*}
hence
\begin{align*}
    \ent{1}{-\ubmk}{\bmzero}{\ubmk}{-}{-} &= \frac{\e}{4} \mu_1 (\ubmk^{\perp} \cdot \bmmu) |\ubmk|^{-1} |\ubmk+\bmmu|^{-1} \left[ |\ubmk+\bmmu|^2 - |\ubmk|^2 + |\bmmu|( |\ubmk+\bmmu| - |\ubmk|  ) \right] ,
\end{align*}
from which we can deduce \eqref{eq:B1mk}. 

Similarly, by \eqref{entcoeff} and \eqref{eq:v20} we have
\begin{align*}
    \entL{0}{\bmzero}{\bmzero}{\bmzero}{+}{-} &= ( \sL_{0}^{[\bmzero]} f_{\bmzero}^{-} , f_{\bmzero}^{+} )  = 
    \im w_{\bmzero}^-(\bmmu) ( f_{\bmzero}^{-} , f_{\bmzero}^{+} ) 
\end{align*}
hence 
formula \eqref{eq:L000p0m}.

\end{proof}

The next lemma states key properties of the entanglement coefficients, highlighting their role in computing the action of jets on the unperturbed basis \( \{f_{\bmk}^\sigma\} \).

\begin{lemma}
\label{actionofL}
Let  $\ent{\ell}{\bmkappa }{\bmj'}{\bmj}{\sigma'}{\sigma}$ be the entanglement coefficients in \eqref{entcoeff}. Then
\begin{itemize}
\item[(i)] for any $\ell\in \N_0$ and $\bmkappa\in \{ n \ubmk: n \in \mathbb{Z} \}$,
\begin{equation}\label{Lsacts}
\begin{aligned}
 \sL_{\ell}^{[\bmkappa]} f_{\bmj}^\sigma & =
 \sum_{\sigma_1 = \pm} \im\,  \sigma_1 \, \td_{{\bmj} + \bmkappa}
 \, 
  \ent{\ell}{\bmkappa}{{\bmj}+\bmkappa}{{\bmj}}{\sigma_1}{\sigma}\; f_{{\bmj}+\bmkappa}^{\sigma_1} 
=   \im \, \td_{{\bmj}+\bmkappa}\, \ent{\ell}{\bmkappa}{\bmj+\bmkappa}{\bmj}{+}{\sigma}\,  f_{\bmj+\bmkappa}^+ 
-
\im \, \td_{{\bmj}+\bmkappa} \,  \ent{\ell}{\bmkappa}{\bmj+\bmkappa}{\bmj}{-}{\sigma} \, f_{\bmj+\bmkappa}^- 
\end{aligned}
\end{equation}
where
\begin{align}\label{eq:dk}
\td_{{\bmk}}& := \left(  |\bmk + \bmmu|^2 \, \Omega(\bmk + \bmmu) \right)^{-1};
\end{align}
\item[(ii)] recalling \eqref{hP}, we have
\begin{align}
\label{hppar}
&\mathcal{P}[\sL_{\ell}^{[\bmkappa]}]  f_{\bmj}^\sigma 
= \sum_{\sigma_1 = \pm }\!\!\!\  (-\sigma_1) \, \td_{\bmj 
+\bmkappa} \, \ent{\ell}{\bmkappa}{\bmj+\bmkappa}{\bmj}{\sigma_1}{\sigma}  \, \Res{\bmj, & \bmj+\bmkappa}{\sigma, &\sigma_1} \; f_{\bmj + \bmkappa}^{\sigma_1} \, ,
\end{align}
where 
\begin{equation}\label{generalresidue}
\Res{\bmj, & \bmj_1}{\sigma, &\sigma_1} := -\dfrac{1}{2\pi} \oint_\Gamma  \dfrac{ {\rm d}\lambda }{(\lambda-\im w_{\bmj}^\sigma)(\lambda-\im w_{\bmj_1}^{\sigma_1}) } 
\end{equation}

with  $w_{\bmj}^{\pm} \equiv w_{\bmj}^{\pm}(\bmmu)$ in \eqref{eq:lambdasigmak}. 
Furthermore 
\begin{equation}\label{revres}
\Res{\bmj_0, &\bmj_1}{\sigma_0, &\sigma_1}  = \bar{ \Res{\bmj_0, &\bmj_1}{\sigma_0, &\sigma_1} }\, .
\end{equation}
\item[(iii)] It holds
\begin{subequations}\label{cBactstot}
\begin{align}\label{cBactspar}
&\big( \sL_{\ell_{2}}^{[\bmkappa_{2}]}\mathcal{P}[\sL_{\ell_1}^{[\bmkappa_1]}]  f_{\bmj}^\sigma , f_{\bmj'}^{\sigma'}\big) = \sum_{\sigma_1 = \pm }\! (- \sigma_1)  \, \td_{\bmj_1}  \  \ent{\ell_1}{\bmkappa_1}{\bmj_1}{\bmj}{\sigma_1}{\sigma} \   \entL{\ell_2}{\bmkappa_2}{\bmj'}{\bmj_1}{\sigma'}{\sigma_1} \, 
 \Res{\bmj, &\bmj_1,}{\sigma, &\sigma_1} 
\end{align}
with $\bmj_1 := \bmj+\bmkappa_1$, and
\begin{align}
 \label{cBactsparbis}
&\big( \sL_{\ell_1}^{[\bmkappa_1]}  f_{\bmj}^{\sigma} , \mathcal{P}[\sL_{\ell_2}^{[\bmkappa_2]}] f_{\bmj'}^{\sigma'} \big)  =\!\!  \sum_{\sigma_1 = \pm }   (-\sigma_1) \, 
\td_{\bmxi_1}\ 
{\ent{\ell_{2}}{\bmkappa_{2}}{\bmxi_1}{\bmj'}{\sigma_1}{\sigma'}} \ 
\entL{\ell_{1}}{\bmkappa_{1}}{\bmxi_1}{\bmj}{\sigma_1}{\sigma} \ 
\Res{\bmj', & \bmxi_1}{\sigma', &\sigma_1}  \, 
\end{align}
\end{subequations}
with $\bmxi_1 = \bmj' +\bmkappa_2$.
\end{itemize}
\end{lemma}

\begin{proof} 
(i)
Since the operator $\sL_{\ell}^{[\bmkappa]}$ shifts the harmonic $\bmj$ to $\bmj+\bmkappa$, there exist $\alpha_{\pm}^\sigma \in \mathbb{C}$ such that
\begin{align}\label{eq:decomp}
\sL_{\ell}^{[\bmkappa]} f_{\bmj}^\sigma = \alpha_+^\sigma  f_{\bmj+\bmkappa}^+ + \alpha_-^\sigma  f_{\bmj+\bmkappa}^-.
\end{align} 
Applying $\cJ(\bmmu) $ to both sides, taking the scalar products with $f_{\bmj+\bmk}^\pm$ and using \eqref{symplfjs} yields
\begin{align*}
\langle \cJ(\bmmu) \sL_{\ell}^{[\bmkappa]} f_{\bmj}^\sigma , f_{\bmj+\bmk}^\pm \rangle &= \mp \, \im    \, \alpha_+^\sigma  \, |\bmj+\bmkappa+\bmmu|^2 \, \Omega(\bmj+\bmkappa+\bmmu) 
\end{align*}
so that, using \eqref{eq:dk} and \eqref{entcoeff}, 
\begin{align*}
    \alpha_+^\sigma&= \frac{\im \langle \cJ(\bmmu) \sL_{\ell}^{[\bmkappa]} f_{\bmj}^\sigma , f_{\bmj+\bmk}^+ \rangle }{|\bmj+\bmkappa+\bmmu|^2 \, \Omega(\bmj+\bmkappa+\bmmu)}= \im\,  \td_{{\bmk+\bmj}} \ent{\ell}{\bmkappa }{\bmj+\bmkappa}{\bmj}{+}{\sigma} ,\\
    \alpha_-^\sigma&= \frac{-\im \langle \cJ(\bmmu) \sL_{\ell}^{[\bmkappa]} f_{\bmj}^\sigma , f_{\bmj+\bmk}^- \rangle }{|\bmj+\bmkappa+\bmmu|^2 \, \Omega(\bmj+\bmkappa+\bmmu)}=-\im\,  \td_{{\bmk+\bmj}} \ent{\ell}{\bmkappa }{\bmj+\bmkappa}{\bmj}{-}{\sigma}.
\end{align*}
Together with \eqref{eq:decomp}, this gives \eqref{Lsacts}.

(ii)
By \eqref{eq:v20}, it holds $( \cL_{\bmmu,0} - \lambda )^{-1} f_{\bmj}^{\sigma} = (\im w_{\bmj}^{\sigma} - \lambda )^{-1} f_{\bmj}^{\sigma}$. Hence, by \eqref{hP} 
\begin{align*}
\cP[ \sL_{\ell}^{[\bmkappa]} ] f_{\bmj}^\sigma &= \frac{1}{2\pi \im} \oint_{\Gamma} ( \cL_{\bmmu,0} - \lambda )^{-1} \sL_{\ell}^{[\bmkappa]} ( \cL_{\bmmu,0} - \lambda )^{-1} f_{\bmj}^\sigma \mathrm{d}\lambda \\
&= \frac{1}{2\pi \im} \oint_{\Gamma} ( \cL_{\bmmu,0} - \lambda )^{-1} \sL_{\ell}^{[\bmkappa]} \frac{1}{\im w_{\bmj}^\sigma -\lambda} f_{\bmj}^\sigma \mathrm{d}\lambda \\
&\stackrel{ \eqref{Lsacts} }{=} \sum_{\sigma_1 = \pm} (-\sigma_1) \td_{\bmj+\bmkappa}  \ent{\ell}{\bmkappa }{\bmj+\bmkappa}{\bmj}{\sigma_1}{\sigma}  \Res{\bmj, & \bmj+\bmkappa}{\sigma, &\sigma_1} \, f_{\bmj+\bmk}^{\sigma_1}.
\end{align*}

(iii) By \eqref{hppar}, as $\{f_{\bmj}\}$ is a basis of eigenvectors, 
\begin{align*}
\big( \sL_{\ell_{2}}^{[\bmkappa_{2}]}\mathcal{P}[\sL_{\ell_1}^{[\bmkappa_1]}]  f_{\bmj}^\sigma , f_{\bmj'}^{\sigma'}\big) &= \sum_{\sigma_1 = \pm }\! (- \sigma_1)  \, \td_{\bmj_1}  \, \ent{\ell_1}{\bmkappa_1}{\bmj_1}{\bmj}{\sigma_1}{\sigma} \Res{\bmj, &\bmj_1,}{\sigma, &\sigma_1} \, \left( \sL_{\ell_2}^{[\bmkappa_2]} f_{\bmj+\bmkappa_1}^{\sigma_1} , f_{\bmj'}^{\sigma'} \right) ,
\end{align*}
and we can deduce \eqref{cBactspar} by \eqref{entcoeff}; similarly, we have
\begin{align*}
\big( \sL_{\ell_1}^{[\bmkappa_1]}  f_{\bmj}^{\sigma} , \mathcal{P}[\sL_{\ell_2}^{[\bmkappa_2]}] f_{\bmj'}^{\sigma'} \big)  &=  \sum_{\sigma_1 = \pm }   (-\sigma_1) \,  
\td_{\bmj'+\bmkappa_2}\,
{\ent{\ell_{2}}{\bmkappa_{2}}{\bmj'+\bmkappa_2}{\bmj'}{\sigma_1}{\sigma'}} \, \Res{\bmj', & \bmj'+\bmkappa_2}{\sigma', &\sigma_1} \,
\left( \sL_{\ell_1}^{[\bmkappa_1]} f_{\bmj}^{\sigma} , f_{\bmj'+\bmkappa_2}^{\sigma_1} \right) ,
\end{align*}
and we can deduce \eqref{cBactsparbis} by \eqref{entcoeff}.
\end{proof}

\section{Proof of Proposition \ref{prop.exp}}\label{sec:proof.prop}
In this section we finally prove Proposition \ref{prop.exp}.
By Taylor expanding the matrix elements $\iota_{ij}(\bmmu, \e)$, $a_{ij}(\bmmu, \e)$ at $\e = 0$
 we get
 \begin{align}
    \iota_{11}(\bmmu, \e) & = \iota_{11}(\bmmu, 0) +\wt\alpha_1(\bmmu) \e +  \mathcal{O}(\e^2)  \\
    \iota_{10}(\bmmu, \e)  & = \wt \iota_{10}(\bmmu, 0) +  \gamma_1(\bmmu) \e + \mathcal{O}(\e^2)   \\
    \iota_{00}(\bmmu, \e) & = \iota_{00}(\bmmu,0) + \wt\alpha_0(\bmmu) \e + \mathcal{O}(\e^2) \\
    a_{11}(\bmmu, \e) & =  \underline{w} \, \iota_{11}(\bmmu, 0) +  \wt\beta_1(\bmmu) \e + \mathcal{O}(\e^2)  \ , \\ 
    a_{10}(\bmmu, \e) & = \wt a_{10}(\bmmu, 0)+ \beta_1(\bmmu) \e + \mathcal{O}(\e^2)  \\
    a_{01}(\bmmu, \e) & = \wt a_{01}(\bmmu, 0)+ \beta_0(\bmmu) \e + \mathcal{O}(\e^2) \ , \\
    a_{00}(\bmmu, \e) & = \underline{w} \, \iota_{00}(\bmmu, 0) +  \wt\beta_0(\bmmu) \e + \mathcal{O}(\e^2) \ .
\end{align}
We start with
\begin{lemma} One has 
\begin{align}
   \wt \iota_{10}(\bmmu, 0) = \wt a_{10}(\bmmu, 0) = \wt a_{01}(\bmmu, 0) = \wt\alpha_1(\bmmu) = \wt\alpha_0(\bmmu) = \wt\beta_1(\bmmu) = \wt\beta_0(\bmmu) = 0 \ , 
\end{align}
proving the Taylor expansions in \eqref{iota11}--\eqref{a00}.
\end{lemma}
\begin{proof}

    First consider $\wt \iota_{10}(\bmmu,0)$: by \eqref{eq:matrI} we have
    \begin{align*}
        \iota_{10}(\bmmu,0) &= \langle \fI(\bmmu,0)f_{\bmzero}^{+} , f_{\ubmk}^{-} \rangle \, \stackrel{ \eqref{fIfL} }{=} \, \langle P_{\bmmu,0}^{\ast}  P_{\bmmu,0} f_{\bmzero}^{+} , f_{\ubmk}^{-} \rangle \, \stackrel{\eqref{eq:SpecProj},\eqref{eq:v20} }{=} 0 .
    \end{align*}
    Next, we compute $\wt a_{10}(\bmmu,0)$: by \eqref{eq:matrA} and since, by Lemma \ref{jetsinFl}, $\fL_0 \in \mathfrak{F}_0$, we have
    \begin{align*}
        \wt a_{10}(\bmmu,0) &= \langle \fL_0 f_{\ubmk}^{-} , f_{\bmzero}^{+} \rangle = 0  \ , 
        \quad 
         \wt a_{01}(\bmmu,0) = \langle \fL_0 f_{\bmzero}^{+} , f_{\ubmk}^{-} \rangle = 0 \ . 
    \end{align*} 
    On the other hand, by Lemma \ref{jetsinFl}, $\mathfrak{I}_1, \mathfrak{L}_1 \in \mathfrak{F}_1$, hence using Lemma \ref{Fj} (iii)
    \begin{align*}
        \wt{\alpha}_1(\bmmu) &
        =  \langle \mathfrak{I}_1 f_{\ubmk}^{-} , f_{\ubmk}^{-} \rangle 
         = 0 \ , 
         \quad 
          \tilde{\alpha}_0(\bmmu)  =  \langle \fI_1 f_{\bmzero}^{+} , f_{\bmzero}^{+} \rangle  = 0 \ , \\
        \wt\beta_1(\bmmu) &= \langle \fL_1 f_{\ubmk}^{-} , f_{\ubmk}^{-} \rangle =  0 \ , 
        \quad 
         \wt\beta_0(\bmmu) = \langle \fL_1 f_{\bmzero}^{+} , f_{\bmzero}^{+} \rangle = 0.
    \end{align*}
\end{proof}

We now compute the remaining coefficients $\gamma_1(\bmmu)$, $\beta_1(\bmmu), \beta_0(\bmmu)$. 
To this aim, using the definitions \eqref{eq:v20}
\eqref{eq:lambdasigmak-restr} and  \eqref{def:undw},  we obtain
\begin{align}
\label{scalar.prod.f}
& (f_{\bmzero}^-, f_{\bmzero}^+)  = \frac12  \left(-|{\bmmu}|^2 + N^2 \right) \ , \quad  
(f_{\ubmk}^-, f_{\ubmk}^+) = \frac12  \left( -|\ubmk + {\bmmu}|^2 + N^2 \right) \ , \\
\label{den.id}
& w_{\bmzero}^- - \uw = - \frac{2 \mu_1 N}{|{\bmmu}|} = - 2 \, \Omega({\bmmu}) , \quad 
w_{\ubmk}^+ - \uw =  \frac{2 (\tm+\mu_1) N}{|\ubmk+{\bmmu}|}  = 2 \, \Omega(\ubmk + \bmmu), \\
& 
 \td_{\ubmk} =\frac{1}{N (\tm + \mu_1) |\ubmk + {\bmmu}|} , \quad \td_{\bmzero} = \frac{1}{N \mu_1 |{\bmmu}|}. \label{eq:dkd0} 
\end{align}
\\

\noindent{\bf Computation of $\gamma_1$:} using Lemma \ref{expansionthm} for the expression of $\fI_{1}$ and  the residue theorem
\begin{align}\notag
     \e \gamma_1 &  = \langle \fI_{1} f_{\ubmk}^-, f_{\bmzero}^+ \rangle = 
\langle P_1 f_{\ubmk}^-, f_{\bmzero}^+ \rangle  + 
\langle  f_{\ubmk}^-, P_1 f_{\bmzero}^+ \rangle 
\\
\notag
& 
\stackrel{\text{Lemma} \ref{Fj} (iii)}{ =}  \langle P_1^{[-\ubmk]} f_{\ubmk}^-, f_{\bmzero}^+ \rangle  + 
\langle  f_{\ubmk}^-, P_1^{[\ubmk]} f_{\bmzero}^+ \rangle 
\stackrel{\eqref{proj.Fl}}{=}
    \langle \cP[\sL_{1}^{[-\ubmk]}] f_{\ubmk}^-, f_{\bmzero}^+ \rangle 
    + \langle  f_{\ubmk}^- , \cP[\sL_{1}^{[\ubmk]}] f_{\bmzero}^+ \rangle \\
    & 
    \stackrel{\text{Lemma} \ref{actionofL}}{=}  \td_{\bmzero} \, \frac{\ent{1}{-\ubmk}{\bmzero}{\ubmk}{-}{-}}{w_{\bmzero}^- - \uw} (f_{\bmzero}^-, f_{\bmzero}^+) -
     \td_{\ubmk} \, \frac{{\ent{1}{\ubmk}{\ubmk}{\bmzero}{+}{+}}}{w_{\ubmk}^+ - \uw} (f_{\ubmk}^-, f_{\ubmk}^+) .
     \label{eq:epsg1}
\end{align}

\vspace{1em}
\noindent{\bf Computation of $\beta_1$ and $\beta_0$:} First, by \eqref{eq:matrA} and \eqref{a10}, 
recall that
$$
\im \e \beta_1 = \langle \fL_1 f_{\ubmk}^-, f_{\bmzero}^+ \rangle \ , \quad
\im \e \beta_0 = \langle \fL_1 f_{\bmzero}^+, f_{\ubmk}^- \rangle  \ . 
$$
Appealing to Lemma \ref{expansionthm} for the expression of $\fL_1$ and to Lemmata \ref{Fj}, \ref{jetsinFl} we obtain
\begin{align}
\notag
    \im \e \beta_1 &  = \langle \fL_{1} f_{\ubmk}^-, f_{\bmzero}^+ \rangle =
\langle \sL_{1} f_{\ubmk}^-, f_{\bmzero}^+ \rangle   + 
\langle \sL_{0} P_1 f_{\ubmk}^-, f_{\bmzero}^+ \rangle
+
\langle \sL_{0}  f_{\ubmk}^-, P_1 f_{\bmzero}^+ \rangle
    \\
        \notag
    & = 
    \langle \sL_{1}^{[-\ubmk]} f_{\ubmk}^-, f_{\bmzero}^+ \rangle + \langle \sL_{0}^{[\bmzero]} \cP[\sL_{1}^{[-\ubmk]}] f_{\ubmk}^- , f_{\bmzero}^+ \rangle + \langle \sL_{0}^{[\bmzero]} f_{\ubmk}^-, \cP[\cL_{1}^{[\ubmk]}] f_{\bmzero}^+ \rangle \\
    \label{beta1.algebric}
    & = \entL{1}{-\ubmk}{\bmzero}{\ubmk}{+}{-} 
    + \td_{\bmzero} \, \frac{\ent{1}{-\ubmk}{\bmzero}{\ubmk}{-}{-}}{w_{\bmzero}^- - \uw} \, \entL{0}{\bmzero}{\bmzero}{\bmzero}{+}{-} 
     -
     \td_{\ubmk} \, \frac{{\ent{1}{\ubmk}{\ubmk}{\bmzero}{+}{+}}}{w_{\ubmk}^+ - \uw} 
    \,  \entL{0}{\bmzero}{\ubmk}{\ubmk}{+}{-} 
\end{align}
where in the last step we used  \eqref{cBactspar}, \eqref{cBactsparbis} and the residue theorem.

Arguing similarly we get 
\begin{align}
\notag
    \im \e \beta_0 &  = \langle \fL_{1} f_{\bmzero}^+, f_{\ubmk}^- \rangle = 
    \langle \sL_{1}f_{\bmzero}^+,  f_{\ubmk}^- \rangle   + 
\langle \sL_{0} P_1 f_{\bmzero}^+ ,f_{\ubmk}^- \rangle
+
\langle \sL_{0} f_{\bmzero}^+ , P_1 f_{\ubmk}^- \rangle
    \\
        \notag
   &  =  \langle \sL_{1}^{[\ubmk]} f_{\bmzero}^+  , f_{\ubmk}^-  \rangle
    + 
    \langle \sL_{0}^{[\bmzero]} \cP[\sL_{1}^{[\ubmk]}] f_{\bmzero}^+  ,  f_{\ubmk}^- \rangle + 
    \langle \sL_{0}^{[\bmzero]}f_{\bmzero}^+ , \cP[\sL_{1}^{[-\ubmk]}] f_{\ubmk}^- \rangle \\
    \label{beta0.algebric}
    & = \entL{1}{\ubmk}{\ubmk}{\bmzero}{-}{+}
    - \td_\ubmk \, \frac{\ent{1}{\ubmk}{\ubmk}{\bmzero}{+}{+}}{w_\ubmk^+ - \uw} \, \entL{0}{\bmzero}{\ubmk}{\ubmk}{-}{+} 
     +
     \td_{\bmzero} \, \frac{{\ent{1}{-\ubmk}{\bmzero}{\ubmk}{-}{-}}}{w_{\bmzero}^- - \uw} 
    \,  \entL{0}{\bmzero}{\bmzero}{\bmzero}{-}{+}
\end{align}
where again we used  \eqref{cBactspar}, \eqref{cBactsparbis} and the residue theorem.\\

\noindent{\bf Computation of $ \beta_1 - \gamma_1 \uw $ and $ \beta_0 - \gamma_0 \uw $:}
By formulas \eqref{beta1.algebric} and  \eqref{eq:epsg1} 
we have
    \begin{align*}
       (\beta_1 - \gamma_1 \uw ) \e &= -\im \entL{1}{-\ubmk}{\bmzero}{\ubmk}{+}{-} + \frac{ \td_{\bmzero} \ent{1}{-\ubmk}{\bmzero}{\bmk}{-}{-} }{  w_{\bmzero}^{-} - \uw } \left[ -\im \entL{0}{\bmzero}{\bmzero}{\bmzero}{+}{-} - \uw \, ( f_{\bmzero}^{-} , f_{\bmzero}^{+}  ) \right] - 
       \frac{ \td_{\ubmk} \ent{1}{\ubmk}{\bmk}{\bmzero}{+}{+} }{  w_{\ubmk}^{+} - \uw } \left[ -\im \entL{0}{\bmzero}{\ubmk}{\ubmk}{+}{-} - \uw \, ( f_{\ubmk}^{-} , f_{\ubmk}^{+} )  \right]  \ .  
   \end{align*}
   Then, using \eqref{eq:L000p0m}, \eqref{eq:L00kpkm}, \eqref{def:undw} and \eqref{scalar.prod.f} we compute 
   \begin{align*}
       -\im \entL{0}{\bmzero}{\bmzero}{\bmzero}{+}{-} - \uw \, ( f_{\bmzero}^{-} , f_{\bmzero}^{+} ) &= \Omega({\bmmu}) (|{\bmmu}|^2 - N^2 ) ,  \\
       -\im \entL{0}{\bmzero}{\ubmk}{\ubmk}{+}{-} - \uw \, ( f_{\ubmk}^{-} , f_{\ubmk}^{+} )  &= 0  \ . 
    \end{align*}
    Therefore, using also identity \eqref{den.id}, 
    \begin{align}\notag
       & (\beta_1 - \gamma_1 \uw ) \e = -\im \entL{1}{-\ubmk}{\bmzero}{\ubmk}{+}{-} - \frac{1}{2} \td_{\bmzero} \, \ent{1}{-\ubmk}{\bmzero}{\ubmk}{-}{-}  (|{\bmmu}|^2-N^2) \nonumber \\
       \notag
       & \stackrel{\eqref{eq:L1mk0+km}, \eqref{eq:B1mk}, \eqref{eq:dkd0}}{=} \frac{\e}{4} (\ubmk^{\perp} \cdot {\bmmu})  \frac{ |\ubmk+{\bmmu}| -|\ubmk| }{N\,|\ubmk|\,  |\ubmk+{\bmmu}|} \left[  ( |\ubmk+{\bmmu}|+|\ubmk| ) |{\bmmu}| - N^2  -\frac{1}{2} \frac{|\ubmk+{\bmmu}| + |\ubmk| + |{\bmmu}| }{|{\bmmu}|} (|{\bmmu}|^2 - N^2) 
       \right] \\
       \notag
       & = \frac{\e}{8} (\ubmk^{\perp} \cdot {\bmmu})  \frac{ \big( |\ubmk+{\bmmu}| -|\ubmk|\big) \, (|\ubmk + {\bmmu}| + |\ubmk| - |{\bmmu}| \big) \, \big(N^2 + |{\bmmu}|^2\big) }{N\,|\ubmk|\, |{\bmmu}|\,  |\ubmk+{\bmmu}|}
   \end{align}
   proving formula \eqref{b1-g1}.
   
Similarly, using \eqref{beta0.algebric} and  \eqref{eq:epsg1} 
     \begin{align*}
       (\beta_0 - \gamma_1 \uw ) \e &= -\im \entL{1}{\ubmk}{\ubmk}{\bmzero}{-}{+} + \frac{ \td_{\bmzero} \ent{1}{-\ubmk}{\bmzero}{\bmk}{-}{-} }{  w_{\bmzero}^{-} - \uw } \left[ -\im \entL{0}{\bmzero}{\bmzero}{\bmzero}{-}{+} - \uw \, ( f_{\bmzero}^{-} , f_{\bmzero}^{+}  ) \right] - 
       \frac{ \td_{\ubmk} \ent{1}{\ubmk}{\bmk}{\bmzero}{+}{+} }{  w_{\ubmk}^{+} - \uw } \left[ -\im \entL{0}{\bmzero}{\ubmk}{\ubmk}{-}{+} - \uw \, ( f_{\ubmk}^{-} , f_{\ubmk}^{+} )  \right] ,
   \end{align*}
and using \eqref{eq:L000m0p}, \eqref{eq:L00kmkp}, \eqref{def:undw} and \eqref{scalar.prod.f} we compute 
   \begin{align*}
       -\im \entL{0}{\bmzero}{\bmzero}{\bmzero}{-}{+} - \uw \, ( f_{\bmzero}^{-} , f_{\bmzero}^{+} ) &=  0 ,  \\
       -\im \entL{0}{\bmzero}{\ubmk}{\ubmk}{-}{+} - \uw \, ( f_{\ubmk}^{-} , f_{\ubmk}^{+} )  &= - \Omega(\ubmk+{\bmmu}) ( |\ubmk+{\bmmu}|^2 -N^2 ) ,
    \end{align*}
yielding
   \begin{align}\notag
       & (\beta_0 - \gamma_1 \uw ) \e \stackrel{\eqref{den.id}}{=} -\im \entL{1}{\ubmk}{\ubmk}{\bmzero}{-}{+} +
       \frac{1}{2} \td_{\ubmk} \, \ent{1}{\ubmk}{\bmk}{\bmzero}{+}{+} (|\ubmk+{\bmmu}|^2-N^2) \nonumber \\
       \notag
       &= \frac{\e}{4} (\ubmk^{\perp} \cdot {\bmmu}) \frac{ |\ubmk| +|{\bmmu}| }{N \, |\ubmk| \, |{\bmmu}| } \left[ - ( |\ubmk|-|{\bmmu}| ) |\ubmk+{\bmmu}| - N^2 +\frac{1}{2} \frac{|\ubmk| - |{\bmmu}| - |\ubmk+{\bmmu}| }{|\ubmk+{\bmmu}|} (|\ubmk+{\bmmu}|^2 - N^2)
       \right] \\ \notag
       & = -\frac{\e}{8} (\ubmk^{\perp} \cdot {\bmmu}) \frac{\big( |\ubmk| +|{\bmmu}|\big) \, \big( |\ubmk + {\bmmu}| + |\ubmk| - |{\bmmu}|\big) \, \big(|\ubmk + {\bmmu}|^2 + N^2 \big)  }{N \, |\ubmk| \, |{\bmmu}| \, |\ubmk + {\bmmu}|}
   \end{align}
 and   proving formula \eqref{b0-g0}.
   
\section{Comparison with physical literature}\label{sec:comparison}

Here we compare our result with the instability results in \cite{bourget2014finite,dauxois2018instabilities} (see also the experimental results in \cite{bourget2013experimental}). Comparing the notation of \cite[Section 3.2]{dauxois2018instabilities} with that of the present paper, we have the following.
\begin{center}
\begin{tabular}{||l | l ||} 
 \hline
 Dauxois et al. \cite{dauxois2018instabilities} & Present paper  \\ [0.5ex] 
 \hline
$\bmk_0 = (\ell_0, m_0)^T$ & $\ubmk = (\tm, \tn)^T$  \\ 
 \hline
 $\bmk_+$ & $\ubmk + {\bmmu}$  \\
 \hline
 $\bmk_-$ & $-{\bmmu}$  \\
 \hline
 $\omega_0$ & $\Omega(\ubmk)$  \\
 \hline
 $\omega_+$ &  $\Omega(\ubmk + {\bmmu})$ \\
 \hline
 $\omega_-$ & $ \Omega({\bmmu})$
 \\ [1ex] 
 \hline
\end{tabular}
\end{center}

Following \cite{dauxois2018instabilities}, we introduce
\begin{align}
    &I_+ ({\bmmu})\coloneqq (\ubmk^\perp \cdot {\bmmu}) \frac{ \Omega(\ubmk+{\bmmu}) (|\ubmk|^2 - |{\bmmu}|^2)+ (\tm+\mu_1) \, N \, (|\ubmk| - |{\bmmu}| )}{2 \Omega(\ubmk+{\bmmu})  |\ubmk+{\bmmu}|^2} , \label{eq:Ip} \\
    &I_-({\bmmu}) \coloneqq - (\ubmk^\perp \cdot {\bmmu}) \frac{ \Omega({\bmmu}) (|\ubmk|^2 - |\ubmk+{\bmmu}|^2) + \mu_1 \, N \, ( |\bmk| - |\bmk+\bmmu| )}{2 \Omega(\bmmu) |\bmmu|^2}  . \label{eq:Im}
\end{align}
In the regime of small Floquet parameter $|\bmmu| \ll 1$, the following expansions are consistent \cite[Section 3.2.2]{dauxois2018instabilities}. 
\begin{proposition}
The functions $\mathtt{e}({\bmmu})$, $I_{+}({\bmmu})$ and $I_{-}({\bmmu})$ defined respectively in \eqref{temu}, \eqref{eq:Ip} and \eqref{eq:Im} satisfy
    \begin{align}\label{D&I}
      {\mathtt{e}({\bmmu})}  = \frac{I_+({\bmmu}) \, I_-({\bmmu})}{{ 4}N^2 |\ubmk|^2} \left( \frac{|\ubmk + {\bmmu}| + |\ubmk| - |{\bmmu}|}{|\ubmk + {\bmmu}| + |\ubmk| + |{\bmmu}|}\right)^2 \, \frac{|\ubmk| + |{\bmmu}|}{|\ubmk| - |{\bmmu}|}  \ . 
   \end{align}
   In addition, for  $\bmmu(\tau)$ as in \eqref{paranear0}, in the regime of small Floquet parameter 
   \begin{equation}\label{limits}
\frac{ \mathtt{e}(\bmmu(\tau)) }{ I_+(\bmmu(\tau)) \, I_-(\bmmu(\tau))} \to   \frac{1}{4N^2 |\ubmk|^2} \ \   \mbox{ as }
\tau \rightarrow 0 \ . 
    \end{equation}
\end{proposition}
\begin{proof}
Observe  that
$$ 
    I_+({\bmmu}) =   (\ubmk^\perp \cdot {\bmmu}) (|\ubmk| - |{\bmmu}|) \frac{ |\ubmk| + |{\bmmu}| + |\ubmk+{\bmmu}| }{2 |\ubmk + {\bmmu}|^2} , \quad 
    I_-({\bmmu}) =   (\ubmk^\perp \cdot {\bmmu}) (|\ubmk+{\bmmu}| - |\ubmk|) \frac{ |\ubmk| + |\ubmk+{\bmmu}| + |{\bmmu}| }{2 |{\bmmu}|^2} .
$$
Then using  \eqref{temu}, formula  \eqref{D&I} follows immediately as well as the  limit in 
\eqref{limits} as $\tau \to 0$.
\end{proof}

\appendix
\section{Proof of Lemma \ref{lem:Rk}}\label{app:Rk}
Define the function
\[
F(x,y) := \frac{\tm + x}{ \sqrt{(\tm+x)^2 + (\tn+y)^2} } + \frac{x}{ \sqrt{x^2 + y^2} } - \frac{\tm}{|\ubmk|},
\]
which is real analytic except at the points \( (x, y)=(0,0) \) and \( (x,y)=(-\tm, -\tn) \).
We define the zero sets
\[
\mathcal{R}_{\ubmk}^\pm := \{ (x,y) \in \mathbb{R}^2 \mid F(x,y) = 0,\ y \gtrless 0 \}.
\]

\paragraph*{\bf Analysis of \( {\mathcal{R}}_{\ubmk}^+ \)}

Fix \( \underline{y} > 0 \). Observe that
\[
\lim_{x \to +\infty} F(x, \underline{y}) = 2 - \frac{\tm}{|\ubmk|} > 0, \quad
\lim_{x \to -\infty} F(x, \underline{y}) = -2 - \frac{\tm}{|\ubmk|} < 0,
\]
so by the intermediate value theorem, for each \( \underline{y} > 0 \), there exists at least one solution \( x \) to the equation \( F(x, \underline{y}) = 0 \).

Furthermore, since \( \partial_x F(x, \underline{y}) > 0 \), such solution is unique. This defines a function \( \varphi_+ \colon (0, +\infty) \to \mathbb{R} \) such that
\[
F(\varphi_+(y), y) = 0.
\]
Since \( (\varphi_+(y), y) \notin \{(0,0), (-\tm, -\tn)\} \), the implicit function theorem implies that \( \varphi_+ \) is real analytic on \( (0, + \infty) \).

\subsubsection*{Asymptotics as \( y \to 0^+ \)}

Let \( \varepsilon > 0 \). Then,
\[
F\left( \left(\frac{\tm \tn }{|\ubmk|^3} \pm \varepsilon \right) y^2, y \right) = \pm \varepsilon y + \mathcal{O}(y^2),
\]
so for \( y > 0 \) sufficiently small, the sign of \( F \) at these points is positive/negative respectively. Thus,
\[
\left( \frac{\tm \tn}{|\ubmk|^3} - \varepsilon \right) y^2 < \varphi_+(y) < \left( \frac{\tm \tn}{|\ubmk|^3} + \varepsilon \right) y^2,
\]
from which we deduce the asymptotic behavior
\[
\varphi_+(y) \sim \frac{\tm \tn}{|\ubmk|^3} y^2 \quad \text{as } y \to 0^+.
\]

\subsubsection*{Asymptotics as \( y \to +\infty \)}

For any \( q > 0 \),
\[
F\left( \frac{\tm}{\sqrt{3 \tm^2 + 4 \tn^2}} \left( \frac{\tn}{2} + y \right) - \frac{\tm}{2} \pm q,\ y \right)
= \pm \frac{1}{4} \left( \frac{3 \tm^2 + 4 \tn^2}{\tm^2 + \tn^2} \right)^{3/2} \cdot \frac{q}{y} + \mathcal{O}(y^{-2}).
\]
Hence, for \( y \to +\infty \),
\[
\varphi_+(y) = \frac{\tm}{\sqrt{3 \tm^2 + 4 \tn^2}} \left( \frac{\tn}{2} + y \right) - \frac{\tm}{2} + \mathcal{O}(1).
\]

\subsubsection*{Monotonicity}

Stationary points of \( \varphi_+ \) correspond to critical points of \( F \), i.e., to solutions of the system
\begin{equation}\label{stat.point}
\begin{cases}
\partial_y F(x,y) = 0 \\
F(x,y) = 0
\end{cases}
\quad \Leftrightarrow \quad
\begin{cases}
\displaystyle \frac{xy}{(x^2 + y^2)^{3/2}} + \frac{(x+\tm)(y+\tn)}{\left((x+\tm)^2 + (y+\tn)^2\right)^{3/2}} = 0 \\
\displaystyle \frac{x}{\sqrt{x^2 + y^2}} + \frac{x+\tm}{\sqrt{(x+\tm)^2 + (y+\tn)^2}} = \frac{\tm}{|\ubmk|}
\end{cases}
\end{equation}
Since \( y > 0 \), \( \tm > 0 \), the first equation implies \( -\tm \leq x \leq 0 \). But the asymptotics show that for \( y > 0 \) small, \( \varphi_+(y) > 0 \), so \( x = \varphi_+(y) > 0 \) contradicts the necessary range for a stationary point. Therefore, \( \varphi_+ \) has no critical points on \( (0, +\infty) \), and hence is strictly increasing.

\smallskip
\paragraph*{\bf Analysis of $ \mathcal{R}_{ \ubmk }^- $} We now consider the set \( F(x,y) = 0 \) with \( y < 0 \). As in the previous case, one finds that for any \( y \neq -\tn \), there exists a unique value \( \varphi_-(y) \) solving 
\[
F(\varphi_-(y), y) = 0,
\]
and the function \( y \mapsto \varphi_-(y) \) is real analytic on the domain \( (-\infty, -\tn) \cup (-\tn, 0) \). Moreover, by similar arguments, one obtains the same asymptotic expansions as in \eqref{varpi.exp}.

We next distinguish between the cases \( 2\tfrac{\tm}{|\ubmk|} < 1 \) and \( 2\tfrac{\tm}{|\ubmk|} > 1 \).

\medskip
\noindent

$\bullet$ Case $2\dfrac{\tm}{|\ubmk|} <1$: In this case
$$
\lim_{x \to -\tm^-} F(x, -\tn) = -1- 2\dfrac{\tm}{|\ubmk|}  < 0, \quad \lim_{x \to - \tm^+} F(x, -\tn) = 1- 2\dfrac{\tm}{|\bmk|} >0 \ , \quad \lim_{x \to +\infty} F(x, -\tn) >0,
$$
and $\partial_x F(x, -\tn)>0$, so there is no solution of $F(x,-\tn) = 0$. However we can compute how it approaches $y=-\tn$, proving the  asymptotic \eqref{varpi.exp2}.
Indeed, consider first  $0<y+\tn\ll 1$; fix  $q \in \R$   and note that  
\begin{align*}
F(-\tm + (1+q)a_*(y+ \tn), y)
& = 
\frac{2 (q+1) \tm}{\sqrt{(4 q^2 + 8 q+1) \tm^2+\tn^2}}-\frac{2 \tm}{\sqrt{\tm^2+\tn^2}} + O(y+n)
\end{align*}
Now, for $q$ sufficiently small, the function 
$$
\frac{2 (q+1) \tm}{\sqrt{(4 q^2 + 8 q+1) \tm^2+\tn^2}}-\frac{2 \tm}{\sqrt{\tm^2+\tn^2}} = 
\frac{2 \tm \left(\tn^2 -3 \tm^2\right)}{|\ubmk|^{3/2}} q + O(q^2)\begin{cases}
    >0 \mbox{ if } q >0 \\
    <0 \mbox{ if } q < 0
\end{cases}
$$
so is the function $F(-\tm + (1+q)a_*(y+ \tn), y)$. This shows that  for any $q>0$ sufficiently small, 
\begin{align*}
   - \tm + (1-q) a_*(\tn+y) < \varphi_-(y) <   - \tm + (1+q) a_* (\tn+y)  \quad \mbox{ as } y \searrow - \tn \ ,
\end{align*}
hence the asymptotic in \eqref{varpi.exp2} holds true for $y > - \tn$.

Now consider $-1 \ll y+\tn<0$; an analogous analysis shows that, for any $q>0$ sufficiently small, 
\begin{align*}
   - \tm - (1-q) a_*(\tn+y) < \varphi_-(y) <   - \tm - (1+q) a_* (\tn+y)  \quad \mbox{ as } y \nearrow - \tn \ ,
\end{align*}
proving the  asymptotic in \eqref{varpi.exp2} also for $y < - \tn$.\\

$\bullet$ Case $2\dfrac{\tm}{|\ubmk|} >1$: In this case, we observe that
\[
\lim_{x \to -\tm^+} F(x, -\tn) = 1 - 2\frac{\tm}{|\ubmk|} < 0, \quad
\lim_{x \to +\infty} F(x, -\tn) = 2 - \frac{\tm}{|\ubmk|} > 0,
\]
so there exists a unique solution \( \varphi_-(-\tn) \) such that \( F(\varphi_-(-\tn), -\tn) = 0 \). By the implicit function theorem, it follows that \( \varphi_-(y) \) is real analytic on \( (-\infty, 0) \).

\medskip
\noindent



{\bf Proof of \eqref{linR}:} Let $\ell \in \Z \setminus\{-1, 0 \}$ and assume that $\ell \ubmk \in \cR_{\ubmk}$.
If $\ell \geq 1$, using the $0$-homogeneity of $\Omega(\bmk)$, we  get
$$
\Omega(\ubmk) = \Omega((1+\ell) \ubmk) + \Omega(\ell \ubmk)  = 2\Omega(\ubmk) 
\ 
\Rightarrow \ \Omega(\ubmk) = 0 
$$
which is not possible since $\tm \neq 0$.
If $\ell \leq -2$, again by homogeneity we get
$
3\Omega(\ubmk) =0$, contradiction.

\section{Proof of Proposition \ref{prop:wBound}}\label{app:de}

\begin{proof}
{\bf Step 1:} {\em There exists a discrete set $\cD_{\ubmk}^+$ such that 
such that
for any $\bmmu \in \cR_\ubmk\setminus \cD_\ubmk^+ $, 
there exists $\alpha=\alpha (\bmmu)>0$ such that:
\begin{align} \label{claimStep1.app}
\inf_{\ell \in \mathbb{Z}, \ell \neq 0} | w^{+}_{\ell \ubmk}( \bmmu ) - w^{+}_{\bmzero}( \bmmu )| &\geq \alpha N.  
\end{align}
}
    We rewrite the difference 
    \begin{align*}
         w^{+}_{\ell \ubmk}( \bmmu ) - w^{+}_{\bmzero}( \bmmu )&= \ell\underline{\bmc} \cdot \ubmk +  \Omega(\ell\ubmk+\bmmu) -  \Omega(\bmmu)\\
         &\stackrel{ \eqref{eq:TrivColl} }{=} \ell \Omega(\ubmk)+ \Omega(\ell\ubmk+\bmmu) -  \Omega(\bmmu) .
    \end{align*}
    As $\bmmu \in \mathcal{R}_{\ubmk}$, then 
    $\Omega (\bmmu)= \Omega(\ubmk)- \Omega(\ubmk+\bmmu)$,
    and the above expression simplifies as
    \begin{align*}
        w^{+}_{\ell \ubmk}( \bmmu ) - w^{+}_{\bmzero}( \bmmu )&= (\ell-1) \Omega (\ubmk)+ \Omega(\ell\ubmk+\bmmu)+\Omega(\ubmk+\bmmu) .
    \end{align*}
As $\abs{\Omega(\ell\ubmk+\bmmu)+\Omega(\ubmk+\bmmu)}\leq 2$ for every $\ell$, there exists $\ell_0>0$ such that 
\begin{equation}
    \abs{ 
    w^{+}_{\ell \ubmk}( \bmmu ) - w^{+}_{\bmzero}( \bmmu )
    }
    \geq 
    \abs{(\ell-1) \Omega (\ubmk)} - \abs{ \Omega(\ell\ubmk+\bmmu)+\Omega(\ubmk+\bmmu) } \geq 1 
    \quad \forall |\ell| \geq \ell_0 \ , 
    \ \forall \bmmu \in \cR_{\ubmk}  . 
\end{equation}
Now consider the finitely many functions $\{\varpi_\ell^+(\bmmu) \}_{0< |\ell| < \ell_0}$, with
$$
\varpi_\ell^+(\bmmu):=  w^{+}_{\ell \ubmk}( \bmmu ) - w^{+}_{\bmzero}( \bmmu )
 =
 (\ell-1) \Omega (\ubmk)+ \Omega(\ell\ubmk+\bmmu)+\Omega(\ubmk+\bmmu) \ .
$$
Notice that 
$ \varpi_\ell^+(\bmmu)$ is real analytic in $\R^2 \setminus \{  - \ell \ubmk, - \ubmk \}$. 
We are going to show that  $\varpi_\ell^+\vert_{\cR_{\ubmk}}$  vanishes only on finitely many values.\\
\paragraph*{\bf Analysis near $\mathbf{0}$} To study the behavior of $\varpi_\ell^+(\bmmu)$ near zero we shall use the parametrization of $\cR_{\ubmk}$ of Lemma 
\ref{lem:Rk}.
In particular, there exists $\tau_0>0$ such that  defining
\begin{equation} \label{eq:mu0}
    \bmmu_0(\tau):= \begin{cases}
        (\varphi_+(\tau), \tau)^T & \mbox{ if } \tau \in (0, \tau_0) \\
        \bmzero  &  \mbox{ if } \tau = 0 \\
       (\varphi_-(\tau), \tau)^T  & \mbox { if }  \tau \in (-\tau_0, 0) 
    \end{cases}
\end{equation}
with $\varphi_\pm$ 
in \eqref{varpi.exp}, 
one has
$$
 \{ \bmmu = (\mu_1,\mu_2)^T \in \cR_{\ubmk} \colon \mu_2 \in B_{\tau_0}(0) \}  = \{  \bmmu_0(\tau)  \colon |\tau| \leq \tau_0   \}  
$$ 
Evaluating $\varpi_\ell^+(\bmmu)$ using this parametrization yields
\begin{equation}
  \varpi_\ell^+\big(  \bmmu_0(\tau)  \big) =  
  (\ell+\sgn(\ell))\Omega(\ubmk) + \mathcal{O}(\tau)  , \ \ \forall |\tau| \leq \tau_0
\end{equation}
which is bounded away from zero if $\tau_0$ is small enough. \\

\paragraph{\bf Analysis near $\pm \infty $}
To study the behaviour near $+ \infty$, we use the parametrization 
$$
(0, \tau_0) \mapsto \bmmu_{+\infty}(\tau) = \frac{1}{\tau}
\begin{pmatrix}
\frac{\tm \tau}{\sqrt{3\tm^2 + 4\tn^2}} \left(\frac{\tn}{2} + \frac{1}{\tau} \right) - \frac{\tm \tau }{2} \\
1
\end{pmatrix}
$$
equivalent to $(\varphi_+(y), y)$ as $y \to \infty$.
In particular, we have
\begin{equation}
  \varpi_\ell^+\big(  \bmmu_{+\infty}(\tau)  \big) =  
  \ell \Omega(\ubmk) + \mathcal{O}(\tau)
 , \ \ \forall |\tau| \leq \tau_0,
\end{equation}
which again is bounded away from zero since $\ell \neq 0$.
Similarly, 
\begin{equation}
  \varpi_\ell^+\big(\bmmu_{-\infty}(\tau) \big) =  
   \ell \Omega(\ubmk) + \mathcal{O}(\tau)
 , \ \ \forall |\tau| \leq \tau_0 
\end{equation}
is bounded away from zero as $\ell \neq 0$.
\\

\paragraph{\bf Analysis in $(\tau_0, {\tau_0}^{-1})$} 
In this interval, the function $y \mapsto \varphi_+(y)$ is real analytic. Recall that 
$\varpi_\ell^+(\bmmu)$ is real analytic in $\R^2 \setminus \{ -\ell \ubmk, -\ubmk \}$. 
However, $-\ell \ubmk \in \cR_{\ubmk}$ implies that $\ell = 0$ or $\ell = 1$, see \eqref{linR}. 
The case $\ell = 0$ is excluded since $|\ell| > 0$, while the case $\ell = 1$ gives a singularity at $-\ubmk$. 
Nevertheless, the curve $\{ (\varphi_+(y), y) \colon y \in (\tau_0, \tau_0^{-1}) \}$ does not pass through this point. 

Hence, the map 
\[
(\tau_0, {\tau_0}^{-1}) \ni y \mapsto \varpi_\ell^+\big( (\varphi_+(y), y) \big)
\]
is real analytic and thus can vanish only at finitely many points.\\

\paragraph{\bf Analysis in $(-{\tau_0}^{-1}, -\tau_0)$} 
We need to distinguish between the cases $2\frac{\tm}{|\ubmk|} < 1$ and $2\frac{\tm}{|\ubmk|} > 1$.

$\bullet$ If $2\frac{\tm}{|\ubmk|} > 1$, the parametrization $(\varphi_-(y), y)$ is real analytic and does not pass through any point of the form $-\ell \ubmk$. 
Therefore, $y \mapsto \varpi_\ell^-\big( (\varphi_-(y), y) \big)$ is again real analytic and can vanish only at finitely many points.

$\bullet$ If $2\frac{\tm}{|\ubmk|} < 1$, the parametrization $(\varphi_-(y), y)$ passes through the point $-\ubmk$ and is only Lipschitz continuous near this point. 
Locally around this point, we can use the following parametrization:
\begin{equation}\label{para-n}
\tau \mapsto \bmmu_{-\tn}(\tau):= 
\begin{pmatrix}
- \tm + \sqrt{\frac{4\tm^2}{(\tn^2 - 3\tm^2)}} \, |\tau| + \mathcal{O}(\tau^2)\\
- \tn + \tau
\end{pmatrix}
\end{equation}

With this parametrization, we obtain
\begin{equation}
\varpi_\ell^+\big( \bmmu_{-\tn}(\tau) \big) =
\begin{cases}
(\ell + 1 + \sgn(\ell - 1)) \Omega(\ubmk) + \mathcal{O}(\tau), & \mbox{if } \ell > 0 \\
\ell \Omega(\ubmk), & \mbox{if } \ell < 0
\end{cases}
\end{equation}
which is bounded away from zero. Since away from $y = -\tn$ the function $y \mapsto \varphi_-(y)$ is real analytic, it can again vanish only at finitely many points.

Hence, the set
\[
\cD_{\ubmk}^+ := \bigcup_{0<|\ell| < \ell_0} \left\{ (\varpi_\ell^+)^{-1}(0) \cap \cR_{\ubmk} \right\}
\]
is discrete. For $\bmmu \in \cR_{\ubmk} \setminus \cD_{\ubmk}^+$, one has
\begin{equation}\label{step1.concl}
\inf_{\ell \in \mathbb{Z}, \, \ell \neq 0} \left| w^{+}_{\ell \ubmk}(\bmmu) - w^{+}_{\bmzero}(\bmmu) \right| = 
\min \left( 
\inf_{|\ell| \geq \ell_0} |\varpi_\ell^+(\bmmu)|, \ 
\min_{0 < |\ell| < \ell_0} |\varpi_\ell^+(\bmmu)| 
\right) \geq \alpha(\bmmu) N > 0,
\end{equation}
proving the claim \eqref{claimStep1.app}.

\vspace{5mm}

{\bf Step 2:} {\em There exists a discrete set $\cD_{\ubmk}^-$ 
such that
for any $\bmmu \in \cR_\ubmk\setminus \cD_\ubmk^- $, 
there exists $\alpha=\alpha (\bmmu)>0$ such that:
\begin{align} \label{claimStep2.app}
\inf_{\ell \in \mathbb{Z}, \ell \neq 1} | w^{-}_{\ell \ubmk}( \bmmu ) - w^{-}_{\ubmk}( \bmmu )| &\geq \alpha N .  
\end{align}
}
Observe that 
    \begin{align*}
       w^{-}_{\ell \ubmk}(\bmmu) -  w^{-}_{\underline{\bmk}}( \bmmu ) &= \underline{\bmc} \cdot ( \ell \ubmk +\bmmu) - \Omega(\ell \ubmk + \bmmu) -  \underline{\bmc} \cdot ( \ubmk +\bmmu) + \Omega(\ubmk+\bmmu) \\
      &= (\ell-1) \Omega(\ubmk)  - \Omega(\ell \ubmk + \bmmu)  + \Omega(\ubmk+\bmmu).
    \end{align*}
    Again $\abs{-\Omega(\ell\ubmk+\bmmu)+\Omega(\ubmk+\bmmu)}\leq 2$ for every $\ell$, so  there exists $\ell_0>0$ such that 
\begin{equation}
    \abs{ 
    w^{+}_{\ell \ubmk}( \bmmu ) - w^{+}_{\bmzero}( \bmmu )
    }
    \geq 
    \abs{(\ell-1) \Omega (\ubmk)} - \abs{- \Omega(\ell\ubmk+\bmmu)+\Omega(\ubmk+\bmmu) } \geq 1 
    \quad \forall | \ell | \geq \ell_0 \ , 
    \ \forall \bmmu \in \cR_{\ubmk}   
\end{equation}
and we restrict to study the zero set of the finitely many  functions $\{\varpi_\ell^-(\bmmu) \}_{|\ell| < \ell_0, \ell \neq 1}$ defined by 
$$
\varpi_\ell^-(\bmmu):=  w^{-}_{\ell \ubmk}( \bmmu ) - w^{-}_{\ubmk}( \bmmu ) \ .
$$

\paragraph{\bf Analysis near $\bmzero$}
We use the same parametrization as above. For $\ell \neq -1$ we argue as before, proving that
$\varpi_\ell^-(\bmmu_0(\tau))$ is bounded away from zero if $\tau_0$ is small enough.

If $\ell =-1$, the zero and first Taylor coefficients vanish, however 
$$
\varpi_{-1}^-(\bmmu_0(\tau)) = 
 -N \frac{\tm \tn^2}{|\ubmk|^5} (\tm^2 - 2\tn^2) \tau^2 + \mathcal{O}(\tau^3).
$$
which does not vanish. \\

\paragraph{\bf Analysis near $\infty$}
With the parametrization 
$$
(0, \tau_0) \mapsto \bmmu_{-\infty}(\tau) = \frac{1}{\tau}
    \begin{pmatrix}
\frac{-\tm \tau}{\sqrt{3\tm^2 + 4\tn^2}} \left(\frac{\tn}{2} + \frac{1}{\tau} \right) - \frac{\tm \tau }{2}\\
1
    \end{pmatrix}
$$
equivalent to $(\varphi_-(y), y)$ as $y \to - \infty$.
Evaluating $\varpi_\ell^-(\bmmu)$ along this parametrization yields
\begin{equation}
\varpi_\ell^-\big(  \bmmu_{-\infty}(\tau)  \big) =  
  (\ell-1) \Omega(\ubmk) + \mathcal{O}(\tau)
 , \ \ \forall |\tau| \leq \tau_0
\end{equation}
which  is bounded away from zero since $\ell \neq 1$. \\

\paragraph{\bf Analysis in $(\tau_0, {\tau_0}^{-1})$}
Also
$ \varpi_\ell^-(\bmmu)$ is real analytic in $\R^2 \setminus \{  - \ell \ubmk, - \ubmk \}$. 
However   $ - \ell \ubmk \in \cR_{\ubmk}$ implies that $\ell = 0$ or $\ell = 1$, see \eqref{linR}. 
The case $\ell = 1$ is excluded in this step,  whereas  $\ell = 0$ gives a singularity at $0$, but the 
 curve $\{ (\varphi_+(y), y) \colon y\in (\tau_0, \tau_0^{-1}) \}$ does not pass through this point. Hence  
$
(\tau_0, {\tau_0}^{-1}) \ni 
 y \mapsto \varpi_\ell^+\big(  (\varphi_+(y), y)  \big) 
$
is real analytic and  
therefore can vanish only in finitely many points. \\

\paragraph{\bf Analysis in $(- {\tau_0}^{-1}, -\tau_0)$}
We need to distinguish the cases $2\frac{\tm}{|\ubmk|} < 1$ and 
$2\frac{\tm}{|\ubmk|} > 1$.\\
$\bullet$ If $2\frac{\tm}{|\ubmk|} > 1$,  the parametrization $(\varphi_-(y), y)$  is real analytic and does not pass through any point of the form  $-\ell\ubmk$.
So again $ y \mapsto \varpi_\ell^-\big(  (\varphi_-(y), y)  \big) $ is real analytic and can vanish only in finitely many points.

$\bullet$ If $2\frac{\tm}{|\ubmk|} < 1$,  the parametrization $(\varphi_-(y), y)$  pass through $-\ubmk$ and it is only Lipschitz near this point. 
We use the parametrization \eqref{para-n} and get 
\begin{equation}
  \varpi_\ell^-\big(  \bmmu_{-\tn}(\tau)  \big) =
  \begin{cases}
  (\ell+1 + \sgn(1-\ell)) \Omega(\ubmk) + \mathcal{O}(\tau)
 , &  \mbox{ if } \ell >0 \  \\
2\Omega(\ubmk) + \cO(\tau) &  \mbox{ if } \ell = 0\\
 (\ell+2)  \Omega(\ubmk)  + \mathcal{O}(\tau) &  \mbox{ if } \ell <0 \  \\
 \dfrac{ \tn (-2 \tn + \sqrt{ \tn^2-3 \tm^2})}{3 \sqrt{-3 \tm^2 + 
  \tn^2} |\ubmk|^2} \Omega(\ubmk) \tau + O(\tau^2) &  \mbox{ if } \ell = -2 
 \end{cases}
\end{equation}
 which is bounded away from zero for $\tau$ small, $\tau \neq 0$.

Since away from $y = - \tn$ the function $y \mapsto \varphi_-(y)$ is real analytic, again it can vanish only in finitely many points.\\

Hence the set
$$
\cD_{\ubmk}^- := \bigcup_{|\ell| < \ell_0 , \atop \ell \neq 1} \left\{(\varpi_\ell^-)^{-1}(0) \cap \cR_{\ubmk} \right\}
$$
is discrete. 
For $\bmmu \in \cR_{\ubmk}\setminus  \cD_{\ubmk}^-$ one argues as in \eqref{step1.concl} to get 
 the  claim \eqref{claimStep2.app}. \\

Setting $\cD_{\ubmk} := \cD_{\ubmk}^+ \cup \cD_{\ubmk}^- \cup \{ \bmzero \}$, we get the thesis.

\end{proof}

\section{ Alternative Fourier-Bloch-Floquet approach} \label{app:FBF}

In view of Proposition  \ref{prop:trav}, the operator 
\begin{equation*} 
\mathcal{L}_{\epsilon} \coloneqq   
(\underline{\bmc} \cdot \nabla) \mathbb{I}_2 +
\begin{pmatrix}
- \gradp \psi_{\epsilon} \cdot \nabla & \left( -N^2 \partial_x - \nabla b_{\epsilon} \cdot \gradp \right) (-\Delta)^{-1} \\
-\partial_x  & - \gradp \psi_{\epsilon} \cdot \nabla  - \nabla \omega_{\epsilon} \cdot \gradp  (-\Delta)^{-1} 
\end{pmatrix}
\end{equation*}
depends on $\bmx$ only via the function $\ubmk \cdot \bmx$. 
It is therefore natural to introduce the change of variables
\begin{equation*}
\bmx = z_1 \frac{\ubmk}{|\ubmk|^2}- z_2 \frac{\ubmk^\perp}{|\ubmk|^2}  =: K \bmz,
\end{equation*}
and to define for every function $f \in L^2(\R^2)$ the map 
$$
\cQ[f](\bmz) :=  f ( K \bmz)  \ .
$$
Under this change of variables, one has 
$$
\cQ^{-1}\grad_{\bmx} \cQ:= \vect{\pa_{x_1}}{\pa_{x_2}} = |\ubmk|^2 \, K \, \grad_{\bmz} \ , \quad i.e. \quad \pa_{x_1} = \tm  \pa_{z_1} - \tn \pa_{z_2} \ , 
\quad \pa_{x_2} = \tn  \pa_{z_1} + \tm \pa_{z_2}
$$
and 
$$
\cQ^{-1}\Delta_{\bmx}\cQ   = |\ubmk|^2 \Delta_{\bmz}
$$
so that, using also $\underline{\bmc}= \frac{N\tm}{|\ubmk|^3}\ubmk$, one has
$$
\wt \cL_\e := \cQ^{-1} \cL_\e \cQ = \frac{N\tm}{|\ubmk|}\pa_{z_1} 
\mathbb{I}_2 +
\begin{pmatrix}
\epsilon \frac{|\ubmk| }{N} \cos(z_1) \pa_{z_2} & \left(-  \frac{N^2}{|\ubmk|^2}( \tm  \pa_{z_1} -\tn \pa_{z_2}) - \epsilon \cos(z_1)  \pa_{z_2}\right) (- \Delta_{\bmz})^{-1} \\
-\tm  \pa_{z_1} + \tn \pa_{z_2}  &
\epsilon \frac{ |\ubmk| }{N} \cos(z_1) \pa_{z_2}
- \epsilon \frac{|\ubmk|}{N} \cos(z_1) \pa_{z_2}   (-\Delta_{\bmz})^{-1}
\end{pmatrix}.
$$
In this way  $\wt \cL_\e$ depends only on the variable $z_1$, and following \cite{jiao2025small} we  perform a Fourier transform in the variable $z_2$ and a Bloch-transform in the variable $z_1$.
Precisely we define for  any $g\in \cS(\R^2)$  its Fourier-Bloch-Floquet transform
$$
\tilde{B}(g)(\al,\mu;x) := \sum_{k\in \Z} \hat g (k+\mu,\al) e^{\im k x}\, ,
$$
where we denote by $\hat g (\xi,\al) = \frac1{(2\pi)^2}\int_{\R^2} e^{-\im (\xi x + \al y)} g(x,y) \mathrm{d} x \mathrm{d} y$ the Fourier transform of $g$, that conjugates
the operator  $\cL_\e$ into the fiber operator 
\begin{equation}
\label{cLame1}
\cL(\al,\mu,\e) :=
    \frac{N\tm}{|\ubmk|}(\pa_{z_1} + \im \mu) 
\mathbb{I}_2 +
\begin{pmatrix}
\epsilon \frac{|\ubmk| }{N} \cos(z_1)\im \alpha  & \left(- \frac{N^2}{|\ubmk|^2}( \tm  (\pa_{z_1} + \im \mu)  - \im \tn \alpha ) - \epsilon \cos(z_1)  \im \alpha \right) (|D + \mu|^2 + \alpha^2)^{-1} \\
-\tm  \pa_{z_1} + \tn \im \alpha   &
\epsilon \frac{1 }{|\ubmk| N} \cos(z_1) \pa_{z_2}
- \epsilon \frac{|\ubmk|}{N} \cos(z_1) \pa_{z_2}   (|D + \mu|^2 + \alpha^2)^{-1} 
\end{pmatrix}
 \ ,
\end{equation}
that we consider, for $(\al, \mu) \in \R^2 \setminus (\{0\} \times \Z)$, acting on $H^1(\T) \times L^2(\T)$ with domain 
 \begin{align*}
     D( \cL(\alpha,\mu,\e) ) &\coloneqq \{ U \in H^1(\T) \times L^2(\T) \colon \cL(\alpha,\mu,\e) U \in H^1(\T) \times L^2(\T) \} \ . 
\end{align*}
Following \cite[pp. 13-14]{jiao2025small}, one can show that
\begin{equation*}
\sigma_{X}(\cL_\epsilon)  \supseteq \bigcup_{(\alpha, \mu) \in \R^2 \setminus (\{ 0 \} \times \Z)} 
\sigma^{p}_{H^1(\T) \times L^2(\T)}(\cL(\alpha, \mu, \epsilon)) 
\end{equation*}
with  $\sigma^p_{H^1(\T) \times L^2(\T)}(\cL(\alpha, \mu, \epsilon))$ denoting the point spectrum of $\cL(\alpha, \mu, \epsilon)$. \\
The operator
\begin{equation}
\cL(\al,\mu,0) :=
    \frac{N\tm}{|\ubmk|}(\pa_{z_1} + \im \mu) 
\mathbb{I}_2 +
\begin{pmatrix}
0  & - \frac{N^2}{|\ubmk|^2}( \tm  (\pa_{z_1} + \im \mu)  - \im \tn \alpha )  (|D + \mu|^2 + \alpha^2)^{-1} \\
-\tm  \pa_{z_1} + \tn \im \alpha &
 0 
\end{pmatrix}
 \ 
\end{equation}
has eigenvalues given by 
$$
\lambda_j^\sigma(\al, \mu) := \im \omega_j^\sigma(\al, \mu) \ , 
\quad 
\omega_j^\sigma(\al, \mu):= \frac{N\tm}{|\ubmk|}(j +  \mu)  
+\sigma  \frac{N}{|\ubmk|} \frac{ \tm (j +\mu) - \tn \alpha }{(|j + \mu|^2 + \alpha^2)^{\frac12}}
$$
which, using $\frac{N\tm}{|\ubmk|} = \underline{\bmc} \cdot \ubmk$ and writing $\bmmu = \mu \ubmk - \alpha \ubmk^\perp$, write as 
$$
\omega_j^\sigma(\al, \mu):= \underline{\bmc} \cdot (j \ubmk +  \bmmu)  
+\sigma  \Omega(j \ubmk + \bmmu)
$$
as in the old formula \eqref{eq:lambdasigmak-restr}. \\

\footnotesize{
\noindent {\bf Acknowledgments.}
R.B. is deeply grateful to Thierry Dauxois for introducing her to the problem and for many stimulating discussions over the years. She also warmly thanks Sylvain Joubaud and Antoine Venaille for their careful reading of a preliminary version of this paper, their valuable comments, and for bringing the work~\cite{akylas2023stability} to her attention.

The authors would like thank the anonymous referees for their valuable suggestions.

Work supported by the PRIN project 2022HSSYPN “Turbulent Effects vs Stability in Equations from Oceanography TESEO", PNRR Italia Domani, funded by the European Union via the program NextGenerationEU, CUP G53D23001790001 and B53D23009300001. A.M. is also supported by the European Union  ERC CONSOLIDATOR GRANT 2023 GUnDHam, Project Number: 101124921.
We also thank GNAMPA support. RB acknowledges support of the Institut Henri Poincaré (UAR 839 CNRS-Sorbonne Université), and LabEx CARMIN (ANR-10-LABX-59-01).
 Views and opinions expressed are however those of the authors only and do not necessarily reflect those of the European Union or the European Research Council. Neither the European Union nor the granting authority can be held responsible for them.
}

\bibliographystyle{abbrv}
\bibliography{Boussinesq_Appr}

\end{document}